\documentclass[11pt]{article}

\usepackage[a4paper,margin=1in]{geometry}
\usepackage{amsmath,amssymb,amsthm}
\usepackage[T1]{fontenc}
\IfFileExists{lmodern.sty}{\usepackage{lmodern}}{}
\IfFileExists{microtype.sty}{\usepackage{microtype}}{}
\usepackage{booktabs,array,longtable,graphicx}
\IfFileExists{placeins.sty}{\usepackage{placeins}}{}
\newif\iffigurepdfs
\IfFileExists{figures/exponent-regions.pdf}{%
  \IfFileExists{figures/mass-frequency-intercepts.pdf}{%
    \IfFileExists{figures/zero-slope-regions.pdf}{\figurepdfstrue}{}%
  }{}%
}{}
\iffigurepdfs\else
  \usepackage{tikz}
  \definecolor{mapup}{HTML}{28694C}
  \definecolor{mapdown}{HTML}{32678B}
  \definecolor{mapfold}{HTML}{A6531E}
\fi
\usepackage[colorlinks=true,linkcolor=blue,citecolor=blue,urlcolor=blue]{hyperref}
\hypersetup{pdftitle={Uniqueness of positive normalized solutions for the double-power nonlinear Schr\"odinger
equation in R\string^3},
pdfkeywords={prescribed mass, double-power nonlinear Schrodinger equation, normalized solutions, uniqueness, mass-frequency curve, zero-frequency bifurcation}}
\newcolumntype{L}[1]{>{\raggedright\arraybackslash}p{#1}}
\newtheorem{theorem}{Theorem}[section]
\newtheorem{lemma}[theorem]{Lemma}
\newtheorem{proposition}[theorem]{Proposition}
\numberwithin{equation}{section}
\newcommand{\R}{\mathbb R}
\newcommand{\dd}{\,\mathrm d}

\title{Uniqueness of positive normalized solutions for the
double-power nonlinear Schr\"odinger equation in $\mathbb{R}^3$}
\author{Huxiao Luo$\thanks{{\small Corresponding author. E-mail: luohuxiao@zjnu.edu.cn (H. Luo).}}$\\
\small\emph{School of Mathematical Sciences, Zhejiang Normal University, Jinhua, Zhejiang, {\rm321004}, China}
}
\date{}
\begin{document}
\maketitle

\begin{abstract}
We study uniqueness up to translations for positive solutions of
\[
 -\Delta u+\lambda u=u^q-u^p \quad\text{in }\mathbb R^3
\]
with prescribed $L^2$ mass. Let $M(\lambda)$ denote the mass at
frequency $\lambda$. For $1<q<p$, every mass admits at most one
solution if and only if $q\le7/3$. In this range, $M'(\lambda)>0$
at every admissible frequency, proving the subcritical and critical
cases of the conjecture of Lewin and Rota Nodari.
For $1<p<q<5$ with $q\ge7/3$, we prove that $M'(\lambda)<0$
for every $\lambda>0$.

The mass has exactly one critical point, a nondegenerate minimum,
in two regimes: $q$ just above $7/3$ with $p>q$ and $p-q$ bounded,
and $q$ just below $7/3$ with $p<q$ and $p$ bounded away from one.
The conclusion is uniform in each regime.

Let $p_\dagger$ be the unique root in $(1,9/5)$ of
\[
 15p^3-38p^2-p+36=0.
\]
Near $(p_\dagger,2p_\dagger-1)$, the equation $M'(0+)=0$ defines
a $C^1$ curve. If $M'(0+)\ge0$, then $M'(\lambda)>0$ for all
$\lambda>0$. If $M'(0+)<0$, the mass has exactly one critical point
on $(0,\infty)$, a nondegenerate global minimum.
\end{abstract}
\begin{quotation}\small\noindent
\end{quotation}

\section{Introduction}
Given $p,q>1$ and $a>0$, we count pairs $(\lambda,u)$ satisfying
\begin{equation}\label{eq:main}
 -\Delta u+\lambda u=u^q-u^p\quad\text{in }\R^3,\qquad
 u>0,\qquad \int_{\R^3}u^2\dd x=a .
\end{equation}
The corresponding time-dependent nonlinear Schr\"odinger equation is
\[
 i\partial_t\psi
 =-\Delta\psi-|\psi|^{q-1}\psi+|\psi|^{p-1}\psi .
\]
For $(p,q)=(5,3)$, the cubic term models self-focusing of a laser beam
and the quintic term opposes focusing at high intensity \cite{KOPV}.
Substituting $\psi(t,x)=e^{i\lambda t}u(x)$ gives the elliptic equation
in \eqref{eq:main}. The conserved mass of this standing wave is $\int_{\R^3}u^2\dd x=a$.

We count translates only once. By a positive finite-energy weak solution we mean a function
\[
 u\in H^1(\mathbb R^3)\cap L^{1+\max\{p,q\}}(\mathbb R^3),
 \qquad u>0 \quad\text{a.e.},
\]
such that
\[
 \int_{\mathbb R^3}\nabla u\cdot\nabla\varphi\,dx
 +\lambda\int_{\mathbb R^3}u\varphi\,dx
 =\int_{\mathbb R^3}(u^q-u^p)\varphi\,dx
 \qquad\text{for every }\varphi\in C_c^\infty(\mathbb R^3).
\]
We count classical solutions in this energy space. Lemma~\ref{lem:regularity}
shows that every positive finite-energy weak solution is classical for
$p>q$ or $1<p<q\le5$. For $1<p<q$ and $q>5$,
Proposition~\ref{prop:frequency} rules out both classical finite-energy
solutions and radial finite-energy weak solutions.

\subsection{The mass at fixed frequency}
First fix $\lambda>0$ and consider \eqref{eq:main} without the mass
constraint. By the symmetry theorem of Gidas, Ni and Nirenberg
\cite[Theorem~2]{GNN}, every positive solution that tends to zero
at infinity is radial about some point and strictly decreasing in the
distance from that point. Translate the center to the origin and write
$r=|x|$. The radial equation is
\begin{equation}\label{eq:radial}
 -u''(r)-\frac{2}{r}u'(r)+\lambda u(r)=u(r)^q-u(r)^p,\qquad r>0.
\end{equation}
A smooth radial profile satisfies $u'(0)=0$.

If $1<q<p$, a positive solution in the prescribed energy space exists precisely for
$0<\lambda<\lambda_*$, and the radial solution centered at the origin
is unique \cite[Section~3 and Theorem~2]{LRN}. The upper endpoint is
\[
 \lambda_*=\max_{s>0}
 \left(\frac{2s^{q-1}}{q+1}-\frac{2s^{p-1}}{p+1}\right).
\]

If $1<p<q<5$, a positive radial solution exists for every $\lambda>0$
and is unique once its center is fixed \cite[Proposition~1.1 and its proof]{FH}.
At $\lambda=0$ the same result gives a unique centered positive radial
profile $u_0\in\dot H^1(\R^3)\cap L^{p+1}(\R^3)$ that tends to zero at
infinity. By \cite[Theorem~1]{LiNi}, every positive decaying solution at
zero frequency is radial about some point. The required sign near zero follows from
\[
 \frac{d}{ds}(s^q-s^p)=s^{p-1}(q s^{q-p}-p)<0
 \qquad\text{for all sufficiently small }s>0.
\]
The profile $u_0$ belongs to $L^2(\R^3)$ exactly when $p<7/3$
\cite[Corollary~1.4]{FH}.

For each admissible $\lambda>0$, let $u_\lambda$ denote the positive
radial solution centered at the origin, and define
\[
 M_{p,q}(\lambda)=\int_{\mathbb R^3}u_\lambda^2\,dx.
\]
For fixed $p,q$ we write $M(\lambda)$. Each positive root of
$M(\lambda)=a$ gives one solution up to translations. If $1<p<q<5$
and $p<7/3$, we also count $u_0$ at $a=a_0=\int_{\R^3}u_0^2\dd x$.

Jeanjean, Zhang and Zhong \cite{JZZ} studied prescribed mass for
$-\Delta u+\lambda u=g(u)$. Under their hypotheses, a connected set of
positive solutions covers all $\lambda>0$. They used the mass limits as
$\lambda\downarrow0$ and $\lambda\to\infty$ to prove existence and,
in some cases, multiplicity. They assume $g(s)>0$ for $s>0$.
For $p\ne q$, our nonlinearity $s^q-s^p$ changes sign at $s=1$. To count the solutions at each
mass, we determine where $M$ increases and where it decreases.

The exponent $7/3$ comes from the single focusing power. For $1<q<5$, let
$Q_q$ be the unique centered positive solution of
$-\Delta Q_q+Q_q=Q_q^q$ \cite{Kwong}. Scaling gives
\[
 \int_{\R^3}\big|\lambda^{1/(q-1)}Q_q(\sqrt{\lambda}\,x)\big|^2\dd x
 =\lambda^{2/(q-1)-3/2}\int_{\R^3}Q_q^2\dd x .
\]
The exponent of $\lambda$ vanishes at $q=7/3$. We write $Q=Q_{7/3}$
and $m_Q=\int Q^2$. Unless specified, integrals are over $\R^3$.

\subsection{Exponent ranges and solution counts}
For $1<q<p$, Lewin and Rota Nodari determined the endpoint behavior
of $M$ and conjectured its shape between the endpoints
\cite[Theorems~3--4 and Conjecture~2]{LRN}. We prove that $M'(\lambda)>0$
on $(0,\lambda_*)$ when $q\le7/3$. This proves the subcritical and
critical cases of the conjecture in dimension three. When $q>7/3$,
the mass diverges at both endpoints, so every mass above its minimum
occurs at least twice. Thus $q\le7/3$ is precisely the condition for
uniqueness at every mass for which a solution exists
(Theorem~\ref{thm:increasing}).

For $1<p<q<5$, the order of the powers reverses. We prove that
$M'(\lambda)<0$ for every $\lambda>0$ when $q\ge7/3$
(Theorem~\ref{thm:powers}). In dimension one, Kfoury, Le Coz and Tsai
\cite[Lemma~26]{KLT} found that for $1<p<q<5$ the mass increases if
$2p+q\le7$ and has a single minimum if $2p+q>7$.

In dimension three, a unique nondegenerate minimum occurs in two
ranges near $7/3$. For $p>q$ we take $q$ just above $7/3$, with
$p-q$ bounded above (Theorem~\ref{thm:uniform-fold}). For $p<q$ we
take $q$ just below $7/3$, with $p$ bounded away from one
(Theorem~\ref{thm:reverse-uniform}). In both ranges $M$ decreases
up to the minimum and increases afterwards. For fixed $1<p<7/3$,
Proposition~\ref{prop:minimum-at-infinity} gives the leading terms of
the minimizing frequency and mass as $q\uparrow7/3$.

For $v(r)=(1+r^2)^{-1/(p-1)}$, $\Delta v$ is a linear combination
of $v^p$ and $v^{2p-1}$. For $q=2p-1$ and $1<p<9/5$, a rescaling
of $v$ gives the zero-frequency profile $u_0$. Solving the linearized equation at $u_0$
gives the exponent $p_\dagger$ where $M'(0+)$ changes sign.
Near $(p_\dagger,2p_\dagger-1)$, a nonnegative initial slope gives
$M'(\lambda)>0$ for all $\lambda>0$. A negative initial slope gives
a unique nondegenerate minimum at positive frequency
(Theorem~\ref{thm:global-birth}).

\begin{figure}[!ht]\centering
\iffigurepdfs
  \includegraphics{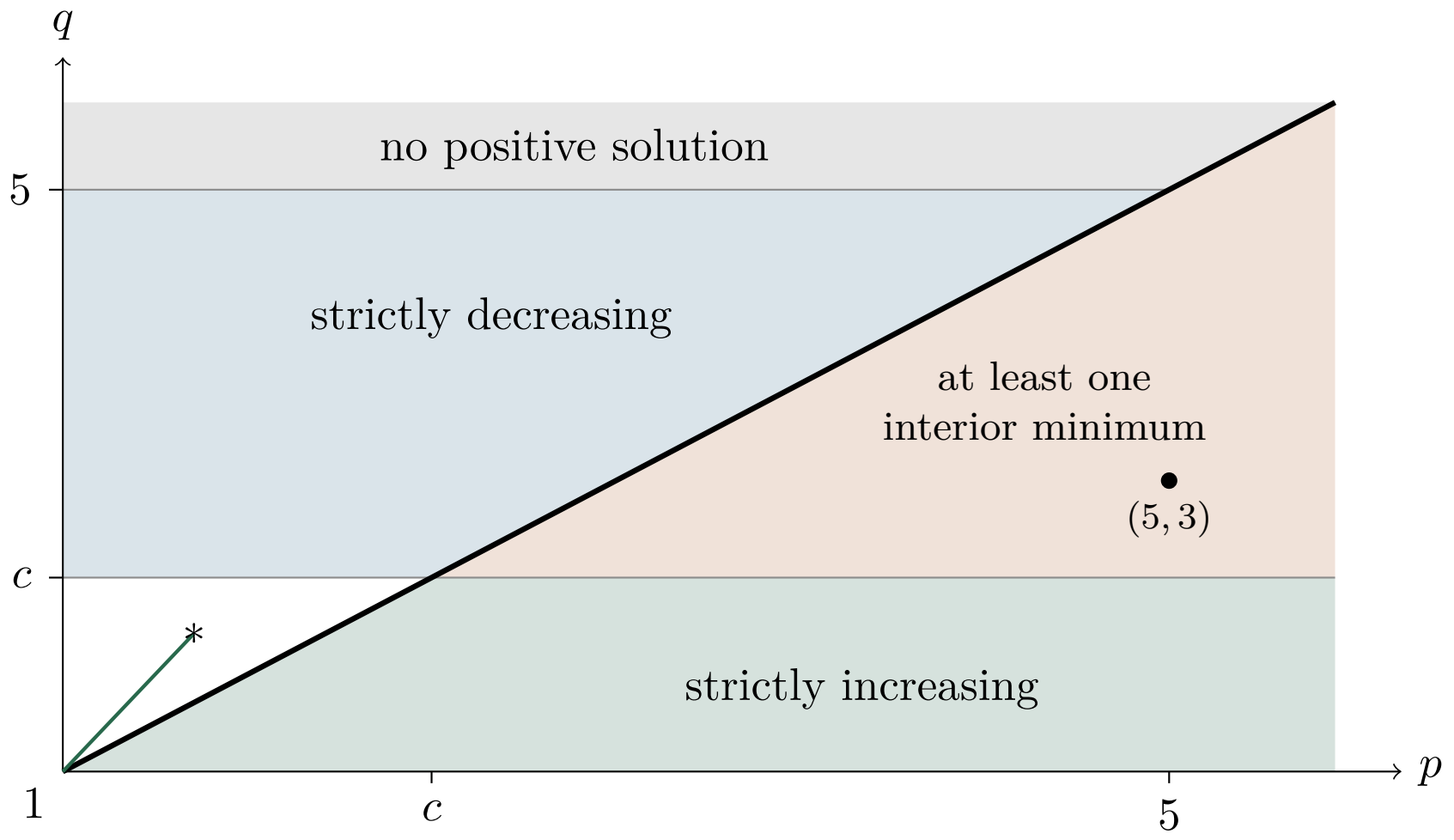}%
\else
\begin{tikzpicture}[x=2.13cm,y=1.12cm,font=\small]
\def\crit{2.3333333}\def\edge{5.60}
\fill[mapup!19] (1,1)--(\crit,\crit)--(\edge,\crit)--(\edge,1)--cycle;
\fill[mapfold!17] (\crit,\crit)--(\edge,\edge)--(\edge,\crit)--cycle;
\fill[mapdown!18] (1,\crit)--(\crit,\crit)--(5,5)--(1,5)--cycle;
\fill[white] (1,1)--(1,\crit)--(\crit,\crit)--cycle;
\fill[black!10] (1,5)--(5,5)--(\edge,\edge)--(1,\edge)--cycle;
\draw[thin,black!45] (1,\crit)--(\edge,\crit);
\draw[thin,black!45] (1,5)--(5,5);
\draw[very thick] (1,1)--(\edge,\edge);
\draw[->] (1,1)--(5.84,1) node[right] {$p$};
\draw[->] (1,1)--(1,5.91) node[above] {$q$};
\foreach \x/\s in {2.3333333/c,5/5}{\draw (\x,1)--++(0,-.08) node[below] {$\s$};}
\foreach \y/\s in {2.3333333/c,5/5}{\draw (1,\y)--++(-.05,0) node[left] {$\s$};}
\node[below left] at (1,1) {$1$};
\node at (3.89,1.57) {strictly increasing};
\node at (2.55,4.12) {strictly decreasing};
\draw[thick,mapup] (1,1)--(1.47414,1.94829);
\node[inner sep=0pt] at (1.47414,1.94829) {$\ast$};
\node[align=center,font=\footnotesize] at (4.55,3.55) {at least one\\interior minimum};
\node at (2.85,5.28) {no positive solution};
\fill[black] (5,3) circle[radius=1.8pt];
\node[below,font=\scriptsize] at (5,2.98) {$(5,3)$};
\end{tikzpicture}
\fi
\caption{The exponent plane, with $c=7/3$. The colors describe the mass
at positive frequencies. On $q=c$, the mass increases for $p>c$ and
decreases for $p<c$. For $p>q>c$, every mass above the minimum occurs
at least twice. The dot marks the cubic--quintic case $(5,3)$
\cite{KOPV,ZZ,ZWZ}. On the green segment $q=2p-1$, $1<p\le p_\dagger$,
the mass increases from $M(0)$ to infinity. Panel F of
Figure~\ref{fig:mass-frequency-global} shows the curves near the star
$(p_\dagger,2p_\dagger-1)$. A zero-frequency solution has finite mass
exactly for $1<p<q<5$ and $p<c$. Its mass is $M_{p,q}(0)$.}\label{fig:global-map}
\end{figure}

\begin{figure}[!htbp]
\centering
\resizebox{.96\linewidth}{!}{%
\iffigurepdfs
  \includegraphics{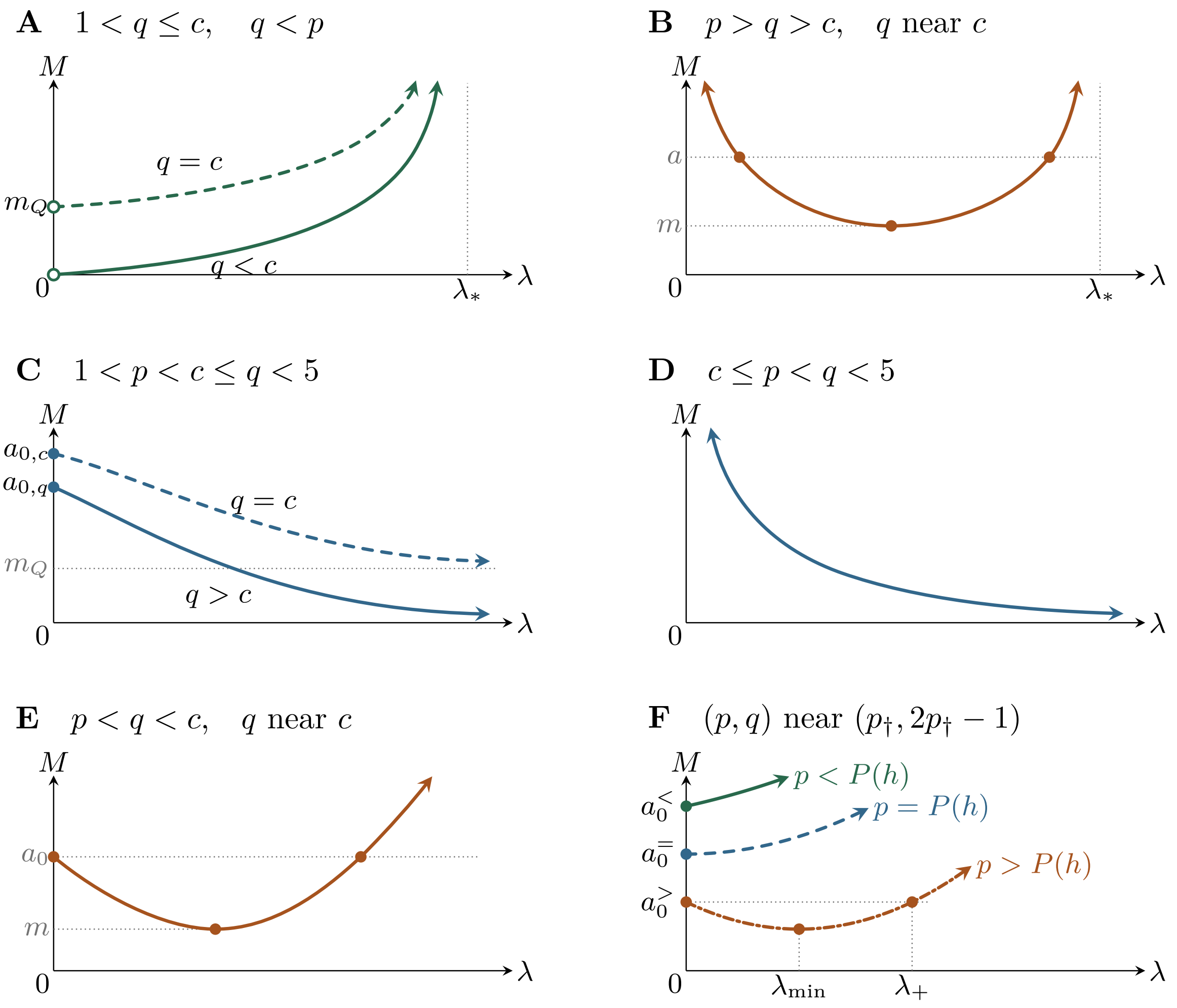}%
\else
\begin{tikzpicture}[x=1cm,y=1cm,font=\small,>=stealth,
  curve/.style={line width=1.05pt,line cap=round},
  level/.style={black!55,thin,densely dotted},
  labeltext/.style={font=\footnotesize,inner sep=1pt}]
\newcommand{\massaxes}[1]{%
  \draw[->,thin] (0,0)--(5.08,0) node[right,labeltext] {$\lambda$};
  \draw[->,thin] (0,0)--(0,2.16) node[above,labeltext] {$M$};
  \node[below left,labeltext] at (0,0) {$0$};
  \node[anchor=west] at (-.56,2.77) {#1};}
\newcommand{\massfinite}{%
  \draw[level] (4.58,0)--(4.58,2.12);
  \node[below,labeltext] at (4.58,0) {$\lambda_*$};}
\newcommand{\massopen}[2]{\draw[draw=#1,fill=white,line width=.85pt]
  #2 circle[radius=1.7pt];}
\newcommand{\massdot}[2]{\fill[#1] #2 circle[radius=1.8pt];}

\begin{scope}
 \massaxes{\textbf{A}\quad $1<q\le c,\quad q<p$}
 \massfinite
 \draw[curve,mapup,->] (0,0) .. controls (1.85,.13) and (3.45,.46) ..
     (3.99,1.36) .. controls (4.15,1.63) and (4.22,1.91) .. (4.25,2.15);
 \draw[curve,mapup,dashed,->] (0,.75) .. controls (1.75,.84) and (3.12,1.18) ..
     (3.63,1.65) .. controls (3.85,1.83) and (3.96,2.00) .. (4.01,2.15);
 \node[left,labeltext] at (0,.75) {$m_Q$};
 \node[labeltext] at (2.10,.07) {$q<c$};
 \node[labeltext] at (1.50,1.20) {$q=c$};
 \massopen{mapup}{(0,0)} \massopen{mapup}{(0,.75)}
\end{scope}
\begin{scope}[shift={(7.00,0)}]
 \massaxes{\textbf{B}\quad $p>q>c$,\quad $q$ near $c$}
 \massfinite
 \draw[level] (0,1.30) node[left,labeltext] {$a$} -- (4.58,1.30);
 \draw[level] (0,.54) node[left,labeltext] {$m$} -- (2.27,.54);
 \draw[curve,mapfold,<->] (.20,2.15) .. controls (.30,1.79) and (.43,1.48) ..
    (.59,1.30) .. controls (1.02,.78) and (1.67,.54) .. (2.27,.54)
    .. controls (2.87,.54) and (3.60,.78) .. (4.02,1.30)
    .. controls (4.18,1.50) and (4.29,1.85) .. (4.34,2.15);
 \massdot{mapfold}{(.59,1.30)} \massdot{mapfold}{(4.02,1.30)}
 \massdot{mapfold}{(2.27,.54)}
\end{scope}

\begin{scope}[shift={(0,-3.85)}]
 \massaxes{\textbf{C}\quad $1<p<c\le q<5$}
 \draw[level] (0,.60) node[left,labeltext] {$m_Q$} -- (4.92,.60);
 \draw[curve,mapdown,->] (0,1.50) .. controls (.94,1.09) and (2.30,.16) ..
     (4.83,.095);
 \draw[curve,mapdown,dashed,->] (0,1.87) .. controls (.96,1.64) and
     (2.60,.735) .. (4.83,.68);
 \node[left,labeltext] at (0,1.87) {$a_{0,c}$};
 \node[left,labeltext] at (0,1.50) {$a_{0,q}$};
 \node[labeltext] at (2.32,1.30) {$q=c$};
 \node[labeltext] at (1.82,.28) {$q>c$};
 \massdot{mapdown}{(0,1.87)} \massdot{mapdown}{(0,1.50)}
\end{scope}
\begin{scope}[shift={(7.00,-3.85)}]
 \massaxes{\textbf{D}\quad $c\le p<q<5$}
 \draw[curve,mapdown,<->] (.27,2.16) .. controls (.45,1.38) and
    (.93,.83) .. (1.72,.55) .. controls (2.55,.26) and (3.61,.14) ..
    (4.84,.10);
\end{scope}

\begin{scope}[shift={(0,-7.70)}]
 \massaxes{\textbf{E}\quad $p<q<c$,\quad $q$ near $c$}
 \draw[level] (0,1.26) node[left,labeltext] {$a_0$} -- (4.70,1.26);
 \draw[level] (0,.46) node[left,labeltext] {$m$} -- (1.79,.46);
 \draw[curve,mapfold,->] (0,1.26) .. controls (.48,.86) and (1.20,.46) ..
    (1.79,.46) .. controls (2.42,.46) and (2.96,.85) .. (3.40,1.26)
    .. controls (3.72,1.58) and (3.98,1.88) .. (4.19,2.15);
 \massdot{mapfold}{(0,1.26)} \massdot{mapfold}{(1.79,.46)}
 \massdot{mapfold}{(3.40,1.26)}
\end{scope}

\begin{scope}[shift={(7.00,-7.70)}]
 \massaxes{\textbf{F}\quad $(p,q)$ near $(p_\dagger,2p_\dagger-1)$}
 \node[left=2pt,labeltext] at (0,1.82) {$a_0^{<}$};
 \node[left=2pt,labeltext] at (0,1.29) {$a_0^{=}$};
 \node[left=2pt,labeltext] at (0,.76) {$a_0^{>}$};
 \draw[level] (0,.76) -- (2.68,.76);
 \draw[level] (1.25,0) -- (1.25,.46);
 \draw[level] (2.50,0) -- (2.50,.76);
 \draw[curve,mapup,->,domain=0:1.14,samples=30]
    plot (\x,{1.82+.22*\x+.06*\x*\x}) node[right,labeltext] {$p<P(h)$};
 \draw[curve,mapdown,dashed,->,domain=0:2.02,samples=30]
    plot (\x,{1.29+.125*\x*\x}) node[right,labeltext] {$p=P(h)$};
 \draw[curve,mapfold,dash pattern=on 3pt off 1.3pt on .8pt off 1.3pt,
       ->,domain=0:3.16,samples=45]
    plot (\x,{0.76-.48*\x+.192*\x*\x}) node[right,labeltext] {$p>P(h)$};
 \massdot{mapup}{(0,1.82)} \massdot{mapdown}{(0,1.29)}
 \massdot{mapfold}{(0,.76)} \massdot{mapfold}{(1.25,.46)}
 \massdot{mapfold}{(2.50,.76)}
 \node[below,labeltext] at (1.25,0) {$\lambda_{\min}$};
 \node[below,labeltext] at (2.50,0) {$\lambda_+$};
\end{scope}
\end{tikzpicture}%
\fi}
\caption{Mass curves for $c=7/3$. Each frequency at mass $a$ gives one
solution up to translations. Open circles mark excluded endpoints.
In A, $\lambda_*$ depends on $(p,q)$. In C, $a_{0,c}$ and $a_{0,q}$ are the
intercepts for $q=c$ and $q>c$. B assumes $0<p-q\le D$ and $0<q-c<e_D$
(Theorem~\ref{thm:uniform-fold}). E assumes $b\le p<q<c$ and
$c-q<\varepsilon_b$, with $1<b<c$ fixed (Theorem~\ref{thm:reverse-uniform}).
In F, $h=q-2p+1$ is fixed. The three intercepts are the masses $M_{p,q}(0)$
for $p<P(h)$, $p=P(h)$ and $p>P(h)$, in schematic vertical order.
The dashed curve starts horizontally. The dash-dotted curve reaches its
minimum at $\lambda_{\min}$ and returns to $a_0^{>}$ at $\lambda_+$
(Theorem~\ref{thm:global-birth}). All three grow without bound.}
\label{fig:mass-frequency-global}
\end{figure}

Figures~\ref{fig:global-map} and~\ref{fig:mass-frequency-global} show
the exponent ranges and their mass curves. Figure~\ref{fig:mass-frequency-local}
shows the transition near $(p_\dagger,2p_\dagger-1)$.
Table~\ref{tab:complete} compares the three-dimensional results with earlier work.

Two questions remain. For $p>q>7/3$, every mass above the minimum admits
at least two solutions, and every sufficiently large mass admits exactly
two \cite[Corollary~5]{LRN}. Does $M'$ have exactly one zero for every
such pair? For $1<p<q<7/3$, solutions exist exactly for
$a\ge\min_{\lambda\ge0}M(\lambda)$, and sufficiently large masses admit
exactly one (Proposition~\ref{prop:explicit-upper-frequency}). The count
at intermediate masses remains open in the unclassified part of the white
triangle in Figure~\ref{fig:global-map}.

\clearpage

\begingroup
\small
\setlength{\LTcapwidth}{\linewidth}
\begin{longtable}{@{}L{.22\linewidth}L{.31\linewidth}L{\dimexpr.47\linewidth-4\tabcolsep\relax}@{}}
\caption{Results in dimension three, with $c=7/3$. We count translates
once and include the profile at $\lambda=0$ when its mass $a_0$ is finite.
Here $m$ denotes the global minimum mass.}\label{tab:complete}\\
\toprule
Exponents & Earlier result & Conclusion in this paper\\\midrule
\endfirsthead
\multicolumn{3}{@{}l}{\tablename~\thetable\ (continued)}\\
\toprule
Exponents & Earlier result & Conclusion in this paper\\\midrule
\endhead
\midrule
\endfoot
\bottomrule
\endlastfoot
$1<q\le c$, $q<p$
& Fixed-frequency uniqueness and positive endpoint slopes \cite{LRN}.
& $M'>0$ for all admissible frequencies. One solution for all $a>0$ if $q<c$, and exactly
for $a>m_Q$ if $q=c$ (Theorem~\ref{thm:increasing}).\\
$p>q>c$
& The endpoint limits give at least two solutions for $a>m$,
and exactly two are known for all sufficiently large $a$
\cite[Theorems~3--4 and Corollary~5]{LRN}. The full count for $(5,3)$
is in \cite{ZZ,ZWZ}.
& For every $D>0$, $0<p-q\le D$ and $0<q-c<e_D$ give one nondegenerate
minimum and counts $0,1,2$ below, at and above $m$
(Theorem~\ref{thm:uniform-fold}).\\
$1<p<q<5$, $q\ge c$
& Strong instability holds at every positive frequency
\cite[Theorem~1.8]{FH}. For $q>c$, \cite[Theorem~1]{BFG}
determines the masses at which the constrained ground-state energy is attained.
& $M'<0$ for all positive frequencies. The unique-solution range is
$(\ell,a_0]$ if $p<c$, or $(0,\infty)$ if $p\ge c$.
Here $\ell=m_Q$ for $q=c$ and $\ell=0$ for $q>c$
(Theorem~\ref{thm:powers}).\\
$b\le p<q<c$, fixed $b>1$, $q$ near $c$
& The zero-frequency profile has finite mass, with strong $H^1$ convergence
at the endpoint \cite{FH}.
& If $0<c-q<\varepsilon_b$, there is a single nondegenerate minimum.
The counts are $0,1,2,1$ on $(0,m),\{m\},(m,a_0],(a_0,\infty)$
(Theorem~\ref{thm:reverse-uniform}).\\
$q=2p-1$, $1<p<5/3$
& The zero-frequency profile is explicit \cite{DPD}.
& $M'(\lambda)>0$ for every $\lambda>0$ exactly when
$p\le p_\dagger$ (Theorem~\ref{thm:global-birth}).
For $p_\dagger<p<5/3$, the initial slope is negative and the
global minimum is attained at positive frequency.\\
Near $(p_\dagger,2p_\dagger-1)$
& The frequency derivative converges on bounded radial intervals
as $\lambda\downarrow0$ \cite[Theorem~4.8]{FHnew}.
& The $C^1$ curve separates global increase from a single global minimum.
For $M'(0+)<0$, the counts are $0,1,2,1$ on $(0,m),\{m\},(m,a_0],(a_0,\infty)$
(Theorem~\ref{thm:global-birth}).\\
\end{longtable}
\endgroup

\subsection{Frequency restrictions and regularity}\label{sec:nonexist}
We use the energy identities to exclude negative frequencies and two exponent ranges.

\begin{proposition}\label{prop:frequency}
A positive classical finite-energy solution of \eqref{eq:main}, or a
positive radial finite-energy weak solution, must satisfy $\lambda\ge0$.
Neither class contains a solution if $p=q$ or if $1<p<q$ and $q\ge5$.
\end{proposition}

\begin{proof}[\normalfont\proofname]
Take a real $\psi\in C_c^\infty(\R^3)$ and use $\psi^2/u$ as a test
function. For a positive classical solution, $u$ is bounded away from
zero on the support of $\psi$. For a radial weak solution we take
$\psi$ supported in an annulus. There $u$ is continuous by the
one-dimensional Sobolev embedding, and the radial equation gives
$u\in C^2$. If $u(r_0)=0$ for some $r_0>0$, then $u'(r_0)=0$ because
$u\ge0$. Uniqueness for the radial initial value problem would then give
$u\equiv0$. Thus the test function is valid in both cases. The identity
\[
 \nabla u\cdot\nabla(\psi^2/u)
 =|\nabla\psi|^2-\left|\nabla\psi-\frac{\psi}{u}\nabla u\right|^2
\]
and $u^{q-1}\ge0$ give
\[
 -\lambda\int\psi^2\le\int|\nabla\psi|^2+\int u^{p-1}\psi^2 .
\]
Choose $\psi$ supported in $1<|x|<2$ with $\|\psi\|_2=1$ and take
$\psi_R(x)=R^{-3/2}\psi(x/R)$. Since $u\in L^{p+1}$, H\"older's
inequality gives
\begin{align*}
 -\lambda
 &\le R^{-2}\|\nabla\psi\|_2^2+
       \|u\|_{p+1}^{p-1}\|\psi_R\|_{p+1}^2\\
 &=R^{-2}\|\nabla\psi\|_2^2+
       R^{-3(p-1)/(p+1)}\|u\|_{p+1}^{p-1}\|\psi\|_{p+1}^2
 \longrightarrow0 .
\end{align*}
Thus $\lambda\ge0$.

Testing the elliptic equation in \eqref{eq:main} with $u$ and with
$x\cdot\nabla u$, respectively, gives
\begin{align}
 \int|\nabla u|^2+\lambda\int u^2&=\int u^{q+1}-\int u^{p+1},\label{eq:energy-test}\\
 \tfrac12\int|\nabla u|^2+\tfrac32\lambda\int u^2
 &=\tfrac3{q+1}\int u^{q+1}-\tfrac3{p+1}\int u^{p+1}.\label{eq:pohozaev}
\end{align}
The first test follows by approximation in the energy space. For the
second insert a cutoff $\chi(x/R)$ supported in $B_{2R}$ and equal to
one on $B_R$. Since $|x\cdot\nabla[\chi(x/R)]|\le C$, the error satisfies
\[
 |\mathrm{error}|\le C\int_{B_{2R}\setminus B_R}
 (|\nabla u|^2+u^2+u^{p+1}+u^{q+1})\longrightarrow0 .
\]
For a radial finite-energy weak solution use a second cutoff at the
origin. The same integral over $B_{2r}\setminus B_r$ tends to zero
as $r\downarrow0$.

When $1<p<q$ and $q\ge5$, eliminate $\int u^{q+1}$ between the energy identity~\eqref{eq:energy-test}
and the Pohozaev identity~\eqref{eq:pohozaev} to obtain
\[
 \left(\frac3{q+1}-\frac12\right)\int|\nabla u|^2
 =\left(\frac32-\frac3{q+1}\right)\lambda\int u^2
  +3\left(\frac1{p+1}-\frac1{q+1}\right)\int u^{p+1}.
\]
The left side is nonpositive. On the right $\lambda\ge0$ and the last
term is strictly positive. This rules out a classical finite-energy
solution and also a radial finite-energy weak solution. The same
conclusion holds for any weak solution for which the Pohozaev identity~\eqref{eq:pohozaev} is valid.

On the diagonal $p=q$, any positive finite-energy weak solution satisfies $-\Delta u+\lambda u=0$
and is smooth by interior regularity. Since $\lambda\ge0$, testing with $u$ gives
\[
 0=\int|\nabla u|^2+\lambda\int u^2\ge\int|\nabla u|^2 .
\]
Thus $u$ is constant and $u\in L^2(\R^3)$ forces $u=0$.

\end{proof}

The next lemma allows us to apply the symmetry and uniqueness theorems
to weak solutions in the ranges used below.

\begin{lemma}\label{lem:regularity}
If $1<q<p$ or $1<p<q\le5$, every positive finite-energy weak solution
is classical. It tends uniformly to zero at infinity if $1<q<p$ or
$1<p<q<5$.
\end{lemma}
\begin{proof}[\normalfont\proofname]
First assume $p>q$, and let $u$ be a positive finite-energy weak solution.
Choose $B$ large enough that $s^p-s^q+\lambda s>0$ for $s>B$.
Testing with $(u-B)_+$ gives
\[
 0=\int|\nabla(u-B)_+|^2
    +\int_{\{u>B\}}(u^p-u^q+\lambda u)(u-B).
\]
Both integrals are nonnegative, so $u\le B$. Thus
$u^q-u^p-\lambda u\in L^2$, since $u\in L^2\cap L^\infty$.
Elliptic regularity gives $u\in H^2$ and then classical regularity.
Uniform continuity and $u\in L^2$ imply $u(x)\to0$ at infinity.

Now assume $1<p<q\le5$. Since $u\in L^6$ and $q-1\le4$,
\[
 -\Delta u=(u^{q-1}-u^{p-1}-\lambda)u,\qquad
 u^{q-1}-u^{p-1}-\lambda\in L^{3/2}_{\rm loc}.
\]
The local Br\'ezis--Kato estimate gives $u\in L^t_{\rm loc}$ for every
finite $t\ge1$ \cite{BK}. Thus, for every $1<t<\infty$,
\[
 u^q-u^p-\lambda u\in L^t_{\rm loc}
 \quad\Longrightarrow\quad u\in W^{2,t}_{\rm loc}.
\]
Taking $t>3$ gives $u\in C^{1,\alpha}_{\rm loc}$ for some $\alpha>0$.
The right side of the equation is then locally H\"older continuous, so
interior elliptic regularity gives $u\in C^2$.

For $q<5$, Proposition~\ref{prop:frequency} gives $\lambda\ge0$. Hence
\[
 -\Delta u\le u^{q-1}u,\qquad
 u^{q-1}\in L^{6/(q-1)}(\R^3),\qquad \frac{6}{q-1}>\frac32.
\]
The interior $L^2$ to $L^\infty$ estimate for nonnegative subsolutions gives
\[
 \sup_{B_1(x)}u\le C\|u\|_{L^2(B_2(x))}\longrightarrow0
 \qquad(|x|\to\infty).
\]
The constant is independent of $x$ because the local norm of
$u^{q-1}$ is bounded by its norm on $\R^3$.
\end{proof}

\section{\texorpdfstring{The uniqueness threshold for $1<q<p$}{The uniqueness threshold for 1 < q < p}}\label{sec:forward-form}
For $1<q<p$, the endpoint masses are known
\cite[Theorems~3--4]{LRN}. We now show when two positive solutions
can have the same mass.

\begin{theorem}\label{thm:increasing}
For $1<q<p$, every prescribed mass $a>0$ admits at most one positive solution of \eqref{eq:main}, up to
translations, if and only if $q\le7/3$.
\begin{enumerate}
\item If $1<q<7/3$, then $M'(\lambda)>0$ for every $0<\lambda<\lambda_*$. Each $a>0$ admits exactly one positive solution.
\item If $q=7/3$, then $M'(\lambda)>0$ for $0<\lambda<\lambda_*$. There is exactly one positive solution for $a>m_Q$ and none for $0<a\le m_Q$.
\item If $q>7/3$, the minimum $m=\min_{0<\lambda<\lambda_*}M(\lambda)$ is attained. Every $a>m$ admits at least two positive solutions.
\end{enumerate}
\end{theorem}

If $M'(\lambda)=0$, the frequency derivative $w=\partial_\lambda u_\lambda$
satisfies both $\int u_\lambda w=0$ and $\langle Lw,w\rangle=0$.
We exclude this possibility for $q\le7/3$ by proving that the
linearized quadratic form is strictly positive on the radial mass constraint.

\subsection{The branch at fixed frequency}\label{sec:ends}
By Lemma~\ref{lem:regularity}, every positive finite-energy weak solution
is classical and decays at infinity. Proposition~\ref{prop:frequency}
gives $\lambda\ge0$. To exclude $\lambda=0$,
take $R$ so large that $u<1$ for $|x|\ge R$. If $\lambda=0$, then
$-\Delta u=u^q(1-u^{p-q})>0$ there. Harmonic comparison on
$R<|x|<T$, followed by $T\to\infty$, gives
\[
 u(x)\ge\min_{|z|=R}u(z)\,
 \frac{|x|^{-1}-T^{-1}}{R^{-1}-T^{-1}}
 \ \longrightarrow\ \frac{R}{|x|}\min_{|z|=R}u(z).
\]
The square of this lower bound has infinite integral in dimension three.
Therefore $\lambda>0$.

By \cite[Theorems~2--4]{LRN}, the centered profile exists precisely on $(0,\lambda_*)$,
satisfies $0<u_\lambda<1$ and $u_\lambda'(r)<0$ for $r>0$, and decays exponentially together with its first
derivatives. The endpoint values
of the mass in dimension three are
\begin{equation}\label{eq:ends}
 M(\lambda_* -)=+\infty,\qquad
 M(0+)=
 \begin{cases}
 0 & 1<q<7/3,\\
 m_Q & q=7/3,\\
 +\infty & q>7/3 .
 \end{cases}
\end{equation}

We next differentiate the centered profile at a fixed admissible frequency.

\begin{lemma}\label{lem:forward-branch}
Fix $p>1$. For $1<q<p$ and $0<\lambda<\lambda_*$, the centered profile
is locally $C^1$ in $(\lambda,q)$ with values in $H^2_{\rm rad}(\R^3)$.
Its frequency derivative satisfies
\[
 L\partial_\lambda u_\lambda=-u_\lambda,\qquad
 L=-\Delta+\lambda-qu_\lambda^{q-1}+pu_\lambda^{p-1}.
\]
\end{lemma}
\begin{proof}[\normalfont\proofname]
Fix an admissible $(\lambda_0,q_0)$, and denote its centered
profile by $U$. On the real radial spaces define
\[
 F(\lambda,q,v)=-\Delta v+\lambda v-|v|^{q-1}v+|v|^{p-1}v,
 \qquad v\in H^2_{\rm rad}(\R^3),\quad F\in L^2_{\rm rad}(\R^3).
\]
For $s>1$ the derivatives of the signed power are
\[
 D_v(|v|^{s-1}v)h=s|v|^{s-1}h,\qquad
 \partial_s(|v|^{s-1}v)=v\big(|v|^{s-1}\log|v|\big).
\]
We set $|v|^{s-1}\log|v|=0$ at $v=0$. On bounded
$v$-intervals the two multipliers are uniformly continuous for $s$ in
a compact subset of $(1,\infty)$. Since $H^2\hookrightarrow L^\infty\cap L^2$,
multiplying an $L^2$ function by either multiplier gives continuous
derivatives into $L^2$. Thus $F$ is $C^1$.

At $U$ the radial linearization is
\[
 L=D_vF(\lambda_0,q_0,U)=-\Delta+\lambda_0+pU^{p-1}-q_0U^{q_0-1}.
\]
For the full-space operator, \cite[Theorem~2]{LRN} gives
\[
 \ker L=\operatorname{span}\{\partial_{x_1}U,\partial_{x_2}U,\partial_{x_3}U\},
 \qquad \ker L\cap H^2_{\rm rad}=\{0\}.
\]
The second equality follows because every translation derivative is odd
under $x\mapsto-x$ and every radial function is even. Multiplication by
$pU^{p-1}-q_0U^{q_0-1}$ is compact from $H^2$ to $L^2$ by local
compactness and uniform decay of the potential outside large balls.
Since $-\Delta+\lambda_0$ is invertible, the Fredholm alternative gives
a bounded radial inverse for $L$. The implicit function theorem now gives
a local $C^1$ solution $v(\lambda,q)$ through $U$.

For $v_-=\max\{-v,0\}$ we have $\|v_-\|_\infty\le\|v-U\|_\infty\to0$.
Testing $F=0$ with $-v_-$ gives
\begin{equation}\label{eq:profile-negative-part}
 \int(|\nabla v_-|^2+\lambda v_-^2)
 =\int(v_-^{q+1}-v_-^{p+1})
 \le\|v_-\|_\infty^{q-1}\int v_-^2 .
\end{equation}
Near $(\lambda_0,q_0)$ we have $\|v_-\|_\infty^{q-1}<\lambda$,
so $v_-=0$. Also $v\ne0$ by $L^2$ convergence to $U$. The strong
maximum principle gives $v>0$, and fixed-frequency uniqueness gives
$v=u_\lambda$. Differentiation at fixed $p,q$ yields
\begin{equation}\label{eq:linearized-frequency}
 w=\partial_\lambda u_\lambda,\qquad Lw=-u_\lambda,\qquad
 L=-\Delta+\lambda-qu_\lambda^{q-1}+pu_\lambda^{p-1}.
\end{equation}
\end{proof}

\subsection{The radial quadratic form}
Fix an admissible frequency $\lambda$, write $u=u_\lambda$, and use
the linearized operator $L$ in \eqref{eq:linearized-frequency}. A radial perturbation
$h$ is tangent to the mass constraint if $\int_{\R^3}uh=0$.
For $h\in H^1(\R^3)$ the quadratic form of $L$ is
\[
 \langle Lh,h\rangle
 =\int_{\R^3}\bigl(|\nabla h|^2+(\lambda-qu^{q-1}+pu^{p-1})h^2\bigr)\,dx .
\]
Write $S(s)=s^{q-1}-s^{p-1}$, so that the profile equation reads
$-\Delta u+\lambda u=uS(u)$.
The next identity uses $\int uh=0$ to express the quadratic form in
terms of a radial integral of $uh$.

\begin{lemma}\label{lem:form}
Let $h$ be real, radial, smooth, compactly supported and tangent to the
mass constraint. For $r>0$, define
\[
 \eta(r)=\frac{2}{r^3u(r)^2}\int_0^r s^2u(s)h(s)\dd s .
\]
Set
\[
 y=-\frac{u'}u,\qquad
 K(r)=1+ry+\frac32r^2y'-\frac14r^2uS'(u).
\]
Then
\begin{equation}\label{eq:form}
\begin{aligned}
 \frac{\langle Lh,h\rangle}{4\pi}
 =\int_0^\infty r^2u^2\bigg[&
   \frac{r^2}{4}(\eta'')^2+K(r)(\eta')^2\\
 &-\frac r4\frac{\dd}{\dd r}
  \left[(7-3q)u^{q-1}+(3p-7)u^{p-1}\right]\eta^2\bigg]\dd r .
\end{aligned}
\end{equation}
\end{lemma}

\begin{proof}[\normalfont\proofname]
\noindent Step 1. We express $h$ in terms of $\eta$ and find the boundary values needed for integration by parts.
Let $D=3u/2+ru'$, the derivative of $t^{3/2}u(tr)$ at $t=1$.
The defining integral of $\eta$ gives
\[
 \left(\tfrac12 r^3u^2\eta\right)'=r^2uh,\qquad
 \left(\tfrac12 r^3u^2\right)'=r^2uD,\qquad
 h=D\eta+\tfrac12ru\eta'.
\]
If $\operatorname{supp}h\subset\overline B_R$, then for $r\ge R$,
\[
 \int_0^r s^2uh\dd s=\int_0^\infty s^2uh\dd s
 =\frac1{4\pi}\int_{\R^3}uh=0.
\]
Hence $\eta=0$ there. At zero the substitution $s=rt$ gives
\[
 \eta(r)=\frac2{u(r)^2}\int_0^1 t^2u(rt)h(rt)\dd t,\qquad
 \eta(0)=\frac{2h(0)}{3u(0)} .
\]
The right side extends smoothly and evenly across zero since $u,h$
are smooth radial functions and $u(0)>0$. Thus
\[
 \eta(r)=\eta(0)+O(r^2),\qquad \eta'(r)=O(r),\qquad \eta''(r)=O(1).
\]

\noindent Step 2. We use the equation for $u$ to compute the terms in $\langle Lh,h\rangle$ that contain $D$.
Let $f_1,f_2$ be smooth real radial functions, and assume that at least
one of them has compact support. By the definition of $L$ and
integration by parts,
\begin{align*}
 \langle L(uf_1),uf_2\rangle
 &=\int_{\mathbb R^3}\nabla(uf_1)\cdot\nabla(uf_2)\dd x
   +\int_{\mathbb R^3}
    (\lambda-qu^{q-1}+pu^{p-1})u^2f_1f_2\dd x \\
 &=\int_{\mathbb R^3}u^2\nabla f_1\cdot\nabla f_2\dd x
  +\int_{\mathbb R^3}(Lu)\,uf_1f_2\dd x .
\end{align*}
For the second equality we used the product identity
\[
 \nabla(uf_1)\cdot\nabla(uf_2)
 =u^2\nabla f_1\cdot\nabla f_2+\nabla u\cdot\nabla(uf_1f_2)
\]
and integrated the last term by parts. The profile equation in \eqref{eq:main} and the definition
$S(s)=s^{q-1}-s^{p-1}$ give
\begin{equation}\label{eq:Lu}
 Lu =(1-q)u^q+(p-1)u^p
   =-u^2S'(u).
\end{equation}
Finally, since all functions are radial,
$\nabla f_1\cdot\nabla f_2=f_1'(r)f_2'(r)$. Passing to radial
coordinates yields
\begin{equation}\label{eq:bilinear}
 \frac{\langle L(uf_1),uf_2\rangle}{4\pi}
 =\int_0^\infty r^2u^2
   \left[f_1'f_2'-uS'(u)f_1f_2\right]\dd r .
\end{equation}

Dividing $h=D\eta+ru\eta'/2$ by $u$ gives
\[
 g:=\frac hu=\frac Du\eta+\frac r2\eta',
 \qquad \frac Du=\frac32-ry,\qquad y=-u'/u .
\]
In \eqref{eq:bilinear}, apply the formula first to $h=ug$ and then to
the pair $D,\,2\eta h-D\eta^2$. Subtracting separates the derivative
and potential terms.
\begin{align*}
 &\frac{\langle Lh,h\rangle-\langle LD,2\eta h-D\eta^2\rangle}{4\pi}\\
 &=\int_0^\infty r^2u^2
 \left[(g')^2-\left(\frac Du\right)'
             \left(2\eta g-\frac Du\eta^2\right)'
       -uS'(u)\left(g-\frac Du\eta\right)^2\right]\dd r .
\end{align*}
The square in the potential term comes from
\[
 g^2-\frac Du\left(2\eta g-\frac Du\eta^2\right)
 =\left(g-\frac Du\eta\right)^2 .
\]
We use $D=\tfrac32u+ru'$ to compute the term containing $LD$. Since $ru'=x\cdot\nabla u$, the identity
\[
 \Delta(x\cdot\nabla u)=x\cdot\nabla(\Delta u)+2\Delta u
\]
and the profile equation in \eqref{eq:main} give
\[
 \begin{aligned}
 L(ru')
 &=x\cdot\nabla\bigl(-\Delta u+\lambda u-u^q+u^p\bigr)
   -2\Delta u\\
 &=-2\Delta u.
 \end{aligned}
\]
Using $D=3u/2+ru'$, the expression for $Lu$ in \eqref{eq:Lu}, and
$\Delta u=\lambda u-u^q+u^p$, we obtain
\begin{equation}\label{eq:LD}
 \frac{LD}{u}
 =-2\lambda+\frac{7-3q}{2}u^{q-1}
            +\frac{3p-7}{2}u^{p-1}.
\end{equation}

For this pairing, $h=D\eta+ru\eta'/2$ and $(r^3u^2/2)'=r^2uD$ give
\[
 r^2u(2\eta h-D\eta^2)
 =r^2uD\eta^2+r^3u^2\eta\eta'
 =\left(\frac{r^3u^2\eta^2}{2}\right)'.
\]
Therefore
\[
 \begin{aligned}
 \frac{\langle LD,2\eta h-D\eta^2\rangle}{4\pi}
 =\int_0^\infty\frac{LD}{u}
      \left(\frac{r^3u^2\eta^2}{2}\right)'\dd r
 =-\int_0^\infty\frac{r^3u^2}{2}
      \left(\frac{LD}{u}\right)'\eta^2\dd r .
 \end{aligned}
\]
The boundary term at zero vanishes because $LD/u$ and $\eta$ are
bounded there and $r^3u^2=O(r^3)$. At infinity it vanishes because
$\eta=0$ for $r\ge R$.

\noindent Step 3. We collect the coefficients of $(\eta')^2$ and $(\eta'')^2$.
Expanding the derivative terms gives
\begin{align*}
 &(g')^2-\left(\frac Du\right)'
          \left(2\eta g-\frac Du\eta^2\right)'\\
 &\qquad=
 \left[g'-\left(\frac Du\right)'\eta\right]^2
 -2\left(\frac Du\right)'\eta'
       \left(g-\frac Du\eta\right).
\end{align*}
Since $g=(D/u)\eta+r\eta'/2$ and $D/u=3/2-ry$, we have
\[
 g'-\left(\frac Du\right)'\eta
   =(2-ry)\eta'+\frac r2\eta'',
 \qquad
 g-\frac Du\eta=\frac r2\eta'.
\]
The derivative terms therefore contribute the integral of
\[
 r^2u^2\left[\left((2-ry)\eta'+\tfrac r2\eta''\right)^2
              +r(ry)'(\eta')^2\right].
\]
On expanding the square, the cross term is
\[
 \int_0^\infty r^3u^2(2-ry)\eta'\eta''\dd r
 =-\frac12\int_0^\infty
   \big[r^3u^2(2-ry)\big]'(\eta')^2\dd r .
\]
Here the boundary term vanishes at zero because $\eta'=O(r)$
and $2-ry$ is bounded. It vanishes at infinity because $\eta'$
is zero for $r>R$.

Since $(r^3u^2)'=r^2u^2(3-2ry)$, the coefficient of
$r^2u^2(\eta')^2$ from these derivative terms is
\begin{align*}
 &(2-ry)^2-\frac12(3-2ry)(2-ry)+\frac32r(ry)'\\
 &\qquad=\frac12(2-ry)+\frac32r(y+ry')
 =1+ry+\frac32r^2y'.
\end{align*}
The remaining potential term contributes
$-r^2uS'(u)/4$ to this coefficient, since
$g-(D/u)\eta=r\eta'/2$. Adding the pairing with $LD$ gives
\[
 \frac{\langle Lh,h\rangle}{4\pi}
 =\int_0^\infty r^2u^2
 \left[
   \frac{r^2}{4}(\eta'')^2+K(r)(\eta')^2
   -\frac r2\left(\frac{LD}{u}\right)'\eta^2
 \right]\dd r .
\]
Substituting the formula for $LD/u$ from \eqref{eq:LD} into the last term gives \eqref{eq:form}.
\end{proof}

The coefficient $K$ in \eqref{eq:form} is positive when $q<3$.

\begin{lemma}\label{lem:K}
Assume $1<q<3$ and $p>q$. Then $K(r)>0$ for $r>0$, and
$y'(r)>0$ wherever $S'(u(r))\le0$. At every $r>0$ we have
\begin{equation}\label{eq:key-K}
 K\ge\frac{3-q}{4}+ry+\frac{5-q}{4}r^2y'.
\end{equation}
\end{lemma}
\begin{proof}[\normalfont\proofname]
We use the sign of $S'(u)$ to control the two terms in $K$ involving
$y'$ and $S'(u)$. Where $S'(u)\le0$, both terms will be nonnegative.
Where $S'(u)\ge0$, integrating the radial equation from $r$ to
infinity gives a lower bound for $y'+r^{-2}$.
With $S(u)=u^{q-1}-u^{p-1}$ and $y=-u'/u$, the radial
equation~\eqref{eq:radial} gives
\begin{equation}\label{eq:ric}
 y'=\left(-\frac{u'}u\right)'=y^2-\frac2r y+S(u)-\lambda .
\end{equation}
Smooth radial symmetry gives $u'(0)=0$ and $u'(r)/r\to u''(0)$.
Passing to $r=0$ in the radial equation~\eqref{eq:radial} yields
\[
 -3u''(0)=u(0)(S(u(0))-\lambda),\qquad
 y'(0)=-\frac{u''(0)}{u(0)}=\frac{S(u(0))-\lambda}{3}\ge0 .
\]
The inequality follows from the maximum at the origin. To see that it
is strict, suppose $S(u(0))=\lambda$. Then $u$ and the constant $u(0)$
solve $-\Delta u=f(u)$ with $f(s)=s(S(s)-\lambda)$. Their nonnegative
difference $v=u(0)-u$ satisfies
\[
 -\Delta v=c(x)v,\qquad
 c(x)=\int_0^1 f'\big(u(x)+t(u(0)-u(x))\big)\dd t .
\]
The coefficient is bounded. Since $v(0)=0$ and
$-\Delta v+\|c\|_\infty v\ge0$, the strong maximum principle gives
$v\equiv0$, contradicting decay of $u$. We have proved
\begin{equation}\label{eq:center-strict}
 S(u(0))-\lambda=-3u''(0)/u(0)=3y'(0)>0 .
\end{equation}
The proof of \eqref{eq:center-strict} uses only the central maximum,
decay of $u$ and local Lipschitz continuity of $s(S(s)-\lambda)$ near $u(0)$.
Thus the same strict inequality also holds for $1<p<q<5$.

If $S'(u(r))\le0$, then $S'(u(s))\le0$ for $0\le s\le r$. Indeed,
\[
 S'(u)=u^{q-2}\big((q-1)-(p-1)u^{p-q}\big)
\]
and $(q-1)-(p-1)u^{p-q}$ increases with $r$ because $u$ decreases. At a first zero
$r_0\in(0,r]$ of $y'$, we would have $y''(r_0)\le0$. Differentiating the equation for $y$ in \eqref{eq:ric} gives
\[
 y''(r_0)=\frac{2y(r_0)}{r_0^2}+S'(u(r_0))u'(r_0)>0
\]
because $y>0$ and $u'<0$. Thus $y'>0$ on $[0,r]$. In this region
\[
 K=1+ry+\tfrac32r^2y'-\tfrac14r^2uS'(u)\ge1,\qquad
 uS'(u)\le0<(q+1)(y'+r^{-2}).
\]
Substituting this estimate for $uS'(u)$ in $K$ gives the lower bound~\eqref{eq:key-K}.

If $S'(u(r))\ge0$, then $S'(u(s))\ge0$ for all $s\ge r$ by the same
monotonicity. Set $v(s)=s u(s)$. The radial equation~\eqref{eq:radial} gives
\[
 v''=(\lambda-S(u))v,\qquad
 \big[(v')^2-(\lambda-S(u))v^2\big]'=S'(u)u'v^2.
\]
Since $v,v'\to0$ at infinity and $-v'/v=y-1/r$, integration gives
\begin{align}
 y'+r^{-2}
 &=\frac{(v')^2-vv''}{v^2}
 =\frac{\int_r^\infty S'(u(s))(-u'(s))s^2u(s)^2\dd s}{r^2u(r)^2}\notag\\
 &\ge \frac{\int_0^{u(r)} S'(t)t^2\dd t}{u(r)^2}
 =\frac{q-1}{q+1}u^{q-1}-\frac{p-1}{p+1}u^{p-1}
 \ge\frac{uS'(u)}{q+1}.\label{eq:ric-exterior}
\end{align}

The first inequality uses $s^2\ge r^2$ and $S'(u)(-u')\ge0$.
The last uses $p+1>q+1$ in the negative term. Substitution gives
\[
 K\ge 1+ry+\tfrac32 r^2y'-\tfrac{q+1}{4}(r^2y'+1)
 =\tfrac{3-q}{4}+ry+\tfrac{5-q}{4}r^2y',
\]
which proves the lower bound~\eqref{eq:key-K}.

It remains to check positivity where $S'(u(r))\ge0$.
If $S(u(r))\ge\lambda$, the equation for $y$ in \eqref{eq:ric} gives $r^2y'\ge(ry)^2-2ry$. Consequently
\[
 K\ge\frac{5-q}{4}
 \left(ry-\frac{3-q}{5-q}\right)^2
 +\frac{3-q}{2(5-q)}>0 .
\]
If $S(u(r))<\lambda$, then $S'(u(s))\ge0$ and $u'(s)<0$ imply
$S(u(s))\le S(u(r))<\lambda$ for $s\ge r$. Hence $v''>0$ on $[r,\infty)$.
Since $v>0$ decays, convexity gives $v'(r)<0$.
The integral estimate~\eqref{eq:ric-exterior} gives
$y'+r^{-2}\ge uS'(u)/(q+1)\ge0$. Hence
\[
 ry=1-rv'/v>1,\qquad r^2y'\ge-1,\qquad
 K\ge\frac{3-q}{4}+ry-\frac{5-q}{4}=ry-\tfrac12>0.
\]

\end{proof}

\subsection{Positivity on the mass constraint}
The sign of the mass derivative follows from strict positivity of the quadratic form on radial functions
tangent to the mass constraint.

\begin{proposition}\label{prop:tangent-positive}
Assume $1<q<p$, $q\le7/3$, and $0<\lambda<\lambda_*$. With
$u=u_\lambda$ and $L$ as in \eqref{eq:linearized-frequency}, we have
\[
 \langle Lh,h\rangle>0
 \qquad\text{for every }0\ne h\in H^1_{\rm rad}(\R^3)
 \text{ with }\int uh=0 .
\]
\end{proposition}
\begin{proof}[\normalfont\proofname]
First take $h$ smooth and compactly supported. Use $\eta$, $y$ and $K$
from Lemma~\ref{lem:form}. In \eqref{eq:form} the coefficient of
$\eta^2$ has the sign of
\[
 -\frac{d}{dr}\big[(7-3q)u^{q-1}+(3p-7)u^{p-1}\big].
\]
Lemma~\ref{lem:K} gives $K>0$ because $q\le7/3<3$.

If $p\ge7/3$, both $7-3q$ and $3p-7$ are nonnegative, and $p>q$
makes at least one of them positive. Since $u'(r)<0$ for $r>0$,
\[
 \frac{d}{dr}\big[(7-3q)u^{q-1}+(3p-7)u^{p-1}\big]<0.
\]
All three terms in \eqref{eq:form} are therefore nonnegative, and
the coefficient of $\eta^2$ is strictly positive for $r>0$.

If $p<7/3$, we write
\[
 (7-3q)u^{q-1}+(3p-7)u^{p-1}
 =(7-3p)S(u)+3(p-q)u^{q-1}.
\]
The last power contributes the strictly positive coefficient
$\tfrac34(q-1)(p-q)r y u^{q-1}$ to the $\eta^2$ term inside the brackets.
We estimate the contribution of $S(u)$ by using
$(y')_+=\max\{y',0\}$. We claim that, almost everywhere on $(0,\infty)$,
\begin{equation}\label{eq:S-derivative}
 (S(u))'\le \big((y')_+\big)'
             +\left(\frac2r-2y\right)(y')_+ .
\end{equation}
Indeed, on $\{y'>0\}$, differentiating the equation for $y$ in \eqref{eq:ric} gives
\[
 (S(u))'=y''+(2/r-2y)y'-2y/r^2
 \le y''+(2/r-2y)y'.
\]
On $\{y'\le0\}$, Lemma~\ref{lem:K} gives $S'(u)>0$, hence
$(S(u))'=S'(u)u'<0$. The right side of the derivative bound~\eqref{eq:S-derivative} vanishes
almost everywhere there, since $((y')_+)'=\mathbf1_{\{y'>0\}}y''$ almost
everywhere for locally $C^1$ functions $y'$.

Multiply the derivative bound~\eqref{eq:S-derivative} by $-r^3u^2\eta^2$ and integrate. The boundary
term at zero vanishes since $y'$ and $\eta$ are bounded and
$r^3u^2=O(r^3)$. At large radii $\eta=0$ by $\int uh=0$. Integration
by parts and $(r^3u^2)'=(3-2ry)r^2u^2$ give
\[
\begin{aligned}
 -\int_0^\infty r^3u^2(S(u))'\eta^2\dd r
 &\ge \int_0^\infty
 \left[(r^3u^2\eta^2)'
       -r^3u^2\left(\frac2r-2y\right)\eta^2\right](y')_+\dd r\\
 &=\int_0^\infty r^2u^2(y')_+
       \big(\eta^2+2r\eta\eta'\big)\dd r\\
 &\ge -\int_0^\infty r^4u^2(y')_+(\eta')^2\dd r .
\end{aligned}
\]
The last inequality follows from $(\eta+r\eta')^2\ge0$. Since $7-3p>0$,
the quadratic form identity~\eqref{eq:form} gives
\[
\begin{aligned}
 \frac{\langle Lh,h\rangle}{4\pi}
 \ge\int_0^\infty r^2u^2\bigg[&
 \frac{r^2}{4}(\eta'')^2
 +\left(K-\frac{7-3p}{4}r^2(y')_+\right)(\eta')^2\\
 &+\frac{3(q-1)(p-q)}{4}r y u^{q-1}\eta^2\bigg]\dd r .
\end{aligned}
\]

It remains to check the coefficient of $(\eta')^2$. Where $y'\le0$,
this coefficient equals $K>0$. Where $y'>0$, the lower bound
\eqref{eq:key-K} gives
\[
 K-\frac{7-3p}{4}r^2y'
 \ge\frac{3-q}{4}+ry+\frac{3p-q-2}{4}r^2y'>0 ,
\]
since $q<3$ and $3p-q-2=2(p-1)+(p-q)>0$.

In both cases the coefficient of $\eta^2$, including the weight $r^2u^2$,
is continuous and strictly positive for $r>0$. Hence for every
$0<a<b<\infty$ there is a constant $c_{a,b}>0$, independent of $h$, such that
\begin{equation}\label{eq:local-tangent}
 \frac{\langle Lh,h\rangle}{4\pi}\ge c_{a,b}\int_a^b\eta(r)^2\dd r .
\end{equation}

For a general radial $h\in H^1$ with $\int uh=0$, take
$\widetilde h_n\in C^\infty_{c,\rm rad}$ converging to $h$ in $H^1$.
Choose a fixed radial test function $\phi$ with $\int u\phi=1$. Then
\[
 h_n=\widetilde h_n-\left(\int u\widetilde h_n\right)\phi,\qquad
 \int uh_n=0,\qquad h_n\longrightarrow h\quad\text{in }H^1,
\]
because $\int u\widetilde h_n\to\int uh=0$. Since the potential is bounded,
\[
 \langle Lh_n,h_n\rangle\longrightarrow\langle Lh,h\rangle.
\]
On each interval $[a,b]\subset(0,\infty)$, the primitives satisfy
\begin{align*}
 \sup_{a\le r\le b}|\eta_n(r)-\eta(r)|
 &\le\frac{2}{\min_{[a,b]}r^3u(r)^2}
       \left(\int_0^b s^2u^2\dd s\right)^{1/2}
       \left(\int_0^b s^2|h_n-h|^2\dd s\right)^{1/2}
 \longrightarrow0 .
\end{align*}
We may therefore pass to the limit in the local lower bound~\eqref{eq:local-tangent}.
If $\langle Lh,h\rangle=0$, then $\int_a^b\eta^2=0$ on every such
interval. By continuity of the primitive,
\[
 \int_0^r s^2u(s)h(s)\dd s=0\quad(r>0)
 \quad\Longrightarrow\quad r^2u(r)h(r)=0\quad\text{a.e.}
\]
As $u>0$, this gives $h=0$.

\end{proof}

\begin{proof}[\normalfont Proof of Theorem~\ref{thm:increasing}]
Suppose $q\le7/3$ and $M'(\lambda)=0$ at an admissible frequency.
Set $u=u_\lambda$ and $w=\partial_\lambda u_\lambda$.
By \eqref{eq:linearized-frequency}, $Lw=-u$,
so $w\ne0$ and
\[
 \int uw=\tfrac12 M'(\lambda)=0,\qquad
 \langle Lw,w\rangle=-\int uw=0 .
\]
Proposition~\ref{prop:tangent-positive} gives $\langle Lw,w\rangle>0$,
a contradiction.
The continuous function $M'$ therefore has one sign on $(0,\lambda_*)$. Since $M(\lambda)\to\infty$
as $\lambda\uparrow\lambda_*$ by the endpoint limits~\eqref{eq:ends}, that sign must be positive.

For $q\le7/3$ the limits in \eqref{eq:ends} now give the exact ranges
\[
 M((0,\lambda_*))=
 \begin{cases}
 (0,\infty),&1<q<7/3,\\
 (m_Q,\infty),&q=7/3 .
 \end{cases}
\]
Each mass in the indicated interval gives one frequency and hence
one positive solution up to translations.

Finally suppose $q>7/3$. By the endpoint limits~\eqref{eq:ends}, $M$ tends to infinity at
both ends of $(0,\lambda_*)$. Thus the continuous function $M$ attains
its minimum $m$ at an interior frequency. For each $a>m$, the intermediate
value theorem gives at least one frequency on each side of a minimizing frequency.
These two frequencies give distinct solutions. This proves part~(3) and
completes the uniqueness criterion.
\end{proof}

\section{\texorpdfstring{One minimum for $p>q$ and $q$ just above $7/3$}{One minimum for p > q and q just above 7/3}}
We now let the lower power $q$ approach $7/3$ from above. The bound
on $q-7/3$ will hold even when $p-q$ tends to zero.

\begin{theorem}\label{thm:uniform-fold}
For every $D>0$ there exists $e_D>0$ such that if
\[
 0<p-q\le D,\qquad 0<q-7/3<e_D,
\]
the mass has exactly one critical frequency $\lambda_c$. We have
\[
 M'<0\quad(0<\lambda<\lambda_c),\qquad
 M''(\lambda_c)>0,\qquad
 M'>0\quad(\lambda_c<\lambda<\lambda_*).
\]
Let $m=M(\lambda_c)$. For every $0<a<m$ there are no positive solutions, for $a=m$ there is one, and for every
$a>m$ there are exactly two, up to translations.
\end{theorem}

We use an increasing parameter $k$ to locate the first zero of $M'$.
For the remaining frequencies we prove strict positivity of the
quadratic form on the mass constraint to exclude any further zero.

\subsection{Rescaling the mass and its derivative}\label{sec:gap-scale}
We need an estimate for $M'$ that is uniform as $p-q$ tends to zero. Write $e=q-7/3>0$ and
$d=p-q\in(0,D]$. Fix an admissible frequency $\lambda>0$.
We introduce an amplitude scale $s>0$ by writing
\[
 u_\lambda(x)=s\,v(\sqrt{\lambda}\,x).
\]
Substituting into the profile equation in \eqref{eq:main} and dividing
by $s\lambda$ gives
\begin{equation}\label{eq:gap-before-normalization}
 -\Delta v+v=\frac{s^{q-1}}{\lambda}v^q-\frac{s^{p-1}}{\lambda}v^p.
\end{equation}
We choose $s$ so that the right-hand side equals $1$ at $v=1$.
Since $d=p-q$, this condition is
\begin{equation}\label{eq:amplitude}
 \frac{s^{q-1}-s^{p-1}}{\lambda}=1
 \qquad\Longleftrightarrow\qquad
 \lambda=s^{q-1}(1-s^d).
\end{equation}
With this choice the right-hand side of \eqref{eq:gap-before-normalization}
takes the form $v^q-(s^{p-1}/\lambda)(v^p-v^q)$.

To choose the smaller root of \eqref{eq:amplitude}, we use
\[
 \frac{\dd}{\dd s}\big[s^{q-1}(1-s^d)\big]
 =s^{q-2}\big[(q-1)-(p-1)s^d\big].
\]
Thus $s\mapsto s^{q-1}(1-s^d)$ vanishes at $0$ and $1$ and has a
single maximum in $(0,1)$. At the origin, the radial equation~\eqref{eq:radial} gives
\begin{equation}\label{eq:central-value}
 u_\lambda(0)^{q-1}\bigl(1-u_\lambda(0)^d\bigr)-\lambda
 =-\frac{3u_\lambda''(0)}{u_\lambda(0)}>0.
\end{equation}
The strict inequality follows from \eqref{eq:center-strict} applied to
$S(s)=s^{q-1}(1-s^d)$.

By \eqref{eq:central-value} and $\lambda>0$, the amplitude equation~\eqref{eq:amplitude} has one root in
$(0,u_\lambda(0))$ and one in $(u_\lambda(0),1)$. We take the smaller
root $s=s_\lambda$. It lies on the increasing part of
$s\mapsto s^{q-1}(1-s^d)$, where $(p-1)s^d<q-1$. For fixed $p,q$,
this is the root that tends to zero as $\lambda\downarrow0$.
In particular, $v(0)=u_\lambda(0)/s>1$. Set
\begin{equation}\label{eq:gap-parameter}
 k=\frac{d s^d}{1-s^d}\in(0,q-1),\qquad
 s^d=\frac{k}{k+d},\qquad
 \lambda(k)=\frac d{k+d}\left(\frac{k}{k+d}\right)^{(q-1)/d}.
\end{equation}
Logarithmic differentiation gives
\begin{equation}\label{eq:gap-log-frequency}
 \frac{\dd\log\lambda}{\dd k}=\frac{q-1-k}{k(k+d)}>0.
\end{equation}
By the amplitude equation~\eqref{eq:amplitude} and the definition
of $k$ in \eqref{eq:gap-parameter},
\[
 \frac{s^{p-1}}{\lambda}=\frac{s^d}{1-s^d}=\frac{k}{d}.
\]
The rescaled profile equation~\eqref{eq:gap-before-normalization} therefore becomes
\begin{equation}\label{eq:gap-profile}
 -\Delta v+v=v^q-kv^q\frac{v^d-1}{d}.
\end{equation}
Let $H(e,d,k)=\int_{\R^3}v^2$. Changing variables in the mass integral
gives $M(\lambda(k))=s^2\lambda^{-3/2}H$. Substituting the formulas
for $s$ and $\lambda$ from \eqref{eq:gap-parameter} and using $q-1=4/3+e$ gives
\[
 M(\lambda(k))=
 \left(1+\frac{k}{d}\right)^{3/2}
 \left(1+\frac{d}{k}\right)^{3e/(2d)} H(e,d,k).
\]
In particular, the explicit factor satisfies
\begin{equation}\label{eq:gap-factor}
 \frac{\dd}{\dd k}\log\!\left(\frac{M(\lambda(k))}{H(e,d,k)}\right)
 =\frac{3}{2(k+d)}-\frac{3e}{2k(k+d)}
 =\frac{3(k-e)}{2k(k+d)} .
\end{equation}

Where $H$ is differentiable in $k$, the product rule and
\eqref{eq:gap-factor} give
\begin{equation}\label{eq:gap-mass-derivative}
 \frac{\dd}{\dd k}M(\lambda(k))
 =\frac{M(\lambda(k))}{k(k+d)H}\,F(e,d,k),
\end{equation}
where
\begin{equation}\label{eq:gap-derivative}
 F(e,d,k)=\tfrac32(k-e)H+k(k+d)H_k,\qquad H_k=\partial_k H .
\end{equation}
Both $M/(k(k+d)H)$ and $\lambda'(k)$ are positive. Dividing
\eqref{eq:gap-mass-derivative} by $\lambda'(k)$ shows that
$M'(\lambda(k))$ has the same sign as $F(e,d,k)$.

\subsection{Exactly one minimum for small \texorpdfstring{$k$}{k}}\label{sec:gap-small}
The formula for $F$ contains $H_k$. To determine the sign of $F_k$
uniformly as $d\downarrow0$, we first differentiate the rescaled profile twice.

\begin{lemma}\label{lem:forward-profile}
Fix $D>0$. For $d\in[0,D]$ and $|e|,|k|<\delta$, with
$\delta>0$ depending only on $D$, equation~\eqref{eq:gap-profile}
has a locally unique radial profile near $Q$ in $H^2$.
The profile is positive, equals $Q$ at $e=k=0$, and is $C^2$ in
$(e,d,k)$ up to $d=0,D$. The same regularity holds for $H=\int v^2$.
For $e,d,k>0$ the profile rescales to $u_{\lambda(k)}$.
\end{lemma}
\begin{proof}[\normalfont\proofname]
At $e=k=0$ the equation reduces to $-\Delta Q+Q=Q^{7/3}$ for every $d$.

For real radial $v\in H^2(\R^3)$ use the signed power
$|v|^{q-1}v$. The potentially singular term in \eqref{eq:gap-profile}
has the extension
\begin{equation}\label{eq:gap-divided-difference}
 \frac{|v|^{q+d-1}v-|v|^{q-1}v}{d}
 =\int_0^1 |v|^{q+td-1}v\log |v|\dd t .
\end{equation}
Use the integral at $d=0$, with value zero where $v=0$. We work near $e=k=0$, with $d\in[0,D]$. Since $q=7/3+e$,
we may restrict to $|e|<1/6$. Then, for every $d\in[0,D]$
and $t\in[0,1]$,
\[
 a=q+td\ge q=\frac73+e\ge\frac{13}{6}=2+\frac16.
\]
For $z\ne0$, the integrand in \eqref{eq:gap-divided-difference} is
$z|z|^{a-1}\log|z|$, where $a=q+td$. We set its value at $z=0$
equal to zero. Each derivative in $a$ adds a factor of $\log|z|$.
Each derivative in $z$ lowers the power of $|z|$ by one.
Thus, for $0<|z|\le1$ and nonnegative integers $i,j$ with $i+j\le2$,
\[
 \left|
 \partial_z^i\partial_a^j\bigl(z|z|^{a-1}\log|z|\bigr)
 \right|
 \le C_D |z|^{a-i}\bigl(1+|\log|z||^{j+1}\bigr).
\]
The constant is uniform because $2+1/6\le a\le 5/2+D$.
Since $i\le2$ and $j+1\le3$, the right-hand side is bounded by
\[
 C_D |z|^{1/6}\bigl(1+|\log|z||^3\bigr)
 \longrightarrow 0
 \qquad\text{as }z\to0.
\]
Thus the integrand and its derivatives of orders one and two extend
continuously to $z=0$, with value zero there.

Let $v$ vary in a fixed $H^2_{\rm rad}$ ball about $Q$.
The embedding $H^2(\R^3)\hookrightarrow L^\infty(\R^3)$ gives a
number $R$ such that $|v(x)|\le R$ for every $v$ in the ball. For $j=0,1,2$, differentiating the
integrand $j$ times in $a$ gives
\[
 v(x)\,|v(x)|^{a-1}\big(\log|v(x)|\big)^{j+1}.
\]
The factor multiplying $v(x)$ extends continuously to zero and is
uniformly continuous for $|z|\le R$ and $13/6\le a\le5/2+D$. As $(v,a)$ varies in $H^2\times\R$, this factor varies
continuously in $L^\infty$ and $v$ varies continuously in $L^2$.
Their product therefore varies continuously in $L^2$.

One derivative in $v$ multiplies a direction $h$ by the derivative of
the scalar integrand in $z$. A mixed derivative in $a$ and $v$ uses
the scalar derivative in $a$ and $z$ in the same way. Two derivatives
in $v$ multiply $h_1h_2$ by the second scalar derivative in $z$.
All three coefficients vary continuously in $L^\infty$ by uniform
continuity on the same closed rectangle, and
\[
 \|h_1h_2\|_2\le\|h_1\|_\infty\|h_2\|_2
 \le C\|h_1\|_{H^2}\|h_2\|_{H^2}.
\]
These bounds give continuity of the derivatives in operator norm.
The fundamental theorem of calculus along a line segment in the
parameters and in $v$ gives the corresponding Fr\'echet derivatives.
Since $a=7/3+e+td$, a derivative in $d$ contributes a factor $t\le1$.
All the bounds are uniform for $0\le t\le1$, so the integral in
\eqref{eq:gap-divided-difference} can be differentiated twice.
The same estimates apply to the signed power $|v|^{q-1}v$ with the
logarithmic factor omitted. Thus the right side of
\eqref{eq:gap-profile} is $C^2$ with values in $L^2_{\rm rad}$.

At $e=k=0$, the derivative of the profile equation in $v$ is the same for every $d$, namely
\[
 L_Q=-\Delta+1-\tfrac73 Q^{4/3}.
\]
Multiplication by $Q^{4/3}$ is compact from $H^2$ to $L^2$,
and the radial kernel of $L_Q$ is zero
\cite[Section~4.1]{LRN}, \cite{Weinstein}.
Since $-\Delta+1:H^2_{\rm rad}\to L^2_{\rm rad}$ is an isomorphism,
the Fredholm alternative gives a bounded inverse for $L_Q$.

Fix $d_0\in[0,D]$. For $d$ close to $d_0$ and $|e|$ small,
$7/3+e+td$ stays uniformly above two for $0\le t\le1$. Thus
\eqref{eq:gap-divided-difference} defines a $C^2$ equation on an open
neighborhood of $(0,d_0,0)$, including when $d_0=0$ or $d_0=D$. The implicit function theorem gives a $C^2$ profile near
$(0,d_0,0)$, unique in a ball about $Q$.

Cover $[0,D]$ by finitely many open intervals in $d$ whose closures
lie in neighborhoods where the implicit function theorem applies.
Choose $\rho>0$ so that $B_{H^2}(Q,\rho)$ lies in every ball of local
uniqueness. On the closure of each interval $I$, continuity and
$v(0,d,0)=Q$ give
\[
 \sup_{d\in\overline I}\|v(e,d,k)-Q\|_{H^2}\longrightarrow0
 \qquad\text{as }(e,k)\to(0,0).
\]
There are only finitely many intervals, so one choice of $\delta>0$
makes every local profile lie in $B_{H^2}(Q,\rho)$ for
$|e|,|k|<\delta$. On an overlap both profiles solve the same equation
in this ball and hence agree. This gives one $C^2$ profile on the
whole parameter range. Since $H=\int v^2$, the mass is also $C^2$.

To prove positivity, let $v_-=\max\{-v,0\}$. Since $Q>0$ and
$H^2(\R^3)$ embeds in $L^\infty(\R^3)$,
\[
 \|v_-\|_\infty\le\|v-Q\|_\infty\le C\|v-Q\|_{H^2}\longrightarrow0
 \quad\text{as }(e,k)\to(0,0),
\]
uniformly for $d\in[0,D]$.

Write $f(v)$ for the right-hand side of the rescaled
equation~\eqref{eq:gap-profile}, using signed powers.
By the integral formula~\eqref{eq:gap-divided-difference},
for $0<|z|\le1$ we have
\[
 \begin{aligned}
 \left|\frac{f(z)}{z}\right|
 &\le |z|^{q-1}
   +|k|\int_0^1 |z|^{q+td-1}|\log|z||\,\dd t \\
 &\le |z|\bigl(1+|\log|z||\bigr),
 \end{aligned}
\]
because $q\ge2$, $d\ge0$ and $|k|\le1$ near $e=k=0$.
The last bound tends to zero as $z\to0$. Choose $\varepsilon>0$
so that $|f(z)|\le |z|/2$ whenever $|z|\le\varepsilon$, and take
$(e,k)$ close enough to $(0,0)$ that $\|v_-\|_\infty\le\varepsilon$.

Testing $-\Delta v+v=f(v)$ with $-v_-$ now gives
\[
 \int_{\mathbb R^3}\bigl(|\nabla v_-|^2+v_-^2\bigr)\,\dd x
 =-\int_{\{v<0\}} f(v)v_-\,\dd x
 \le \frac12\int_{\mathbb R^3}v_-^2\,\dd x .
\]
Thus $v_-=0$, so $v\ge0$. Also $v(0)>0$ by uniform closeness to $Q$.
Since $f(z)/z$ extends continuously to zero, the strong maximum
principle gives $v>0$ in $\mathbb R^3$.

For $e,d,k>0$
in this neighborhood, rescaling back gives the unique centered
solution at $\lambda(k)$. Thus all sufficiently small positive
$k$ correspond to admissible frequencies.

\end{proof}

The sign of $F_k$ at $e=k=0$ determines the shape of the mass for small $k$.

\begin{lemma}\label{lem:forward-small}
For every $D>0$ there are $k_0,e_D>0$ such that, for $0<d\le D$ and
$0<e<e_D$, the function $F(e,d,\cdot)$ has exactly one zero $k_c$ in
$(0,k_0]$. It is negative before $k_c$ and positive after $k_c$ up to
$k_0$, and $M''(\lambda(k_c))>0$.
\end{lemma}
\begin{proof}[\normalfont\proofname]
At $k=0$ the definition~\eqref{eq:gap-derivative} gives
$F(e,d,0)=-3eH(e,d,0)/2<0$ for $e>0$. We now estimate $F_k$
by first computing its value at $e=k=0$.

For $d>0$, differentiating \eqref{eq:gap-profile} in $k$ at $e=k=0$ gives
\[
 d L_Q v_k=-(Q^{7/3+d}-Q^{7/3}).
\]
We need
$dH_k=2d\int Qv_k$, which we can compute without solving for $v_k$.
For each $\omega>0$, the scaled profile $\omega^{3/4}Q(\sqrt\omega\,x)$
solves $-\Delta v+\omega v=v^{7/3}$.
Differentiating this last equality in $\omega$ at $\omega=1$ gives
\begin{equation}\label{eq:critical-dilation}
 Z=\tfrac34Q+\tfrac12x\cdot\nabla Q,\qquad
 L_QZ=-Q,\qquad L_Q=-\Delta+1-\tfrac73Q^{4/3}.
\end{equation}
Integration by parts gives
\[
 dH_k=-2d\int (L_QZ)v_k
     =2\int Z(Q^{7/3+d}-Q^{7/3}).
\]
For $a>1$, integration of $x\cdot\nabla(Q^{a+1})$ gives
\begin{equation}\label{eq:dilation-moment}
 2\int ZQ^a
 =\tfrac32\int Q^{a+1}+\frac1{a+1}\int x\cdot\nabla(Q^{a+1})
 =\frac{3(a-1)}{2(a+1)}\int Q^{a+1}.
\end{equation}
For the term of power $7/3$, use $L_QQ=-\tfrac43Q^{7/3}$ and
the dilation identity $L_QZ=-Q$ in \eqref{eq:critical-dilation} to obtain
\begin{equation}\label{eq:critical-moment}
 2\int ZQ^{7/3}=-\tfrac32\int QL_QZ=\tfrac32m_Q.
\end{equation}
Differentiating \eqref{eq:gap-derivative} at $e=k=0$ gives
\[
 \begin{aligned}
 F_k(0,d,0)
 &=\tfrac32m_Q+dH_k(0,d,0)\\
 &=\tfrac32m_Q+2\int ZQ^{7/3+d}-2\int ZQ^{7/3}
 =2\int ZQ^{7/3+d}.
 \end{aligned}
\]
The cancellation in the last equality follows from \eqref{eq:critical-moment}. Hence
\begin{equation}\label{eq:gap-positive}
 F_k(0,d,0)=\frac{3(4/3+d)}{2(10/3+d)}
                  \int Q^{10/3+d}>0.
\end{equation}
At $d=0$ the formula holds by continuity and equals $3m_Q/2$.

The derivative $F_k$ is continuous in $(e,d,k)$ by the $C^2$ dependence
of $H$ in Lemma~\ref{lem:forward-profile}. By \eqref{eq:gap-positive}, $F_k(0,d,0)$ has a
positive minimum on $[0,D]$. Compactness then gives $k_0,e_D,c_D>0$ such that
\[
 F_k(e,d,k)\ge c_D\qquad
 (0\le e\le e_D,\ 0\le d\le D,\ 0\le k\le k_0).
\]

Since $H(e,d,0)$ is bounded there,
\[
 F(e,d,k_0)\ge -\tfrac32 e H(e,d,0)+c_Dk_0>0
\]
if $e_D$ is chosen so that
$\tfrac32 e_D\sup H(e,d,0)<c_Dk_0$, with the supremum taken over the
closed parameter rectangle above. For $e>0$, $F$ therefore increases from a
negative value at $0$ to a positive value at $k_0$. It has exactly one
zero $k_c$ between them. Since $\lambda'(k)>0$, we have
\begin{equation}\label{eq:fold-small-sign}
 M'(\lambda(k))<0\quad(0<k<k_c),\qquad
 M'(\lambda(k_c))=0,\qquad
 M'(\lambda(k))>0\quad(k_c<k\le k_0).
\end{equation}
Set $\lambda_c=\lambda(k_c)$. Since $M'(\lambda_c)=0$, the chain rule gives
\[
 \left.\frac{d^2}{dk^2}M(\lambda(k))\right|_{k=k_c}
 =M''(\lambda_c)\lambda'(k_c)^2 .
\]
When we differentiate \eqref{eq:gap-mass-derivative} at $k_c$, the
term containing $F(e,d,k_c)$ vanishes. We obtain
\[
 M''(\lambda_c)=
 \left(\frac{\dd k}{\dd\lambda}(\lambda_c)\right)^2
 \frac{M(\lambda_c)F_k(e,d,k_c)}{k_c(k_c+d)H(e,d,k_c)}>0.
\]
\end{proof}

\subsection{Positivity at the remaining frequencies}\label{sec:gap-outer}
We now exclude stationary points beyond the interval in
Lemma~\ref{lem:forward-small}.

\begin{proposition}\label{prop:forward-outer}
Fix $D>0$ and $k_0$ as in Lemma~\ref{lem:forward-small}.
For $0<d\le D$, $e>0$ sufficiently small in terms of $D$, and every
admissible $k\ge k_0$, the linearized operator $L$ at the rescaled
profile $v$ satisfies
\begin{equation}\label{eq:fold-tangent-positive}
 \langle Lh,h\rangle>0\qquad
 (0\ne h\in H^1_{\rm rad},\ \int vh=0).
\end{equation}
\end{proposition}
\begin{proof}[\normalfont\proofname]
\noindent Step 1. We write the quadratic form in the rescaled variables and bound $K$ from below.
We have $k_0\le k<q-1$ by \eqref{eq:gap-parameter}. Decrease $e_D$
so that $7/3<q=7/3+e\le5/2$.
With $u_\lambda(x)=s\,v(\sqrt\lambda\,x)$, the rescaled coefficient is
\[
 N(z)=\frac{S(sz)}{\lambda}
     =z^{q-1}\left(1-k\frac{z^d-1}{d}\right).
\]
Then \eqref{eq:gap-profile} reads $-\Delta v+v=vN(v)$, and its
linearized operator is $L=-\Delta+1-N(v)-vN'(v)$.
The radial equation for $v$ is
\begin{equation}\label{eq:gap-radial}
 -(r^2v')'=r^2v\big(N(v)-1\big),\qquad r>0.
\end{equation}
For a smooth compactly supported radial
$h$ with $\int vh=0$, recall the three quantities in the radial form,
\[
 y=-\frac{v'}v,\qquad
 K=1+ry+\tfrac32r^2y'-\tfrac14 r^2vN'(v),\qquad
 \eta(r)=\frac{2}{r^3v(r)^2}\int_0^r t^2v(t)h(t)\dd t .
\]
Applying the radial form identity~\eqref{eq:form} after the change of radius gives
\begin{equation}\label{eq:scaled-radial-form}
 \frac{\langle Lh,h\rangle}{4\pi}
 =\int_0^\infty r^2v^2\left[
   \tfrac14r^2(\eta'')^2+K(r)(\eta')^2
   +\tfrac34 r y v^{q-1} B(v)\eta^2\right]\dd r ,
\end{equation}
where the last term follows from
\begin{equation}\label{eq:scaled-dilation}
\begin{aligned}
 &-\frac r4\frac{\dd}{\dd r}
 \left[(7-3q)(1+k/d)v^{q-1}+(3p-7)(k/d)v^{p-1}\right]\\
 &\qquad=\tfrac34r y v^{q-1}B(v),\\
 & B(v)=k(p-1)v^d+ek(p-1)\frac{v^d-1}{d}+e(k-q+1).
\end{aligned}
\end{equation}
To compute $B(v)$, use
$7-3q=-3e$, $3p-7=3(e+d)$ and $v'=-yv$ in the first line.
After removing the factor $3ryv^{q-1}/4$, the remaining expression is
\[
 -e(q-1)(1+k/d)+(e+d)(p-1)(k/d)v^d.
\]
Collecting the terms with $v^d-1$ gives \eqref{eq:scaled-dilation}.

To transfer the estimate for $K$ to the rescaled profile, write
$u_\lambda(R)=s v(r)$, $r=\sqrt\lambda R$. Then
\[
 -\frac{u_\lambda'(R)}{u_\lambda(R)}=\sqrt\lambda\,y(r),\qquad
 R^2\left(-\frac{u_\lambda'}{u_\lambda}\right)'(R)=r^2 y'(r).
\]
For the potential term, $s^{q-1}/\lambda=1+k/d$ and
$s^{p-1}/\lambda=k/d$ give
\[
 (q-1)(sv)^{q-1}-(p-1)(sv)^{p-1}=\lambda v N'(v).
\]
Since $R^2\lambda=r^2$, each term of $K$ is unchanged by the scaling.
In the proof of Lemma~\ref{lem:K} the three lower bounds
were $1$, $(3-q)/(2(5-q))$, and $ry-1/2$ with $ry>1$. Their minimum is
at least $1/10$ for $7/3\le q\le5/2$. We also retain the more precise
bound $K\ge ry-1/2$ whenever $N'(v)\ge0$ and $N(v)<1$.

\noindent Step 2. We estimate the exterior integral where $v$ is small.
First we give a uniform lower bound for $v(2)$. Since
$sv=u_\lambda<1$, we have $N(v)=S(sv)/\lambda\ge0$. For $0<v\le1$,
\[
 v^{2-q}N'(v)=(q-1)\left(1-k\frac{v^d-1}{d}\right)-kv^d
 \ge q-1-k>0.
\]
Since $N(1)=1$, this gives $N(v)\le1$ for $v\le1$.
Also $v(0)=u_\lambda(0)/s>1$ by the choice of the smaller root in
\eqref{eq:amplitude}. For $v\ge1$, the definition gives $N(v)\le v^{q-1}$. In both
regions $v(N(v)-1)\le v(0)^q$, so the rescaled radial equation~\eqref{eq:gap-radial} gives
\[
 -(r^2v')'\le r^2v(0)^q
 \quad\Longrightarrow\quad
 -v'(r)\le rv(0)^q/3
 \quad\Longrightarrow\quad
 v(r)\ge v(0)-r^2v(0)^q/6 .
\]
Choose $r_1=v(0)^{-(q-1)/2}\le1$. Then
$r_1v(r_1)\ge(5/6)v(0)^{(3-q)/2}\ge5/6$.

Write $w=rv$. The rescaled radial equation~\eqref{eq:gap-radial},
the bound $N(v)\ge0$ and the decay of $w,w'$ at infinity give
\[
 w''=(1-N(v))w\le w,\qquad
 \big[e^{-r}(w'+w)\big]'\le0,\qquad w'+w\ge0.
\]
In particular $e^r w(r)$ is increasing. It follows that
\[
 v(2)=w(2)/2\ge e^{r_1-2}w(r_1)/2\ge (5/12)e^{-2}.
\]

For $0<z\le1$, $(1-z^d)/d\le|\log z|$. Hence
$0\le N(z)\le z^{q-1}(1+k|\log z|)$ there. Since
$4/3\le q-1\le3/2$ and $k<q-1$, choose $b>0$ so small that
$2b<(5/12)e^{-2}$ and $N(z)\le1/2$ for $0<z\le2b$, for all our
parameters. Radial decrease now gives $r>2$ wherever $v(r)\le2b$.
On this entire exterior interval $w''\ge w/2$, so
$e^{-r/\sqrt2}(w'+w/\sqrt2)$ increases to zero. Combining
$w'+w/\sqrt2\le0$ with $w'+w\ge0$ gives
\[
 1/\sqrt2\le -w'/w\le1,\qquad
 1/\sqrt2+1/r\le y\le1+1/r .
\]
Set $W=r^3v^2=rw^2$. For $v\le2b$ we have $r>2$, so
\[
 \frac{W'}{W}=\frac1r+2\frac{w'}w\le-\frac1{\sqrt2},
 \qquad
 K\ge ry-\frac12\ge\frac r{\sqrt2}.
\]
In particular, $r^2v^2K\ge W/\sqrt2$.

To control $\int_{v\le b}W\eta^2\dd r$, choose
a smooth function $\chi$ of $v$ equal to one on $v\le b$ and zero on
$v\ge2b$, with derivative in $v$ bounded by $C/b$. Its radial derivative
is bounded by $C$ since $|v'|=yv\le (1+1/r)2b$ on $b\le v\le2b$.
On the support of $\chi$, $W\le-\sqrt2 W'$. Integration by parts gives
\begin{align*}
 \int_0^\infty W(\chi\eta)^2\dd r
 &\le-\sqrt2\int_0^\infty W'(\chi\eta)^2\dd r
 =2\sqrt2\int_0^\infty W\chi\eta\,(\chi\eta)'\dd r\\
 &\le2\sqrt2\left(\int_0^\infty W(\chi\eta)^2\dd r\right)^{1/2}
             \left(\int_0^\infty W((\chi\eta)')^2\dd r\right)^{1/2}.
\end{align*}
Here $\chi$ vanishes near zero and $\eta$ vanishes at large radii by
$\int vh=0$. Squaring the estimate and using the support of $\chi'$
gives
\[
 \int W(\chi\eta)^2\le8\int W((\chi\eta)')^2
 \le16\int_{v\le2b} W(\eta')^2
       +16\|\chi'\|_\infty^2\int_{b\le v\le2b}W\eta^2 .
\]
All three integrals are in $\dd r$. Since $\chi=1$ on $v\le b$,
\begin{equation}\label{eq:gap-exterior}
 \int_{v\le b}W\eta^2\dd r\le C\left[
 \int_{v\le2b} W(\eta')^2\dd r+\int_{b\le v\le2b}W\eta^2\dd r\right].
\end{equation}

\noindent Step 3. We use the exterior estimate to absorb the negative part of the quadratic form.
Recall that $W=r^3v^2$. The coefficient of $\eta^2$ in the full
integral \eqref{eq:scaled-radial-form} is $(3/4)W yv^{q-1}B(v)$, where
\[
 B(v)=k(p-1)v^d+ek(p-1)\frac{v^d-1}{d}+e(k-q+1).
\]
For $v\ge b$, $v^d\ge b^D$ and $(v^d-1)/d\ge\log b$. Since
$k\ge k_0>0$ and $p-1,d,k$ are bounded above,
\[
 B(v)\ge k(p-1)b^D-Ce(k+1)\ge c k
\]
after decreasing $e_D$. Thus the coefficient of $\eta^2$ is
nonnegative on $v\ge b$. On $b\le v\le2b$ it is at least $cW$, since
$y\ge1/\sqrt2$ and $v^{q-1}\ge b^{3/2}$. For $v\le b$ we have
\[
 B(v)\ge-e(q-1)-ek(p-1)|\log v|,\qquad
 yv^{q-1}(1+|\log v|)\le C .
\]
Hence the negative part of the integral containing $\eta^2$ is at most
\[
 Ce\int_{v\le b}W\eta^2
 \le Ce\left(\int_{v\le2b}W(\eta')^2+
             \int_{b\le v\le2b}W\eta^2\right)
\]
by the exterior estimate~\eqref{eq:gap-exterior}. All integrals here are in $\dd r$.

On $v\le2b$ we have $W\le\sqrt2\,r^2v^2K$. On $b\le v\le2b$
the coefficient of $\eta^2$ is at least $cW$, with $c>0$ independent
of $e,d,k$. Thus the two positive integrals in
\eqref{eq:scaled-radial-form}, after subtracting the exterior error, give
\[
 \frac{\langle Lh,h\rangle}{4\pi}
 \ge (1-Ce)\int_0^\infty r^2v^2K(\eta')^2\dd r
      +(c-Ce)\int_{b\le v\le2b}W\eta^2\dd r .
\]
We have also dropped the nonnegative integral containing $\eta''$.
Choose $e_D$ so that $Ce_D\le\min(1,c)/2$. After decreasing $c$,
we obtain
\[
 \frac{\langle Lh,h\rangle}{4\pi}
 \ge \frac12\int_0^\infty r^2v^2K(\eta')^2\dd r
       +c\int_{b\le v\le2b}W\eta^2\dd r .
\]

To pass to radial $H^1$, approximate $h$ by smooth compactly supported
$h_n$ with $\int vh_n=0$, using the mass correction in the proof of
Proposition~\ref{prop:tangent-positive}. Then the forms converge and $\eta_n\to\eta$
locally uniformly on $(0,\infty)$. Along a subsequence $h_n\to h$ a.e.,
and
\[
 \tfrac12 r\eta_n'=h_n/v-(3/2-ry)\eta_n
 \longrightarrow h/v-(3/2-ry)\eta=\tfrac12 r\eta'\quad\text{a.e.}
\]
Fatou's lemma extends the two nonnegative integrals in the lower bound
to $h$. If the form is zero, $K>0$ implies $\eta'=0$ a.e. The primitive
$\eta$ is locally absolutely continuous, so it is constant. The integral
over $b<v<2b$ forces the constant to vanish, since $v$ decreases from
$v(0)>1$ to zero. The identity $h/v=(3/2-ry)\eta+r\eta'/2$ now gives $h=0$.

\end{proof}

\begin{proof}[\normalfont Proof of Theorem~\ref{thm:uniform-fold}]
Choose $e_D$ small enough for Lemma~\ref{lem:forward-small} and
Proposition~\ref{prop:forward-outer}. The lemma gives one nondegenerate
minimum for $0<k\le k_0$ and $M'(\lambda(k_0))>0$.
Suppose that $M'(\lambda)=0$ at an admissible frequency with
$k\ge k_0$. Put $w=\partial_\lambda u_\lambda$ in the original variables.
Then $\int u_\lambda w=0$, and the original linearized equation~\eqref{eq:linearized-frequency} reads
$L_{\rm orig}w=-u_\lambda$. Keep the two scale factors $s$ and
$\sqrt\lambda$ fixed at this frequency and set $h(x)=w(x/\sqrt\lambda)$.
Changing variables in $L_{\rm orig}w=-u_\lambda$ and in the mass integral gives
\[
 Lh=-(s/\lambda)v,\qquad
 \int vh=\frac{\lambda^{3/2}}s\int u_\lambda w=0.
\]
Since $v>0$, the equation for $Lh$ implies $h\ne0$. Taking the inner
product with $h$ gives
\[
 \langle Lh,h\rangle=-\frac{s}{\lambda}\int vh=0,
\]
contrary to \eqref{eq:fold-tangent-positive}. Thus $M'$ never vanishes
on the interval $\lambda(k_0)\le\lambda<\lambda_*$.
$M'$ is positive at $\lambda(k_0)$ by \eqref{eq:fold-small-sign}.
Continuity gives $M'>0$ on this entire interval.

By \eqref{eq:ends}, $M$ tends to infinity both at $0$ and at $\lambda_*$.
It decreases to $m=M(\lambda_c)$ on the first interval and increases
from $m$ on the second. Each $a>m$ is attained once on each side of
$\lambda_c$. The mass $m$ is attained only at $\lambda_c$, and a mass
below $m$ is never attained. Each frequency gives one positive solution
up to translations.
\end{proof}

\section{\texorpdfstring{Strict decrease for $1<p<q<5$ and $q\ge7/3$}{Strict decrease for 1 < p < q < 5 and q >= 7/3}}\label{sec:variational}
For $1<p<q<5$, there is a unique centered positive profile at every
$\lambda>0$, and the branch converges to the zero-frequency profile
as $\lambda\downarrow0$ \cite{FH,FHnew}. When $q\ge7/3$, we use the
mass-preserving dilation to prove strict decrease on the whole branch.

\begin{theorem}\label{thm:powers}
Let $1<p<q<5$ and $q\ge7/3$. Then $M'(\lambda)<0$ for every $\lambda>0$.
Put $\ell=m_Q$ when $q=7/3$ and $\ell=0$ when $q>7/3$.
If $p<7/3$, a mass admits a positive solution exactly when $\ell<a\le a_0$.
If $p\ge7/3$, the condition is $\ell<a<\infty$. In both cases the solution is unique up to translations.
\end{theorem}

We estimate $M'$ by testing the second variation of the action in
a direction tangent to its constraint. The test combines the frequency
derivative with a dilation that preserves the mass.
The linearization and endpoint estimates below also apply when $q<7/3$.
For these estimates take $1<p<q<5$, fix $\lambda>0$, and write $u=u_\lambda$.
Let $L=-\Delta+\lambda+pu^{p-1}-qu^{q-1}$ be the radial linearized operator.

\begin{lemma}\label{lem:rad}
The operator $L:H^2_{\rm rad}\to L^2_{\rm rad}$ is invertible.
On $1<p<q<5$, $\lambda>0$, the centered profile depends jointly
and continuously differentiably on $(p,q,\lambda)$ in $H^2_{\rm rad}$.
Its frequency derivative satisfies
\begin{equation}\label{eq:reverse-linearized}
 w=\partial_\lambda u,\qquad Lw=-u,\qquad
 L=-\Delta+\lambda+pu^{p-1}-qu^{q-1}.
\end{equation}
\end{lemma}
\begin{proof}[\normalfont\proofname]
The potential $pu^{p-1}-qu^{q-1}$ is bounded and tends to zero at infinity.
Multiplication by this potential is compact from radial $H^2$ to radial
$L^2$ by the local compactness and exterior estimate in Section~\ref{sec:ends}.
Since $-\Delta+\lambda$ is invertible for $\lambda>0$, it suffices by
the Fredholm alternative to check that the radial kernel of $L$ is zero.
This follows from \cite[Proposition~2.3]{FHnew}.

For fixed $p,q$, the $H^1$ frequency derivative and its equation are
given in \cite[Proposition~2.4 and (2.2)]{FHnew}. To vary the powers
as well, use the $C^1$ map from the proof of Lemma~\ref{lem:forward-branch},
\[
 (s,v)\longmapsto |v|^{s-1}v,\qquad
 (1,\infty)\times H^2_{\rm rad}\longrightarrow L^2_{\rm rad}.
\]
The implicit function theorem gives a local $C^1$ branch in
$(p,q,\lambda)$. Every nearby profile $v$ is nonzero by $H^2$
continuity. The same continuity gives $\|v_-\|_\infty^{q-1}<\lambda$.
By \eqref{eq:profile-negative-part},
\[
 \int|\nabla v_-|^2+
 (\lambda-\|v_-\|_\infty^{q-1})\int v_-^2\le0.
\]
Thus $v_-=0$. The strong maximum principle gives $v>0$, and uniqueness
at the given frequency gives $v=u_\lambda$. Differentiating the profile
equation now gives \eqref{eq:reverse-linearized}.
\end{proof}

\subsection{An upper bound for the mass derivative}\label{sec:reverse-slope}
Fix $\lambda>0$ and write $u=u_\lambda$. Let $\Phi$ be the action
and $I$ its derivative along $t\mapsto tv$ at $t=1$.
\begin{equation}\label{eq:action-constraint}
\begin{aligned}
 \Phi(v)&=\tfrac12\int(|\nabla v|^2+\lambda v^2)
        +\tfrac1{p+1}\int|v|^{p+1}-\tfrac1{q+1}\int|v|^{q+1},\\
 I(v)&=\left.\frac{d}{dt}\Phi(tv)\right|_{t=1}=\langle \Phi'(v),v\rangle
       =\int(|\nabla v|^2+\lambda v^2+|v|^{p+1}-|v|^{q+1}).
\end{aligned}
\end{equation}
The mass-preserving dilation gives an upper bound for $M'$ for every
$1<p<q<5$.

\begin{proposition}\label{prop:slope-bound}
Set $D=3u/2+x\cdot\nabla u$. Then $I'(u)D<0$ and
\begin{equation}\label{eq:generalbound}
 \tfrac12 M'(\lambda)\le
 \frac{M(\lambda)^2}{(I'(u)D)^2}\langle LD,D\rangle .
\end{equation}
The quadratic form on the right is
\begin{equation}\label{eq:curv}
 \langle LD,D\rangle=
 \frac{7-3q}{2}\int|\nabla u|^2
 -\frac{9(p-1)(q-p)}{4(p+1)}\int u^{p+1}.
\end{equation}
\end{proposition}
\begin{proof}[\normalfont\proofname]
By \cite[Theorem~2.1]{FH}, $u$ minimizes $\Phi$ among nonzero
functions satisfying $I(v)=0$. Since $\Phi'(u)=0$, differentiating
$I(v)=\langle \Phi'(v),v\rangle$ gives $I'(u)h=\langle Lu,h\rangle$.
Eliminating $\int u^{q+1}$ by $I(u)=0$ gives
\[
 I'(u)u=-(q-1)\int(|\nabla u|^2+\lambda u^2)-(q-p)\int u^{p+1}<0.
\]
If $I'(u)h=0$, the implicit function theorem and $I'(u)u<0$ give
\[
 I(u+th+\sigma(t)u)=0,\qquad
 \sigma(0)=0,\qquad \sigma'(0)=-\frac{I'(u)h}{I'(u)u}=0 .
\]
The action has a minimum at $t=0$ along this curve. Since
$\sigma'(0)=0$ and $\Phi'(u)=0$, its second derivative at zero is
$\Phi''(u)[h,h]=\langle Lh,h\rangle$. Hence
\begin{equation}\label{eq:tangent}
 \langle Lh,h\rangle\ge0\qquad\text{whenever }I'(u)h=\langle Lu,h\rangle=0 .
\end{equation}

We construct a direction to which the tangent space inequality~\eqref{eq:tangent} applies. Start with
$u_s(x)=s^{3/2}u(sx)$, whose mass is independent of $s$. Thus
\[
 D=\left.\partial_s u_s\right|_{s=1}=\tfrac32 u+x\cdot\nabla u,
 \qquad \int uD=0.
\]
Changing variables in the action gives
\[
 \Phi(u_s)=\tfrac{s^2}{2}\int|\nabla u|^2+\tfrac{\lambda}{2}M
       +\frac{s^{3(p-1)/2}}{p+1}\int u^{p+1}
       -\frac{s^{3(q-1)/2}}{q+1}\int u^{q+1}.
\]
Since $\Phi'(u)=0$, setting the first derivative in $s$ to zero at $s=1$ gives
\begin{equation}\label{eq:mass-dilation}
 \frac{3(q-1)}{2(q+1)}\int u^{q+1}
 =\int|\nabla u|^2+\frac{3(p-1)}{2(p+1)}\int u^{p+1}.
\end{equation}
Differentiate $\Phi(u_s)$ once more at $s=1$ and use
\eqref{eq:mass-dilation} to eliminate $\int u^{q+1}$. This gives
\[
 \left.\frac{d^2}{ds^2}\Phi(u_s)\right|_{s=1}
 =\frac{7-3q}{2}\int |\nabla u|^2
 -\frac{9(p-1)(q-p)}{4(p+1)}\int u^{p+1}.
\]
Since $\Phi'(u)=0$, this second derivative is $\langle LD,D\rangle$.
Differentiating the constraint $I(u_s)$ at $s=1$ gives
\begin{equation}\label{eq:cross}
\begin{aligned}
 I'(u)D
 &=2\int|\nabla u|^2+\frac{3(p-1)}{2}\int u^{p+1}-\frac{3(q-1)}{2}\int u^{q+1}\\
 &=-(q-1)\int|\nabla u|^2-\frac{3(p-1)(q-p)}{2(p+1)}\int u^{p+1}<0.
\end{aligned}
\end{equation}

For $w=\partial_\lambda u_\lambda$, self-adjointness and the linearized equation~\eqref{eq:reverse-linearized} give
\[
 \begin{aligned}
 I'(u)w&=\langle Lu,w\rangle=-\int u^2=-M,\\
 \langle Lw,D\rangle&=-\int uD=0,\qquad
 \langle Lw,w\rangle=-\int uw=-M'(\lambda)/2 .
 \end{aligned}
\]
To make $w$ tangent to $I=0$, add a multiple of $D$. Since $I'(u)D<0$,
the choice $h=w+\frac{M}{I'(u)D}D$ gives $I'(u)h=-M+M=0$.
The mixed term in $\langle Lh,h\rangle$ vanishes. Thus the tangent space inequality~\eqref{eq:tangent} gives
\[
 0\le\langle Lh,h\rangle=-\tfrac12M'(\lambda)
             +\frac{M^2}{(I'(u)D)^2}\langle LD,D\rangle .
\]
Rearranging gives \eqref{eq:generalbound}. The denominator is nonzero by
\eqref{eq:cross}.
\end{proof}

\subsection{\texorpdfstring{Endpoint masses for $p<q<5$}{Endpoint masses for the higher focusing power}}\label{sec:reverseends}
At large frequency a rescaling gives a regular equation near $Q_q$.
Use the scaling
\[
 t=\lambda^{-(q-p)/(q-1)},\qquad
 u_\lambda(x)=\lambda^{1/(q-1)}v_t(\sqrt\lambda\,x).
\]
The profile equation in \eqref{eq:main} and the mass become
\begin{equation}\label{eq:scaled}
 -\Delta v_t+v_t+t v_t^p=v_t^q,\qquad
 M(\lambda)=\lambda^{2/(q-1)-3/2}\int v_t^2.
\end{equation}
As $\lambda\to\infty$, $t\to0$, so we study the rescaled equation
near its solution $Q_q$ at $t=0$. We use a weighted norm so the
profile stays positive even at large radii.

\begin{lemma}\label{lem:weighted}
Fix $1<b<5$ and $d>1$. Let $Q_b$ be the centered positive solution of $-\Delta Q_b+Q_b=Q_b^b$. Near $(q,p,t)=(b,d,0)$ the rescaled equation~\eqref{eq:scaled}, written as
\[
 -\Delta v+v=v^q-tv^p,
\]
has a positive radial branch through $Q_b$ such that $v$ and $\int v^2$ are $C^2$ in $(q,p,t)$. For $t>0$ and
$p\ne q$, whenever rescaling produces an allowed frequency, this branch gives the centered profile.
\end{lemma}
\begin{proof}[\normalfont\proofname]
For $p<2$, the second derivative of $v^p$ contains $v^{p-2}$, which
grows as $v$ tends to zero. We choose a norm that controls each
perturbation by the decay of the profile. Set $h(r)=e^{-r}/(1+r)$ and
choose $\tau>0$ with $1+2\tau<\min(b,d)$. For continuous radial
functions define
\[
 \|v\|_X=\sup_{r\ge0}|v(r)|/h(r),\qquad
 \|f\|_Y=\sup_{r\ge0}e^{(1+\tau)r}|f(r)|.
\]
In each case take the space of functions with finite norm.

We first show that $G=(-\Delta+1)^{-1}$ maps $Y$ compactly into $X$.
For $f\in Y$, the radial solution $v=Gf$ satisfies
\[
 -(rv)''+rv=rf,\qquad (rv)(0)=0,\qquad
 \lim_{r\to\infty}rv(r)=0 .
\]
The Green function on the half-line gives
\begin{equation}\label{eq:radial-inverse}
 Gf(r)=\frac{e^{-r}}{r}\int_0^r s\sinh(s)f(s)\dd s
       +\frac{\sinh(r)}{r}\int_r^\infty s e^{-s}f(s)\dd s .
\end{equation}
At $r=0$ use the continuous limit. To prove compactness of $G:Y\to X$,
take $\|f\|_Y\le1$. With $h=e^{-r}/(1+r)$, the radial inverse formula~\eqref{eq:radial-inverse} gives for $r\ge1$
\[
 \left|\frac{Gf(r)}{h(r)}-\int_0^\infty s\sinh(s)f(s)\dd s\right|
 \le \frac C r+C(1+r)^2e^{-\tau r}.
\]
Indeed the first prefactor after division by $h$ is $1+1/r$, the omitted
tail is bounded by $C\int_r^\infty se^{-\tau s}\dd s$, and the second
term after division is bounded by $C(1+r)^2e^{-\tau r}$. On $[0,R]$,
the integral formula~\eqref{eq:radial-inverse} and $(r^2(Gf)')'=r^2(Gf-f)$ yield
\[
 \sup_{[0,R]}|Gf|+\sup_{[0,R]}|(Gf)'|\le C_R .
\]
Arzel\`a--Ascoli gives a subsequence converging locally uniformly.
The numbers $\int_0^\infty s\sinh(s)f(s)\dd s$ are bounded, so they
also converge along a subsequence. The estimate for $r\ge R$ makes
$Gf/h$ uniformly Cauchy outside $[0,R]$ once $R$ is large. This proves
compactness in $X$.

We next compare $Q_b$ with $h$ to control the derivatives of $v^s$.
For $0<\mu<1$, choose $R$ so large that $1-Q_b^{b-1}\ge\mu^2$
on $r\ge R$. Comparison with a multiple of $e^{-\mu r}/r$ gives
$Q_b(r)\le C e^{-\mu r}/r$ there. Choose $\mu$ close enough to one
that $b\mu>1+\tau$. Then $Q_b^b\in Y$ and \eqref{eq:radial-inverse} gives
\[
 Q_b=G(Q_b^b),\qquad
 \lim_{r\to\infty}\frac{Q_b(r)}{h(r)}
 =\int_0^\infty s\sinh(s)Q_b(s)^b\,ds>0 .
\]
Positivity on bounded intervals now gives $ch\le Q_b\le Ch$ on
$[0,\infty)$. The same bounds, with adjusted constants, hold for $v$
in a small $X$-ball about $Q_b$.

For $ch\le v\le Ch$ and $s$ near $b$ or $d$ we have
\[
 |v^{s-2}z_1z_2|\le C h^s\|z_1\|_X\|z_2\|_X,\qquad
 |\log v|\le C(1+r).
\]
Each exponent derivative contributes one logarithm. Shrink the exponent
neighborhood so that $s\ge1+2\tau$. Then
\[
 \sup_{r\ge0}e^{(1+\tau)r}h(r)^s(1+r)^2<\infty .
\]
Taylor's formula on $ch(r)\le v(r)\le Ch(r)$ now gives two continuous
derivatives of $(v,s)\mapsto v^s$ into $Y$.

We now apply the implicit function theorem to the profile equation.
At $(q,p,t,v)=(b,d,0,Q_b)$, the equation and its derivative in $v$ are
\[
 v-G(v^q-tv^p)=0,\qquad z\longmapsto z-G(bQ_b^{b-1}z).
\]
Multiplication by $bQ_b^{b-1}$ maps $X$ continuously into $Y$.
Since $G:Y\to X$ is compact, the derivative is the identity minus a
compact operator. If $z=G(bQ_b^{b-1}z)$, then $bQ_b^{b-1}z\in Y\subset L^2$ and
\[
 z\in H^2_{\rm rad},\qquad (-\Delta+1-bQ_b^{b-1})z=0.
\]
The radial kernel is zero \cite[Section~4.1]{LRN}, so $z=0$. The
Fredholm alternative gives a bounded inverse. The implicit function
theorem now gives a $C^2$ branch in $X$. Each profile is positive
because $v\ge ch$. Also,
\[
 \|v\|_2\le\|h\|_2\|v\|_X,
\]
so the mass is $C^2$. At each admissible frequency the rescaled
solution is the unique centered positive profile.
\end{proof}

The behavior at the two ends determines which mass levels occur.

\begin{proposition}\label{prop:reverse-ends}
For $1<p<q<5$ the endpoint masses are
\[
 \lim_{\lambda\downarrow0}M(\lambda)=
 \begin{cases}a_0,&p<7/3,\\ +\infty,&p\ge7/3,\end{cases}
 \qquad
 \lim_{\lambda\to\infty}M(\lambda)=
 \begin{cases}+\infty,&q<7/3,\\m_Q,&q=7/3,\\0,&q>7/3.\end{cases}
\]
\end{proposition}
\begin{proof}[\normalfont\proofname]
As $\lambda\downarrow0$, $u_\lambda\to u_0$ in
$\dot H^1\cap L^{p+1}$ by \cite[Corollary~1.4 and Theorem~1.6]{FH}.
When $p<7/3$ the convergence holds in $H^1$, so $M(\lambda)\to a_0$.
When $p\ge7/3$, $u_0\notin L^2$. On each ball $B_R$, H\"older's
inequality gives
\[
 \|u_\lambda-u_0\|_{L^2(B_R)}
 \le |B_R|^{1/2-1/(p+1)}\|u_\lambda-u_0\|_{p+1}\longrightarrow0 .
\]
It follows that
\[
 \liminf_{\lambda\downarrow0} M(\lambda)\ge\int_{B_R}u_0^2
 \qquad\text{for every }R>0 .
\]
Since $u_0\notin L^2$, letting $R\to\infty$ gives $M(\lambda)\to\infty$.

As $\lambda\to\infty$, the parameter $t$ in \eqref{eq:scaled} tends
to zero. Lemma~\ref{lem:weighted} gives $\int v_t^2\to\int Q_q^2>0$.
The exponent $2/(q-1)-3/2$ in \eqref{eq:scaled} is positive, zero or
negative according as $q<7/3$, $q=7/3$ or $q>7/3$. This proves the
three limits at infinity.
\end{proof}

\begin{proof}[\normalfont Proof of Theorem~\ref{thm:powers}]
For $q\ge7/3$, the gradient term in \eqref{eq:curv} is nonpositive.
The other term is strictly negative because $(p-1)(q-p)>0$, including
when $q=7/3$. Hence $\langle LD,D\rangle<0$, and \eqref{eq:generalbound}
gives $M'(\lambda)<0$ at every positive frequency.

With $\ell=m_Q$ for $q=7/3$ and $\ell=0$ for $q>7/3$, the endpoint
limits in Proposition~\ref{prop:reverse-ends} and strict decrease give
\[
 M((0,\infty))=
 \begin{cases}
 (\ell,a_0),&p<7/3,\\
 (\ell,\infty),&p\ge7/3 .
 \end{cases}
\]
Each mass in this range is attained at exactly one positive frequency.
For $p<7/3$, the zero-frequency solution supplies the upper endpoint
$a=a_0$. For $p\ge7/3$, that profile has infinite mass.
Fixed-frequency uniqueness gives the claimed counts.
\end{proof}

\section{\texorpdfstring{One minimum for $p<q$ and $q$ just below $7/3$}{One minimum for p < q and q just below 7/3}}\label{sec:strip}
We next follow the minimum that appears when the higher focusing power
$q$ decreases below $7/3$.

\begin{theorem}\label{thm:reverse-uniform}
For each $1<b<7/3$ there is $\varepsilon_b>0$ such that if
$b\le p<q<7/3$ and $7/3-q<\varepsilon_b$, then $M'$ has a unique zero
$\lambda_c>0$. We have
\[
 M'<0\quad(0<\lambda<\lambda_c),\qquad M''(\lambda_c)>0,\qquad
 M'>0\quad(\lambda>\lambda_c).
\]
Put $m=M(\lambda_c)$. The solution counts on
$(0,m),\{m\},(m,a_0],(a_0,\infty)$ are $0,1,2,1$.
\end{theorem}

For $1<b<7/3$ and $b\le p<q<7/3$, write $e=7/3-q$ and $d=q-p$.
In this range the rescaled parameter $k$ decreases as the frequency
increases. We locate the minimum at large frequencies and then use
Proposition~\ref{prop:slope-bound} at smaller frequencies.

\begin{lemma}\label{lem:reverse-large}
For each $1<b<7/3$ there are $k_0,\varepsilon_b>0$ such that,
if $b\le p<q<7/3$ and $0<e<\varepsilon_b$, there is exactly one
stationary frequency $\lambda_c$ in $[\lambda(k_0),\infty)$,
where $\lambda(k)$ is given by \eqref{eq:reverse-k}.
We have $M'<0$ on $[\lambda(k_0),\lambda_c)$, $M'>0$ on
$(\lambda_c,\infty)$, and $M''(\lambda_c)>0$.
\end{lemma}
\begin{proof}[\normalfont\proofname]
For each $\lambda>0$ there is exactly one $s>1$ satisfying
$s^{q-1}-s^{p-1}=\lambda$. Indeed the left side increases from zero
to infinity on $s>1$. Set $k=d/(s^d-1)$. Then
$s^d=1+d/k$, and $\lambda=s^{q-1}d/(k+d)$. Hence
\begin{equation}\label{eq:reverse-k}
 \lambda(k)=\frac d{k+d}\left(1+\frac dk\right)^{(q-1)/d},\qquad
 \frac{\dd\log\lambda}{\dd k}=-\frac{q-1+k}{k(k+d)}<0 .
\end{equation}
As $k$ increases from zero to infinity, $\lambda$ decreases from
infinity to zero.

Write $u_\lambda(x)=s\,v(\sqrt\lambda\,x)$. Since
$s^{q-1}/\lambda=1+k/d$ and $s^{p-1}/\lambda=k/d$, we have
\begin{equation}\label{eq:reverse-profile}
 -\Delta v+v=v^q+k\frac{v^q-v^{q-d}}d.
\end{equation}
Changing variables in the integral of $u_\lambda^2$ gives
\[
 M(\lambda(k))=\left(1+\frac kd\right)^{3/2}
                \left(1+\frac dk\right)^{3e/(2d)}H(e,d,k),
 \qquad H(e,d,k)=\int v^2 .
\]
The logarithmic derivative of the explicit factor is
$3(k-e)/(2k(k+d))$. Thus
\begin{equation}\label{eq:reverse-mass-derivative}
 \frac{\dd M(\lambda(k))}{\dd k}
 =\frac{M(\lambda(k))}{k(k+d)H}F(e,d,k),\qquad
 F(e,d,k)=\tfrac32(k-e)H+k(k+d)H_k .
\end{equation}
Since $\lambda'(k)<0$ in \eqref{eq:reverse-k}, $M'(\lambda)$ has the opposite sign to $F$.

We need $H$ and its first two derivatives uniformly for
$0\le d\le7/3-b$. At $e=k=0$ the profile is $Q=Q_{7/3}$ for every $d$.
For positive $v$ we can write
\[
 \frac{v^q-v^{q-d}}d=\int_0^1 v^{q-td}\log v\dd t.
\]
The integral also defines the quotient at $d=0$. Take
$\tau=(b-1)/4$. Since $q=7/3-e$ and $0\le d\le7/3-b$,
\[
 q-td\ge b-|e|\ge1+3\tau\qquad(0\le t\le1,\ |e|\le\tau).
\] In the weighted ball about $Q$ in
Lemma~\ref{lem:weighted}, $v$ is comparable to $e^{-r}/(1+r)$ and
\[
 \sup_{r\ge0}e^{(1+\tau)r}v(r)^{q-td}
                 (1+|\log v(r)|^3)\le C.
\]
Write $h=e^{-r}/(1+r)$ and $a=q-td$. Each perturbation $z\in X$
satisfies $|z|\le h\|z\|_X$. Since $v\ge ch$, we have
\[
 \left|\frac{z}{v}\right|\le C\|z\|_X .
\]
For $i+j\le2$, differentiating $v^a\log v$ $i$ times in $v$ and
$j$ times in $a$ therefore gives
\[
 \big|D_v^i\partial_a^j(v^a\log v)[z_1,\ldots,z_i]\big|
 \le C v^a(1+|\log v|^3)\prod_{\ell=1}^i\|z_\ell\|_X .
\]
For $i=0$ the product equals one. The preceding weighted bound and
the power estimates in Lemma~\ref{lem:weighted} give two continuous
derivatives into $Y$, uniformly for $0\le d\le7/3-b$ and $0\le t\le1$.

At $e=k=0$, the radial linearization is $L_Q=-\Delta+1-(7/3)Q^{4/3}$.
The radial inverse used in Lemma~\ref{lem:weighted} gives a local $C^2$
branch by the implicit function theorem. All local branches equal $Q$
at $e=k=0$. The finite-cover argument in Lemma~\ref{lem:forward-profile}
gives a common neighborhood in $(e,k)$. Each profile satisfies $v\ge ch$
and therefore rescales to the unique centered positive solution.

To count the zeros of $F$ for small $k$, we compute $F_k(0,d,0)$.
For $e=k=0$ and $d>0$, differentiate the rescaled equation~\eqref{eq:reverse-profile} in $k$. Recall from
\eqref{eq:critical-dilation} that $Z=3Q/4+x\cdot\nabla Q/2$ and $L_QZ=-Q$.
Then
\begin{align*}
 dL_Qv_k&=Q^{7/3}-Q^{7/3-d},\\
 dH_k&=-2d\int ZL_Qv_k=2\int Z(Q^{7/3-d}-Q^{7/3}),\\
 F_k(0,d,0)&=\tfrac32 m_Q+dH_k=2\int ZQ^{7/3-d}
 =\frac{3(4/3-d)}{2(10/3-d)}\int Q^{10/3-d}.
\end{align*}
Here $2\int ZQ^{7/3}=3m_Q/2$ by \eqref{eq:critical-moment}. The last
equality follows from \eqref{eq:dilation-moment} with $a=7/3-d$.

Since $d\le7/3-b$, we have $4/3-d\ge b-1>0$. The formula extends continuously to $d=0$,
so the derivative has a positive minimum on the entire closed
$d$-interval. Choose $k_0>0$ and then an upper bound for $e$ so that
$F_k\ge c>0$ on $0\le k\le k_0$ uniformly in $d$. For $e>0$ we have $F(e,d,0)=-3eH(e,d,0)/2<0$. Integrating
$F_k\ge c$ from $0$ to $k_0$ gives
\[
 F(e,d,0)<0,\qquad F(e,d,k_0)\ge F(e,d,0)+ck_0>0
\]
once $\tfrac32 e\sup_d H(e,d,0)<ck_0$. The bound on $H$ allows
one choice of the upper bound for $e$ for all $d$. Thus $F$ has exactly
one zero $k_c$ in $(0,k_0)$, with $F_k(e,d,k_c)>0$. Set
$\lambda_c=\lambda(k_c)$. Differentiating \eqref{eq:reverse-mass-derivative}
at $k_c$ eliminates the derivative of its prefactor because $F(e,d,k_c)=0$.
The chain rule then gives
\[
 M''(\lambda_c)=
 \frac{1}{\lambda'(k_c)^2}
 \frac{M(\lambda_c)F_k(e,d,k_c)}{k_c(k_c+d)H(e,d,k_c)}>0 .
\]
Since $\lambda'(k)<0$, the signs of $F$ on either side of $k_c$ give
\[
 M'(\lambda)>0\quad(\lambda>\lambda_c),\qquad
 M'(\lambda)<0\quad(\lambda(k_0)\le\lambda<\lambda_c).
\]
\end{proof}

At smaller frequencies, keeping the scale at $k_0$ fixed gives one
energy estimate for the whole interval.

\begin{lemma}\label{lem:reverse-low}
Under the assumptions of Lemma~\ref{lem:reverse-large}, after decreasing
$\varepsilon_b$ if necessary, $M'(\lambda)<0$ for
$0<\lambda\le\lambda(k_0)$.
\end{lemma}
\begin{proof}[\normalfont\proofname]
Put $\lambda_0=\lambda(k_0)$ and keep the amplitude
$s_0=(1+d/k_0)^{1/d}$ and the spatial scale
$\sqrt{\lambda_0}$ fixed as $\lambda$ varies. With
$u_\lambda(x)=s_0v(\sqrt{\lambda_0}\,x)$ the profile equation in \eqref{eq:main} becomes
\begin{equation}\label{eq:fixed-scale}
 -\Delta v+\rho v=(1+k_0/d)v^q-(k_0/d)v^p,\qquad
 \rho=\lambda/\lambda_0\in(0,1].
\end{equation}
We need to bound $\int|\nabla v|^2$ by a constant times
$\int v^{p+1}$, independently of $\rho,e,d$. First we bound the
gradient integral alone. In these variables the action in
\eqref{eq:action-constraint} becomes
\[
 \Phi_{\rm sc}(v)=\frac{\sqrt{\lambda_0}}{s_0^2}
       \Phi\bigl(s_0 v(\sqrt{\lambda_0}\,\cdot)\bigr).
\]
For fixed nonzero $\psi\in C_c^\infty$ with $\psi\ge0$, this gives
\begin{align*}
 \Phi_{\rm sc}(t\psi)
 &=\frac{t^2}{2}\int(|\nabla\psi|^2+\rho\psi^2)
 -\frac{t^{q+1}}{q+1}\int\psi^{q+1}\\
 &\hspace{10mm}
 +\frac{k_0}{d}\int\left(\frac{(t\psi)^{p+1}}{p+1}
                           -\frac{(t\psi)^{q+1}}{q+1}\right).
\end{align*}
The difference in the last integral is controlled uniformly as $d=q-p$
approaches zero. Indeed
\[
 \frac{\partial}{\partial a}\frac{z^{a+1}}{a+1}
 =\frac{z^{a+1}}{a+1}\left(\log z-\frac1{a+1}\right).
\]
This derivative is positive for $z\ge2$ and uniformly bounded for
$0\le z\le2$, $b\le a\le7/3$. The mean value theorem yields
\[
 \frac1d\left(\frac{z^{p+1}}{p+1}-\frac{z^{q+1}}{q+1}\right)
 \le C\mathbf1_{\{0<z\le2\}} .
\]
On the ray $t\psi$, the constraint function
$t\,\frac{d}{dt}\Phi_{\rm sc}(t\psi)$ is positive for small $t>0$
because its quadratic term is positive. It is negative for large
$t$ because $q+1>p+1>2$ and the coefficient of $t^{q+1}$ is negative.
Thus at least one $t\psi$ belongs to the constraint. The fixed change of
variables multiplies every action value by the same positive constant,
so $v$ minimizes $\Phi_{\rm sc}$ on the rescaled constraint.

Dividing the Pohozaev identity for \eqref{eq:fixed-scale} by three gives
\[
 \frac16\int|\nabla v|^2+\frac{\rho}{2}\int v^2
 +\frac{k_0}{d(p+1)}\int v^{p+1}
 -\frac{1+k_0/d}{q+1}\int v^{q+1}=0 .
\]
The Pohozaev identity gives $\Phi_{\rm sc}(v)=\frac13\int|\nabla v|^2$.
Since $v$ minimizes the rescaled action on the constraint and the ray
$t\psi$ meets that constraint,
\[
 \begin{aligned}
 \frac13\int|\nabla v|^2
 &=\Phi_{\rm sc}(v)\le\sup_{t>0}\Phi_{\rm sc}(t\psi)\\
 &\le\sup_{t\ge0}\left[C(1+t^2)
       -\frac{t^{q+1}}{q+1}\int\psi^{q+1}\right]\le C.
 \end{aligned}
\]
The constants are uniform since $q\ge b>1$ and $\psi$ is fixed.

Testing the fixed-scale equation~\eqref{eq:fixed-scale} with $v$ gives
\[
 \int|\nabla v|^2+\rho\int v^2
 =\int v^{p+1}\left(v^d+k_0\frac{v^d-1}{d}\right).
\]
For every $\eta>0$ and all $z>0$,
\[
 z^{p+1}\left(z^d+k_0\frac{z^d-1}{d}\right)
 \le C_\eta z^{p+1}+\eta z^6 .
\]
For $0<z\le1$ the bracket is at most one. For $z\ge1$,
\[
 \frac{z^d-1}{d}=\int_0^1 z^{td}\log z\dd t\le z^d\log z,\qquad
 \frac{z^{q+1}(1+k_0\log z)}{z^6}
 \le z^{-8/3}(1+k_0\log z)\longrightarrow0 .
\]
Thus $\eta z^6$ suffices for all large $z$, uniformly in $e,d$. On the
remaining bounded interval the bracket is uniformly bounded by the
integral formula. Testing with $v$ and using Sobolev's inequality gives
\[
 \int|\nabla v|^2\le C_\eta\int v^{p+1}
            +\eta C_{\rm Sob}\left(\int|\nabla v|^2\right)^3 .
\]
The uniform gradient bound allows us to choose $\eta>0$ such that
\[
 \eta C_{\rm Sob}\left(\sup\int|\nabla v|^2\right)^2\le\tfrac12 .
\]
Moving the last term to the left gives
\[
 \tfrac12\int|\nabla v|^2\le C_\eta\int v^{p+1}.
\]

We now use \eqref{eq:curv} at the original scale, with $u=u_\lambda$
and $D=3u/2+x\cdot\nabla u$.
Since $u(x)=s_0v(\sqrt{\lambda_0}\,x)$,
\[
 \int|\nabla u|^2=s_0^2\lambda_0^{-1/2}\int|\nabla v|^2,\qquad
 \int u^{p+1}=s_0^{p+1}\lambda_0^{-3/2}\int v^{p+1}.
\]
By \eqref{eq:reverse-k}, $(q-p)s_0^{p-1}/\lambda_0=k_0$.
Substitution in \eqref{eq:curv} gives
\[
 \frac{\langle LD,D\rangle}{s_0^2\lambda_0^{-1/2}}
 =\frac{3e}{2}\int|\nabla v|^2
  -\frac{9k_0(p-1)}{4(p+1)}\int v^{p+1}
 \le \left(3C_\eta e-\frac{9k_0(b-1)}{4(b+1)}\right)\int v^{p+1}.
\]
The inequality uses $\int|\nabla v|^2\le2C_\eta\int v^{p+1}$ and the increase of
$(p-1)/(p+1)$ for $p\ge b>1$.
The bracket is negative for all sufficiently small $e>0$, uniformly in
$d$ and $0<\lambda\le\lambda_0$. Now \eqref{eq:generalbound} gives
$M'(\lambda)<0$ on this whole frequency interval.

\end{proof}

\begin{proof}[\normalfont Proof of Theorem~\ref{thm:reverse-uniform}]
Lemmas~\ref{lem:reverse-large} and~\ref{lem:reverse-low} cover
$(0,\infty)$ and give $M''(\lambda_c)>0$ and
\[
 M'<0\quad(0<\lambda<\lambda_c),\qquad
 M'>0\quad(\lambda>\lambda_c).
\]
Here $p<q<7/3$, so $M(0)=a_0$ and $M(+\infty)=+\infty$.
Including the zero-frequency profile on the decreasing branch gives
\[
 M([0,\lambda_c))=(m,a_0],\qquad M((\lambda_c,\infty))=(m,\infty).
\]
Every $m<a\le a_0$ occurs on both branches. In particular,
$a=a_0$ occurs at zero and at one positive frequency above $\lambda_c$.
For $a>a_0$ only the increasing branch contributes. The mass $m$ occurs
only at $\lambda_c$, and no frequency gives a mass smaller than $m$.

\end{proof}

With $p$ fixed, we can also compute how fast $\lambda_c$ tends to
infinity and how $M(\lambda_c)$ approaches $m_Q$ as $q\uparrow7/3$.

\begin{proposition}\label{prop:minimum-at-infinity}
Fix $1<p<c=7/3$. As $q\uparrow c$, write $e=c-q$ and let
$\lambda_c$ denote the unique minimizing frequency from
Theorem~\ref{thm:reverse-uniform}. Then
\begin{align*}
 \lambda_c&\sim
 \left[\frac{(c-p)(p-1)\int Q^{p+1}}
 {(p+1)m_Q\,e}\right]^{(c-1)/(c-p)},\\
 M(\lambda_c)&=m_Q+\frac{3m_Q}{2(c-p)}\,e\log(1/e)+O(e).
\end{align*}
The constant in $O(e)$ may depend on $p$.
\end{proposition}
\begin{proof}[\normalfont\proofname]
For fixed $p$, set $H(q,t)=\int v_t^2$ in \eqref{eq:scaled}.
Lemma~\ref{lem:weighted} gives $H\in C^2$ near $(c,0)$ and
$H(c,0)=m_Q$. To compute $H_t(c,0)$, let
$L_Q=-\Delta+1-cQ^{c-1}$ and $z=\partial_t v_t|_{(c,0)}$.
Differentiation of the rescaled equation yields $L_Qz=-Q^p$.
Use $Z=3Q/4+x\cdot\nabla Q/2$ from \eqref{eq:critical-dilation}.
Since $L_Q Z=-Q$ and $L_Q z=-Q^p$, self-adjointness and
\eqref{eq:dilation-moment} with $a=p$ give
\[
 H_t(c,0)=2\int Qz=-2\int Z L_Q z
 =2\int ZQ^p=\frac{3(p-1)}{2(p+1)}\int Q^{p+1}>0 .
\]

With $t=\lambda^{-(q-p)/(q-1)}$ the mass is
$t^{-\sigma}H(q,t)$, where $\sigma=3e/[2(q-p)]$.
Stationarity is equivalent to $tH_t(q,t)-\sigma H(q,t)=0$.
The left side vanishes at $(q,t)=(c,0)$ and its derivative in $t$
there is $H_t(c,0)>0$. The implicit function theorem gives a
root tending to zero. Taylor's formula gives
\[
 t=\sigma\,\frac{m_Q}{H_t(c,0)}+O(e^2)>0 .
\]
Indeed, differentiability gives $t=O(e)$. Since $H$ is $C^2$,
\[
 H_t(q,t)=H_t(c,0)+O(e),\qquad H(q,t)=m_Q+O(e).
\]
By Theorem~\ref{thm:reverse-uniform}, this root gives the global
minimum. Substituting $t$ in $\lambda_c=t^{-(q-1)/(q-p)}$ gives the
stated asymptotic formula, since
\[
 \left(\frac{q-1}{q-p}-\frac{c-1}{c-p}\right)\log t
 =O(e|\log e|)\longrightarrow0.
\]
Finally,
\[
 t^{-\sigma}=1+\frac{3e}{2(c-p)}\log(1/e)+O(e).
\]
Multiplying by $H(q,t)=m_Q+O(e)$ proves the mass formula.
\end{proof}

\section{\texorpdfstring{The zero-frequency transition near $(p_\dagger,2p_\dagger-1)$}{The zero-frequency transition near (p-dagger, 2p-dagger - 1)}}\label{sec:birth-main}
For $1<p<q<7/3$, the mass starts at $a_0=M_{p,q}(0)$ and tends to
infinity as $\lambda\to\infty$. Near a particular point on $q=2p-1$,
the sign of the initial slope determines the shape of the entire curve.

\begin{theorem}\label{thm:global-birth}
Let $p_\dagger$ be the unique root in $(1,9/5)$ of
\[
 15p^3-38p^2-p+36=0 .
\]
Write $h=q-2p+1$. In a neighborhood of $(p_\dagger,2p_\dagger-1)$,
the initial slope $M_{p,q}'(0+)$ vanishes on a single $C^1$ curve
\[
 p=P(h),\qquad q=2P(h)-1+h,\qquad P(0)=p_\dagger .
\]
Count solutions up to translations, including the zero-frequency
solution at its mass $a_0=M_{p,q}(0)$.
\begin{enumerate}
\item If $p\le P(h)$, then $M'_{p,q}(\lambda)>0$ for every $\lambda>0$.
The initial slope is positive for $p<P(h)$ and zero for $p=P(h)$.
There is exactly one solution for $a\ge a_0$ and none for $0<a<a_0$.
\item If $p>P(h)$, then $M'_{p,q}(0+)<0$. The mass has exactly one
stationary frequency $\lambda_{\min}>0$. It decreases strictly up to
$\lambda_{\min}$ and increases strictly for $\lambda>\lambda_{\min}$, with
$M''_{p,q}(\lambda_{\min})>0$. Put $m=M_{p,q}(\lambda_{\min})<a_0$.
The number of solutions is
\[
\begin{array}{c|cccc}
 a &(0,m)&\{m\}&(m,a_0]&(a_0,\infty)\\ \hline
 \text{number of solutions}&0&1&2&1 .
\end{array}
\]
At $a=a_0$ the two frequencies are $0$ and $\lambda_+>\lambda_{\min}$.
As $(p,q)\to(p_\dagger,2p_\dagger-1)$ with $p>P(h)$,
the minimizing frequency $\lambda_{\min}$ tends to zero.
\end{enumerate}

On the whole line segment $q=2p-1$, $1<p<5/3$, one has
\[
 M'_{p,\,2p-1}(\lambda)>0\ \text{for every }\lambda>0
 \quad\Longleftrightarrow\quad p\le p_\dagger .
\]
For $p_\dagger<p<5/3$, the initial slope is negative and $M$ attains
a global minimum below $a_0$ at a positive frequency.
\end{theorem}

On $q=2p-1$ with $p_\dagger<p<5/3$, there is exactly one positive
stationary frequency whenever $p$ is close to either endpoint. Part~(2) of Theorem~\ref{thm:global-birth} gives
uniqueness near $p_\dagger$, and Theorem~\ref{thm:reverse-uniform} gives
uniqueness near $5/3$. For intermediate $p$, the number of stationary frequencies remains open.

On $q=2p-1$ the explicit zero-frequency profile gives a formula for
$M'(0+)$. Its sign changes at $p_\dagger$. We prove that $M''$ stays
positive at small positive frequencies near $(p_\dagger,2p_\dagger-1)$.
This determines the number of small stationary frequencies. To cover
the remaining frequencies we prove strict increase on the line up to
$p=p_\dagger$ and an upper bound for every frequency with $M'\le0$.

\begin{figure}[!htbp]
\centering
\iffigurepdfs
  \includegraphics{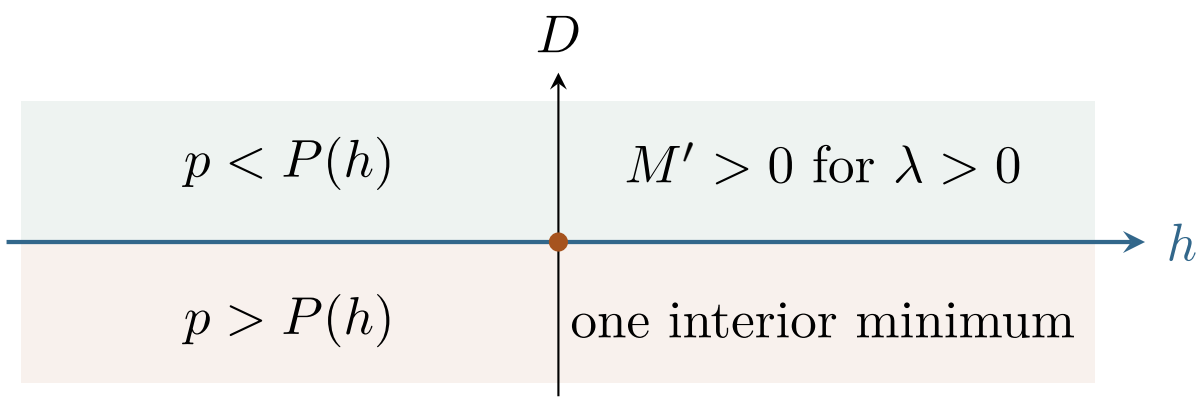}%
\else
\begin{tikzpicture}[x=1cm,y=1cm,>=stealth,font=\small]
 \fill[mapup!8] (-3.55,0) rectangle (3.55,0.93);
 \fill[mapfold!8] (-3.55,-0.93) rectangle (3.55,0);
 \draw[thick,mapdown,->] (-3.65,0)--(3.88,0) node[right] {$h$};
 \draw[thin,->] (0,-1.02)--(0,1.12) node[above] {$D$};
 \fill[mapfold] (0,0) circle[radius=1.8pt];
 \node at (1.75,.52) {$M'>0$ for $\lambda>0$};
 \node at (1.75,-.52) {one interior minimum};
 \node at (-1.78,.52) {$p<P(h)$};
 \node at (-1.78,-.52) {$p>P(h)$};
\end{tikzpicture}
\fi
\caption{Exponents near $(p_\dagger,2p_\dagger-1)$ in the coordinates
$h=q-2p+1$ and $D=M'_{p,q}(0+)$ (Theorem~\ref{thm:global-birth}).
The horizontal axis $D=0$ is the transition curve $p=P(h)$. On this
axis the mass starts horizontally and increases for every $\lambda>0$.
Above the axis it starts with positive slope and keeps increasing.
Below the axis it starts with negative slope and has one minimum.
The vertical axis $h=0$ is the line $q=2p-1$. The corresponding mass
curves appear in panel F of Figure~\ref{fig:mass-frequency-global}.}
\label{fig:mass-frequency-local}
\end{figure}

Figure~\ref{fig:mass-frequency-local} shows the exponents corresponding
to the three curves in panel F of Figure~\ref{fig:mass-frequency-global}.
For $p>P(h)$, the mass $a_0$ occurs at $\lambda=0$ and again at
$\lambda=\lambda_+$.

\subsection{The initial slope and its sign on the explicit line}\label{sec:initial}
For $1<p<q<5$ with $p<7/3$, the centered profile $u_0$ belongs to
$L^2$. Write $w_0$ for the radial solution of
\begin{equation}\label{eq:zero-linearized}
 (-\Delta+pu_0^{p-1}-qu_0^{q-1})w_0=-u_0
\end{equation}
whose behavior at infinity is specified in \cite[Theorem~1.1]{FHnew}.
At fixed exponents, \cite[Theorem~4.8]{FHnew} gives
$\partial_\lambda u_\lambda\to w_0$ in $C^2_{\rm loc}$ as $\lambda\downarrow0$.
For $9/5\le p<7/3$ the slope is $-\infty$ in dimension three
\cite[Theorem~1.2(i) and Proposition~2.6]{FHnew}.
Our mass $M=\int u^2$ is twice the mass used in \cite{FHnew}.
For $1<p<9/5$, we compute the initial slope by estimating the mass
derivative outside a large ball.

\begin{proposition}\label{prop:endder}
Assume $1<p<q<5$ and $p<9/5$. Then
\begin{equation}\label{eq:zero-slope}
 M'(0+)=\lim_{\lambda\downarrow0}M'(\lambda)=2\int u_0w_0\in\R .
\end{equation}
\end{proposition}
\begin{proof}[\normalfont\proofname]
For $\lambda>0$, set $u=u_\lambda$ and $w=\partial_\lambda u_\lambda$.
By \cite[Theorem~1.6]{FH}, local elliptic regularity, and
\cite[Theorem~4.8]{FHnew}, for each $R>0$ we have
\[
 u_\lambda\to u_0\quad\text{uniformly on }[0,R],\qquad
 \partial_\lambda u_\lambda\to w_0\quad\text{in }C^2([0,R]).
\]
Thus we may pass to the limit on bounded intervals in
\[
 M'(\lambda)=8\pi\int_0^\infty r^2u(r)w(r)\,dr .
\]
To control the integral for large $r$, we first bound $u$ by a power of $r$.
Since $u_0(r)\to0$ as $r\to\infty$, choose $R>1$ such that
\[
 u_0(R)^{q-p}<\min\{1/2,p/q\}.
\]
For this fixed $R$, we have $u_\lambda(R)\to u_0(R)$ as
$\lambda\downarrow0$, so the same strict inequality holds for
$u_\lambda(R)$ when $\lambda>0$ is sufficiently small.
Since $u_\lambda$ decreases with $r$ and $q>p$, we obtain
\[
 u_\lambda(r)^{q-p}\le u_\lambda(R)^{q-p}<\min\{1/2,p/q\}
 \qquad (r\ge R).
\]
The profile equation in \eqref{eq:main} then gives
$-\Delta u+u^p/2\le0$. Set $k=1/(p-1)$, so that $2kp=2k+2$. The power
$v(r)=Cr^{-2k}$ satisfies
\[
 -\Delta v+v^p/2=C\big(C^{p-1}/2-2k(2k-1)\big)r^{-2k-2}.
\]
Take $C$ large enough that $-\Delta v+v^p/2\ge0$ and $v(R)\ge u(R)$.
Both functions vanish at infinity. A positive interior maximum of $u-v$
would give
\[
 0\ge-\Delta(u-v)+\tfrac12(u^p-v^p)>0 .
\]
Thus $u(r)\le Cr^{-2/(p-1)}$ on $r\ge R$.

We next use $u^{2-p}$ to bound $|w|$. By \eqref{eq:reverse-linearized},
\[
 Lw=-u,\qquad L=-\Delta+\lambda+pu^{p-1}-qu^{q-1}.
\]
The choice of $R$ gives $\lambda+pu^{p-1}-qu^{q-1}\ge0$ for $r\ge R$.
For $0<\theta<1$,
\[
 -\Delta(u^\theta)=-\theta u^{\theta-1}\Delta u+
         \theta(1-\theta)u^{\theta-2}|\nabla u|^2 .
\]
Substituting $\Delta u=\lambda u-u^q+u^p$ from the profile equation in \eqref{eq:main} gives
\begin{align}
 L(u^\theta)&=(1-\theta)\lambda u^\theta
  +u^{p-1+\theta}\big[(p-\theta)-(q-\theta)u^{q-p}\big]\notag\\
 &\qquad+\theta(1-\theta)u^{\theta-2}|\nabla u|^2 .
 \label{eq:power-comparison}
\end{align}
Choose $\theta=2-p\in(0,1)$. Then $p-1+\theta=1$ and
$p-\theta=2(p-1)$. Increase $R$ until $(q-\theta)u^{q-p}\le p-1$.
It follows that $L(u^{2-p})\ge(p-1)u$, the same power of $u$ that
appears in $Lw=-u$.

At the fixed radius $R$, $w(R)\to w_0(R)$ and $u(R)\to u_0(R)>0$.
Take $C$ so that $Cu(R)^{2-p}\ge|w(R)|$ and $C(p-1)\ge1$. Then
\[
 L(Cu^{2-p}\pm w)\ge\big(C(p-1)-1\big)u\ge0,\qquad
 Cu^{2-p}\pm w\longrightarrow0\quad(r\to\infty).
\]
The exterior maximum principle gives $|w|\le Cu^{2-p}$. Consequently
\[
 r^2|uw|\le Cr^2u^{3-p}\le Cr^{2-2(3-p)/(p-1)},\qquad
 2-\frac{2(3-p)}{p-1}<-1\ \Longleftrightarrow\ p<9/5.
\]
The local uniform convergence of $u$ and $w$ controls the integral on
$[0,R]$. On $[R,\infty)$ the preceding bound for $r^2|uw|$ is
integrable, uniformly in small $\lambda>0$. Applying dominated
convergence to $M'(\lambda)=8\pi\int_0^\infty r^2uw\,dr$ gives
\[
 \lim_{\lambda\downarrow0}M'(\lambda)
 =8\pi\int_0^\infty r^2u_0(r)w_0(r)\,dr
 =2\int_{\R^3}u_0w_0\,dx .
\]
For $\lambda>0$, continuity at zero and the mean value theorem give
$\xi_\lambda\in(0,\lambda)$ such that
\[
 \frac{M(\lambda)-M(0)}{\lambda}=M'(\xi_\lambda)
 \longrightarrow 2\int_{\R^3}u_0 w_0\,dx .
\]
This proves the formula for the right derivative at zero.

\end{proof}

On the line $q=2p-1$ we can compute the initial slope explicitly.

The algebraic profile on the line $q=2p-1$ is, up to constants and spatial
scaling, the optimizer of del Pino and Dolbeault \cite[Theorem~1]{DPD}.
Fixing the coefficients of \eqref{eq:main} determines both constants and
hence the zero-frequency mass. If $q=2p-1$ and $1<p<7/3$, put $k=1/(p-1)>3/4$ and
\[
 C=\left(\frac{2(k+1)}{2k-1}\right)^{k},\qquad b=\frac{k+1}{k(2k-1)^2}.
\]
Then
\begin{equation}\label{eq:explicit-profile}
 u_0(r)=C(1+br^2)^{-k},\qquad
 a_0=C^2b^{-3/2}\pi^{3/2}\frac{\Gamma(2k-3/2)}{\Gamma(2k)} .
\end{equation}
In fact direct differentiation gives
\[
 -\Delta u_0=2k(1-2k)b C^{1-p}u_0^p
             +4k(k+1)b C^{1-q}u_0^q .
\]
The coefficients are $-1$ and $1$, respectively, since
$C^{p-1}=2(k+1)/(2k-1)$ and $q-1=2(p-1)$. To compute the mass and
the slope below, substitute $t=br^2$ in the radial integral. For $l>3/2$,
\begin{equation}\label{eq:beta-moment}
 \begin{aligned}
 J(l)&:=\int_{\R^3}(1+b|x|^2)^{-l}\dd x\\
 &=2\pi b^{-3/2}\int_0^\infty t^{1/2}(1+t)^{-l}\dd t
 =\pi^{3/2}b^{-3/2}\frac{\Gamma(l-3/2)}{\Gamma(l)} .
 \end{aligned}
\end{equation}
Taking $l=2k>3/2$ gives $a_0=C^2J(2k)$ in \eqref{eq:explicit-profile}.

We now use the explicit profile to find where $2\int u_0w_0$ changes
sign along $q=2p-1$.

\begin{proposition}\label{prop:threshold}
Let $q=2p-1$ and $1<p<9/5$. The polynomial $15p^3-38p^2-p+36$
has a unique zero $p_\dagger$ in this interval, with $7/5<p_\dagger<3/2$.
The right derivative $M'(0+)$ is positive for $p<p_\dagger$, zero for
$p=p_\dagger$ and negative for $p>p_\dagger$.
\end{proposition}

\begin{proof}[\normalfont\proofname]
By \eqref{eq:zero-slope}, $M'(0+)=2\int u_0w_0$.
Use $k=1/(p-1)>5/4$ and $u_0=C(1+br^2)^{-k}$ from
\eqref{eq:explicit-profile}. We compute $w_0$ by solving
$L_0w_0=-u_0$, where $L_0=-\Delta+pu_0^{p-1}-qu_0^{q-1}$.

For any real $j$ put $f_j(r)=(1+br^2)^{-j}$. Differentiating twice
in the radial Laplacian and substituting $u_0=C f_k$ gives
\[
 b^{-1} L_0 f_j=
 [2j(1-2j)+2(k+1)(2k-1)] f_{j+1}
 +[4j(j+1)-4(k+1)(k+2)] f_{j+2}.
\]
Substituting $j=k-1,k,k+1$ gives
\[
\begin{aligned}
 b^{-1}L_0 f_{k-1}&=4(3k-2)f_k-8(2k+1)f_{k+1},\\
 b^{-1}L_0 f_k&=2(2k-1)f_{k+1}-8(k+1)f_{k+2},\\
 b^{-1}L_0 f_{k+1}&=-4(k+1)f_{k+2}.
\end{aligned}
\]
The coefficient of $f_k$ in $b^{-1}L_0 f_{k-1}$ fixes the first
coefficient of $w_0$. Cancelling the terms in $f_{k+1}$ and $f_{k+2}$
then gives the other two. The resulting solution is
\begin{equation}\label{eq:explicit-w}
 w_0=\frac Cb\left[-\frac{f_{k-1}}{4(3k-2)}
 -\frac{(2k+1)f_k}{(2k-1)(3k-2)}
 +\frac{2(2k+1)f_{k+1}}{(2k-1)(3k-2)}\right].
\end{equation}
Indeed, the function on the right is regular at the origin and tends to zero
since $k>1$. The limit of $\partial_\lambda u_\lambda$ also decays,
by the comparison $|w_0|\le Cu_0^{2-p}$ in
Proposition~\ref{prop:endder}. Their difference is a regular decaying
radial solution of $L_0v=0$. To apply \cite[Proposition~A.12]{FHnew},
recall that a regular radial solution is determined by $v(0)$. Indeed,
if $v(0)=0$, integrating
\[
 (r^2v')'=r^2(pu_0^{p-1}-qu_0^{q-1})v
\]

twice gives
\[
 \sup_{[0,r]}|v|\le
 \frac{r^2}{6}\sup_{[0,r]}|pu_0^{p-1}-qu_0^{q-1}|\sup_{[0,r]}|v| .
\]
For small $r$ this forces $v=0$, and ordinary uniqueness extends the
identity to every radius. Thus all regular solutions are multiples of
the one with central value one. Proposition~A.12 of \cite{FHnew} shows that this solution
has a nonzero coefficient of the nondecaying mode at infinity.
The regular decaying homogeneous solution is therefore zero, which proves
\eqref{eq:explicit-w}.

Multiply the three terms in \eqref{eq:explicit-w} by $2u_0$ and
integrate. Dividing by $a_0/b=C^2J(2k)/b$ gives
\begin{align*}
 \frac{b M'(0+)}{a_0}
 &=-\frac{1}{2(3k-2)}\,\frac{J(2k-1)}{J(2k)}
   -\frac{2(2k+1)}{(2k-1)(3k-2)}\\
 &\quad+\frac{4(2k+1)}{(2k-1)(3k-2)}\,\frac{J(2k+1)}{J(2k)} .
\end{align*}
In the beta integral~\eqref{eq:beta-moment}, the identity
$\Gamma(z+1)=z\Gamma(z)$ gives
\[\frac{J(2k-1)}{J(2k)}=\frac{2k-1}{2k-5/2},\qquad \frac{J(2k+1)}{J(2k)}=\frac{2k-3/2}{2k}.\]
Thus
\begin{equation}\label{eq:explicit-slope}
 M'(0+)=\frac{a_0}{b}\,
 \frac{12k^3-32k^2+7k+15}{k(2k-1)(3k-2)(4k-5)}.
\end{equation}
All denominator factors are positive for $k>5/4$. The sign is
therefore the sign of
\[
 (p-1)^3(12k^3-32k^2+7k+15)=15p^3-38p^2-p+36 .
\]
Write $f(p)=15p^3-38p^2-p+36$. Its second derivative
$f''(p)=90p-76$ is positive on $[1,9/5]$. On $[1,3/2]$,
\[
 f'(p)\le f'(3/2)=-55/4<0,\qquad
 f(7/5)=32/25>0,\qquad f(3/2)=-3/8<0 .
\]
Thus $f$ has exactly one zero on $[1,3/2]$, and its derivative there
is negative. On $[3/2,9/5]$, convexity gives
\[
 f(p)\le\max\{f(3/2),f(9/5)\}<0,
\]
since $f(9/5)=-36/25$.
The factor multiplying $f(p)$ in \eqref{eq:explicit-slope} is smooth
and positive. This proves the signs of $M'(0+)$ and gives
$d'(p_\dagger)<0$ for $d(p)=M_{p,\,2p-1}'(0+)$.
\end{proof}

On $q=2p-1$, the sufficient condition for a negative slope in
\cite[Theorem~1.3]{FHnew} becomes
\[
 2p-1\ge\frac{43-15p}{3(5-p)}
 \quad\Longleftrightarrow\quad p\ge 4-\frac{\sqrt{57}}{3},
 \qquad 1<p<5/3 .
\]
The exact threshold in Proposition~\ref{prop:threshold} is smaller.
Indeed, at $p=37/25$ the cubic equals $-276/3125$, while
$(43-15p)/(3(5-p))-(2p-1)=8/825>0$.

\subsection{The curvature at the change of sign}\label{sec:birth}

To study the curvature near zero frequency, we first express $M''$ at
positive frequency in terms of $u_\lambda$ and its first derivative.

\begin{lemma}\label{lem:curvature-identity}
Let $1<p<q<5$, $\lambda>0$, $u=u_\lambda$, and
$w=\partial_\lambda u_\lambda$. Then
\begin{equation}\label{eq:second}
 M''(\lambda)=6\int w^2
    +2\int\big[p(p-1)u^{p-2}-q(q-1)u^{q-2}\big]w^3 .
\end{equation}
Both $w$ and $\partial_\lambda^2 u_\lambda$ belong to $H^2_{\rm rad}$.
\end{lemma}
\begin{proof}[\normalfont\proofname]
For any fixed $\lambda>0$, put
\[
 t=\lambda^{-(q-p)/(q-1)},\qquad
 u_\lambda(x)=\lambda^{1/(q-1)}v_t(\sqrt\lambda\,x),\qquad
 h(r)=e^{-r}/(1+r).
\]
By \cite[(2.4)]{FHnew}, $ch\le v_t\le Ch$ for this fixed
frequency. Use the weighted spaces $X,Y$ and the operator
$G=(-\Delta+1)^{-1}$ from Lemma~\ref{lem:weighted}, with
$1+2\tau<\min(p,q)$. A sufficiently small $X$-ball about $v_t$
satisfies $(c/2)h\le v\le(C+c/2)h$. The power estimates in that
lemma therefore make $(v,p,q,t)\mapsto v-G(v^q-tv^p)$ a $C^2$ map
into $X$. Its derivative in $v$ sends $z$ to
\[
 z-G\big((qv_t^{q-1}-tpv_t^{p-1})z\big).
\]
Multiplication by the coefficient maps $X$ continuously into $Y$,
and $G:Y\to X$ is compact. If the displayed derivative vanishes,
its equation gives $z\in H^2_{\rm rad}$, since $Y\subset L^2$ and
$G$ maps $L^2$ into $H^2$. Rescaling $z$ to $z(\sqrt\lambda\,x)$
then gives a radial $H^2$ solution of
\[
 (-\Delta+\lambda+pu_\lambda^{p-1}-qu_\lambda^{q-1})\,
 z(\sqrt\lambda\,x)=0 .
\]
Lemma~\ref{lem:rad} gives $z=0$. The Fredholm alternative gives a
bounded inverse on $X$, so the implicit function theorem gives a local
$C^2$ profile in $(p,q,t)$. Every profile in the chosen ball is
positive. Uniqueness at each frequency then gives the centered
solution after rescaling.

Fix $p,q$ and $\lambda>0$. Write $u=u_\lambda$,
$w=\partial_\lambda u_\lambda$, $z=\partial_\lambda^2u_\lambda$,
and $r=\sqrt{\lambda}\,|x|$. The $C^2$ dependence in the weighted
norm gives
\[
 |v_t(r)|+|\partial_t v_t(r)|+|\partial_t^2 v_t(r)|
 \le C\frac{e^{-r}}{1+r}.
\]
To estimate the spatial derivatives, write the rescaled equation
\eqref{eq:scaled} in radial form and differentiate it in $t$:
\[
 \begin{aligned}
 (r v_t)'' &=
 r\bigl(v_t+t v_t^p-v_t^q\bigr),\\
 (r\partial_t v_t)'' &=
 r\bigl[(1+tpv_t^{p-1}-qv_t^{q-1})\partial_t v_t+v_t^p\bigr].
 \end{aligned}
\]
Here the primes denote derivatives in $r$. Both right-hand sides
are $O(e^{-r})$ as $r\to\infty$. For either $f=rv_t$ or
$f=r\partial_t v_t$, the bound on $f''$ implies that $f'$ has a
finite limit. This limit is zero because $f(r)\to0$. Hence
\[
 f'(r)=-\int_r^\infty f''(s)\,ds=O(e^{-r}).
\]
It follows that $\partial_r v_t$ and $\partial_r\partial_t v_t$
are $O(e^{-r}/r)$. The identity
$(rv_t)''=r\partial_r^2v_t+2\partial_r v_t$ gives the same bound
for $\partial_r^2v_t$.

Now differentiate $u_\lambda(x)=\lambda^{1/(q-1)}v_t(r)$ at fixed $x$.
Since $\partial_\lambda r=r/(2\lambda)$, one frequency derivative
introduces at most one factor of $r$, and two introduce at most $r^2$.
Using $v_t(r)\ge c e^{-r}/(1+r)$, we obtain
\[
 |w(x)|\le C(1+r)u(x),\qquad
 |z(x)|\le C(1+r)^2u(x).
\]
In particular, $w,z\in L^2(\R^3)$ and, for $s=p,q$,
\[
 u^{s-2}w^2\le C(1+r)^2 u^s\in L^2(\R^3).
\]
Let $L=-\Delta+\lambda+pu^{p-1}-qu^{q-1}$. Differentiating
\eqref{eq:main} once and twice in $\lambda$ at fixed $p,q$ gives
\[
 \begin{aligned}
 Lw&=-u,\\
 Lz&=-2w-\big[p(p-1)u^{p-2}-q(q-1)u^{q-2}\big]w^2 .
 \end{aligned}
\]
Here $\partial_\lambda^2(\lambda u)=2w+\lambda z$.
Both right sides are in $L^2$ and the potential is bounded. Thus
$\Delta w,\Delta z\in L^2$ and $w,z\in H^2$. Integration by parts gives
\[
 \int uz=-\int (Lw)z=-\int wLz
 =2\int w^2+\int \big[p(p-1)u^{p-2}-q(q-1)u^{q-2}\big]w^3 .
\]
Substitution in $M''=2\int(w^2+uz)$ gives \eqref{eq:second}.
\end{proof}

We next control $u_\lambda$ and $\partial_\lambda u_\lambda$ as both exponents vary and $\lambda\downarrow0$.

\begin{lemma}\label{lem:joint-endpoint}
Let the exponents vary in a sufficiently small closed rectangle about
$(p_\dagger,2p_\dagger-1)$ contained in $1<p<11/7$, $p<q<2$.
As $(p,q,\lambda)\to(\bar p,\bar q,0)$ with $\lambda>0$, the centered
profiles converge in $C^1_{\rm loc}$ to $u_0$ at $(\bar p,\bar q)$,
and $\partial_\lambda u_\lambda$ converges locally uniformly to $w_0$
at the same pair. There exist $C,R>0$, uniform for the exponents in
the rectangle and all sufficiently small $\lambda>0$, such that for
$r\ge R$,
\[
 u_\lambda(r)\le Cr^{-2/(p-1)},\qquad
 |\partial_\lambda u_\lambda(r)|\le C u_\lambda(r)^{2-p}.
\]
\end{lemma}
\begin{proof}[\normalfont\proofname]
Write $u=u_\lambda$ and $w=\partial_\lambda u_\lambda$. Fix a
limiting exponent pair $(\bar p,\bar q)$ in the rectangle.

For $0<\lambda\le1$, take a fixed nonzero nonnegative
$\psi\in C_c^\infty$. For each $(p,q,\lambda)$, let $\Phi$ and $I$ be
the functionals in \eqref{eq:action-constraint} at these parameter values.
The profile $u_\lambda$ minimizes $\Phi$ among nonzero functions with $I=0$.
Since $q-p$ and $p-1$ have positive lower bounds on the rectangle,
the powers of $s$ in \eqref{eq:action-constraint} give
\[
 I(s\psi)>0\quad(s\text{ small}),\qquad
 I(s\psi)<0\quad(s\text{ large}),\qquad
 \Phi(s\psi)\le Cs^2+Cs^{p+1}-cs^{q+1}.
\]
Here $c>0$ since $\int\psi^{q+1}/(q+1)$ has a positive minimum.
The upper bound is nonpositive for all sufficiently large $s$, uniformly
on the rectangle. The ray meets $I=0$, so minimization gives
\[
 \Phi(u_\lambda)\le \sup_{s>0}\Phi(s\psi)\le C.
\]
Eliminating $\int u^{q+1}$ by $I(u)=0$ gives
\[
 \left(\tfrac12-\tfrac1{q+1}\right)\int(|\nabla u|^2+\lambda u^2)
 +\left(\tfrac1{p+1}-\tfrac1{q+1}\right)\int u^{p+1}=\Phi(u)\le C .
\]
Both coefficients have a positive minimum on the rectangle. Hence,
using the Sobolev inequality,
\[
 \int|\nabla u|^2+\int u^{p+1}\le C,\qquad \|u\|_6\le C .
\]
The largest $q$ on the closed rectangle is less than $2$. Fix $t>3$
so close to $3$ that $tq<6$ on the rectangle. Since $1<p<q$ and
$0<\lambda\le1$, the uniform $L^6$ bound on $u$ gives
\[
 \|u^q-u^p-\lambda u\|_{L^t(B_{2R})}+\|u\|_{L^t(B_{2R})}\le C_R .
\]
Interior elliptic estimates give $\|u\|_{W^{2,t}(B_R)}\le C_R$.
As $t>3$, the profiles are precompact in $C^1$ on every smaller ball.

At the maximum $r=0$, the radial equation~\eqref{eq:radial} gives
$u(0)^q-u(0)^p\ge\lambda u(0)>0$, so $u(0)\ge1$. For each $r>0$,
Cauchy--Schwarz gives the uniform decay bound
\[
 |u(r)|\le\int_r^\infty |u'(s)|\dd s
 \le\left(\int_r^\infty s^2|u'(s)|^2\dd s\right)^{1/2}
     \left(\int_r^\infty s^{-2}\dd s\right)^{1/2}
 \le C r^{-1/2}.
\]
Let $u_n\to U$ in $C^1_{\rm loc}$ along a subsequence with
$(p_n,q_n,\lambda_n)\to(\bar p,\bar q,0)$. Passing to the limit in
the equation and the pointwise bounds gives
\[
 -\Delta U=U^{\bar q}-U^{\bar p},\qquad
 U(0)\ge1,\qquad 0\le U(r)\le Cr^{-1/2}\quad(r>0).
\]
The strong maximum principle gives $U>0$. Fatou's lemma gives
\[
 \int\bigl(|\nabla U|^2+U^{\bar p+1}\bigr)
 \le\liminf_{n\to\infty}\int\bigl(|\nabla u_n|^2+u_n^{p_n+1}\bigr)
 \le C .
\]
The unique centered positive solution in this class is $u_0$.
Every subsequence has a further subsequence with this limit, so
$u_\lambda\to u_0$ in $C^1_{\rm loc}$ as the parameters approach
$(\bar p,\bar q,0)$.

We next bound $w=\partial_\lambda u_\lambda$ uniformly for small $\lambda$.
The uniform decay of $u$, the lower bounds on $p-1,q-p$ and
Proposition~\ref{prop:endder} give one large $R$ such that, for $r\ge R$,
\[
 u(r)\le Cr^{-2/(p-1)},\quad
 \lambda+pu^{p-1}-qu^{q-1}\ge0,\quad L(u^{2-p})\ge(p-1)u .
\]
For the exterior supersolution $Cr^{-2/(p-1)}$ in
Proposition~\ref{prop:endder}, a single $C$ works for the whole rectangle.
Indeed $1/(p-1)$ is bounded, $p-1$ has a positive lower bound, and
the boundary values $u(R)R^{2/(p-1)}$ are uniformly bounded at the
fixed radius $R$. The last
bound is \eqref{eq:power-comparison} with $\theta=2-p$. Since
$u(R)\to u_0(R)>0$ and $Lw=-u$, exterior comparison gives
\begin{equation}\label{eq:wbar}
 |w(r)|\le C(1+|w(R)|)u(r)^{2-p}\qquad(r\ge R).
\end{equation}
We bound the boundary value and the interior at once. If
$\|w_n\|_{L^\infty(B_R)}\to\infty$ along a sequence, put
$v_n=w_n/\|w_n\|_{L^\infty(B_R)}$. Then
\[
 \|v_n\|_{L^\infty(B_R)}=1,\qquad
 L_n v_n=-\frac{u_n}{\|w_n\|_{L^\infty(B_R)}}\longrightarrow0
 \quad\text{locally uniformly}.
\]
By \eqref{eq:wbar}, $|v_n(r)|\le 2C u_n(r)^{2-p_n}$ for $r\ge R$
and all large $n$. Thus $v_n$ is bounded also on every larger ball.
The local convergence of $u$ and of the exponents gives local uniform convergence of the potentials. Interior elliptic
estimates on these larger balls give a subsequence converging uniformly
on compact sets, including $\overline B_R$, to a radial $v$ satisfying
\[
 \begin{gathered}
 (-\Delta+\bar p\,u_0^{\bar p-1}-\bar q\,u_0^{\bar q-1})v=0,
 \qquad \|v\|_{L^\infty(B_R)}=1,\\
 |v(r)|\le 2C u_0(r)^{2-\bar p}\longrightarrow0
 \quad(r\to\infty).
 \end{gathered}
\]
By \cite[Proposition~A.12]{FHnew}, the only regular decaying radial
solution of this homogeneous equation is zero. This contradicts
$\|v\|_{L^\infty(B_R)}=1$. We have therefore proved
\[
 \|w\|_{L^\infty(B_R)}\le C,\qquad
 |w(r)|\le C u(r)^{2-p}\quad(r\ge R),
\]
where the second bound follows from \eqref{eq:wbar}.

Local elliptic estimates now give a locally uniformly convergent
subsequence of $w$. Its limit solves \eqref{eq:zero-linearized} with
source $-u_0$ and tends to zero at infinity. The difference of any two
such limits solves the homogeneous equation and decays, so it vanishes
by the same proposition. At fixed exponents the limit is $w_0$
by \cite[Theorem~4.8]{FHnew}. Thus the limit is also $w_0$ when the
exponents vary. The estimates are uniform on the closed rectangle. Any sequence violating them would have a subsequence whose exponents
converge to a pair in the rectangle, where the argument just given
applies.
\end{proof}

We can now take the zero-frequency limit in the curvature identity.

\begin{lemma}\label{lem:birth-curvature}
There is a constant $B>0$ such that
\[
 M_{p,q}'(\lambda)\longrightarrow0,\qquad
 M_{p,q}''(\lambda)\longrightarrow B
\]
as $(p,q)\to(p_\dagger,2p_\dagger-1)$ and $\lambda\downarrow0$ jointly.
\end{lemma}

\begin{proof}[\normalfont\proofname]
Take $(\bar p,\bar q)=(p_\dagger,2p_\dagger-1)$, and use the rectangle
and the exterior bounds from Lemma~\ref{lem:joint-endpoint}.
Write $u=u_\lambda$ and $w=\partial_\lambda u_\lambda$.
On $r\ge R$, use $u\le1$ and $|w|\le Cu^{2-p}$ to obtain
\begin{align*}
 |w|^2+u^{p-2}|w|^3+u^{q-2}|w|^3
 &\le C\big(u^{4-2p}+u^{p-2+3(2-p)}\big)
 \le C u^{4-2p},\\
 |uw|&\le C u^{3-p}\le C u^{4-2p}.
\end{align*}
Because the rectangle is compactly contained in $p<11/7$,
\[
 u^{4-2p}\le C r^{-2(4-2p)/(p-1)}\le C r^{-3-\delta}
 \quad(r\ge R)
\]
for some $\delta>0$. The radial integral of the bound over $[R,\infty)$ is finite since
$\int_R^\infty r^{-1-\delta}\dd r<\infty$.
On compact balls $u\to u_0>0$ and $w\to w_0$ uniformly. Dominated
convergence therefore gives
\[
 \begin{aligned}
 \lim M'(\lambda)&=2\int u_0w_0=0,\\
 B:=\lim M''(\lambda)&=6\int w_0^2
 +2\int\big[\bar p(\bar p-1)u_0^{\bar p-2}
            -\bar q(\bar q-1)u_0^{\bar q-2}\big]w_0^3 .
 \end{aligned}
\]
The limit of $M'$ is zero by the explicit slope formula~\eqref{eq:explicit-slope} at $p=p_\dagger$.

We use the explicit profiles on $q=2p-1$ to show that $B>0$.
Take $p=p_\dagger$, $q=2p-1$ and $k=1/(p-1)>2$, since
$p_\dagger<3/2$ by Proposition~\ref{prop:threshold}.
Recall $C^{p-1}=2(k+1)/(2k-1)$, $b=(k+1)/(k(2k-1)^2)$ and
$a_0=\int u_0^2$. With $y=(1+br^2)^{-1}$ the profiles from \eqref{eq:explicit-profile} and \eqref{eq:explicit-w} become
\[
 u_0=Cy^k,\qquad w_0=\frac{u_0}{b}H(y),\qquad
 H(y)=-\frac{1}{4(3k-2)y}
      -\frac{2k+1}{(2k-1)(3k-2)}(1-2y).
\]
Substitute these profiles into the mass curvature formula~\eqref{eq:second} at $u_0,w_0$ and
divide by $a_0/b^2$. Since $u_0^{p-1}=C^{p-1}y$ and
$u_0^{q-1}=C^{q-1}y^2$, the result is
\[
 \frac{b^2 B}{a_0}
 =\int_{\R^3}\left\{6H(y)^2
 +2H(y)^3\left[\frac{2(k+1)(2k-1)}{k}y
             -\frac{8(k+1)(k+2)}{k}y^2\right]\right\}
             \frac{u_0^2}{a_0}\dd x .
\]
The powers of $y$ in braces range from $-2$ to $5$. To evaluate each
one let
\[
 J_j=a_0^{-1}\int u_0^2 y^j=J(2k+j)/J(2k).
\]

The beta integral~\eqref{eq:beta-moment} gives
\[
 J_j=\frac{\Gamma(2k+j-3/2)\Gamma(2k)}
           {\Gamma(2k+j)\Gamma(2k-3/2)},\qquad
 J_0=1,\quad J_{j+1}=\frac{2k+j-3/2}{2k+j}\,J_j .
\]
The smallest gamma argument is $2k-7/2>0$. Expand the braces and
replace each $y^j$ by $J_j$, using the recurrence twice backwards
and five times forwards. With common denominator
$k^2(2k-1)^2(3k-2)^3(4k-7)(4k-5)$, the numerator is
\begin{align*}
 &432k^7-3456k^6+9216k^5-7580k^4-3939k^3+6455k^2+150k-1260\\
 &\quad=(12k^3-32k^2+7k+15)(36k^4-192k^3+235k^2+62k-60)
             +72(8k^2-5k-5).
\end{align*}
At $p=p_\dagger$ the first factor on the right is zero by
\eqref{eq:explicit-slope}. Hence
\[
 \frac{b^2 B}{a_0}
 =\frac{72(8k^2-5k-5)}
 {k^2(2k-1)^2(3k-2)^3(4k-7)(4k-5)}>0 ,
\]
because $k>2$.

\end{proof}

\subsection{The transition at small frequencies}
To apply the implicit function theorem to $M'_{p,q}(0+)=0$, we
first differentiate the initial slope with respect to the exponents.

\begin{lemma}\label{lem:initial-slope-regularity}
Near $(p_\dagger,2p_\dagger-1)$, the function $D(p,q)=M'_{p,q}(0+)$
is $C^1$. For $s=p,q$, its derivative $\partial_sD$ is the continuous
limit of $\partial_s M'_{p,q}(\lambda)$ as $\lambda\downarrow0$.
\end{lemma}
\begin{proof}[\normalfont\proofname]
Work on a closed rectangle about $(p_\dagger,2p_\dagger-1)$ contained
in $1<p<3/2$, $p<q<2$.

\noindent Step 1. We compute the derivatives of the slope with respect to $p$ and $q$ at positive frequency.
Fix $\lambda>0$. Write $u=u_\lambda$, $w=\partial_\lambda u$ and
\[
 L=-\Delta+\lambda+pu^{p-1}-qu^{q-1}.
\]
Let $s$ denote the exponent
being varied, with the other exponent held fixed, and put
$z=\partial_s u$. The derivative equations for the two choices are
\begin{equation}\label{eq:exponent-linearized}
 L(\partial_p u)=-u^p\log u,\qquad
 L(\partial_q u)=u^q\log u,\qquad Lw=-u .
\end{equation}
The weighted $C^2$ argument in the proof of Lemma~\ref{lem:curvature-identity}
justifies $\partial_s w$ and integration by parts at fixed $\lambda>0$.
Differentiating the frequency equation $Lw=-u$ in \eqref{eq:exponent-linearized} with respect to $s$ gives
\[
 L(\partial_s w)=-z-(\partial_s L)w,\qquad
 \int u\,\partial_s w=\int wz+\int w^2\,\partial_s(pu^{p-1}-qu^{q-1}).
\]
Consequently the derivative of $2\int uw$ is
\begin{align}
 \partial_s M'_{p,q}(\lambda)
 &=4\int wz+2\int w^2
   \left.\frac{\partial}{\partial s}(p v^{p-1}-qv^{q-1})\right|_{v=u}
   \notag\\
 &\quad+2\int w^2\big[p(p-1)u^{p-2}-q(q-1)u^{q-2}\big]z .
 \label{eq:birth-exponent}
\end{align}
In the first line $v$ is held fixed. That derivative equals
$u^{p-1}(1+p\log u)$ for $s=p$, and $-u^{q-1}(1+q\log u)$ for $s=q$.

\noindent Step 2. We bound $\partial_p u$ and $\partial_q u$ as $\lambda\downarrow0$.
Let $(p,q,\lambda)\to(\bar p,\bar q,0)$ in the rectangle. Lemma~\ref{lem:joint-endpoint} gives $u\to u_0$, $w\to w_0$ locally and
\[
 u(r)\le Cr^{-2/(p-1)},\qquad |w(r)|\le C u(r)^{2-p}\qquad(r\ge R).
\]
For $z=\partial_s u$ use $\theta=3/4$ in \eqref{eq:power-comparison}.
Since $q-p$ is bounded below and $u$ is uniformly small at large radii,
\[
 L(u^{3/4})\ge c u^{p-1/4}
       \ge c' u^s|\log u|=c'|Lz|,\qquad r\ge R .
\]
The second bound follows from $s\ge p$ and boundedness of
$u^{1/4}|\log u|$ on $(0,1]$. The exterior potential is nonnegative
and $u(R)\to u_0(R)>0$. Comparing both signs as in \eqref{eq:wbar} gives
\[
 |z(r)|\le C(1+|z(R)|)\,u(r)^{3/4},\qquad r\ge R .
\]
If $\|z_n\|_{L^\infty(B_R)}\to\infty$, normalize by this norm as in
the proof of Lemma~\ref{lem:joint-endpoint}. The normalized source tends locally
to zero since $u_n^s\log u_n$ is locally bounded. The exterior bound
controls the normalized solutions on larger balls, so interior estimates
give a regular radial limit $v$ with
\[
 L_0 v=0,\qquad \|v\|_{L^\infty(B_R)}=1,\qquad
 |v(r)|\le2C u_0(r)^{3/4}\longrightarrow0 .
\]
This contradicts \cite[Proposition~A.12]{FHnew}. Thus $z$ is bounded
on compact sets and $|z|\le Cu^{3/4}$ on the exterior. Each local limit
solves the corresponding exponent derivative equation from \eqref{eq:exponent-linearized} at $(\bar p,\bar q,0)$ and
decays. The difference of two limits solves $L_0v=0$ and decays, so
the limit is unique. The estimates are uniform on a smaller closed
rectangle by the same argument applied to any convergent sequence of
exponent pairs.

\noindent Step 3. We pass to $\lambda=0$ in the two slope derivatives.
For $r\ge R$, use $u\le1$, $|w|\le Cu^{2-p}$ and $|z|\le Cu^{3/4}$.
The terms in \eqref{eq:birth-exponent} then satisfy
\begin{align*}
 |wz|&\le Cu^{11/4-p},\\
 |w|^2|z|(u^{p-2}+u^{q-2})&\le Cu^{2(2-p)+3/4+p-2}=Cu^{11/4-p},\\
 |w|^2u^{s-1}(1+s|\log u|)
 &\le Cu^{3-p}(1+|\log u|)\le Cu^{11/4-p}.
\end{align*}
For the logarithmic term we used that $u^{1/4}(1+|\log u|)$ is
bounded on $(0,1]$. Since $p<3/2$, the exterior decay estimate gives
\[
 u^{11/4-p}\le C r^{-2(11/4-p)/(p-1)}\le C r^{-5},
 \qquad \int_R^\infty r^2 r^{-5}\,dr<\infty .
\]
On each fixed ball $u,w,z$ converge locally uniformly, and $u_0>0$.
Dominated convergence in \eqref{eq:birth-exponent} gives a joint
limit of each exponent derivative at $(\bar p,\bar q,0)$.
The bound $|uw|\le Cu^{3-p}$ gives, in the same way,
\[
 M'_{p,q}(\lambda)\longrightarrow
 2\int u_0w_0=D(\bar p,\bar q)
 \quad\text{as }(p,q,\lambda)\to(\bar p,\bar q,0),\quad\lambda>0 .
\]

For $G=M'$, $\partial_p M'$ or $\partial_q M'$, let $G(p,q,0)$
denote the limit just obtained. To see that these boundary values are
continuous, fix an interior pair $(p,q)$ and $\varepsilon>0$. Joint
convergence gives $\delta>0$ such that
\[
 |G(p',q',\lambda)-G(p,q,0)|<\varepsilon
 \quad\text{if }|p'-p|+|q'-q|<\delta,\quad 0<\lambda<\delta .
\]
Letting $\lambda\downarrow0$ at fixed $(p',q')$ gives the same bound
with $\le\varepsilon$ for $G(p',q',0)$. Thus $G$ is continuous up to $\lambda=0$. On each compact rectangle
in this neighborhood the extension is uniformly continuous.

For a segment inside the exponent rectangle, the fundamental theorem
of calculus at $\lambda>0$ and uniform convergence on the segment give
\[
 \begin{aligned}
 D(p+\xi,q)-D(p,q)
 &=\lim_{\lambda\downarrow0}\int_0^\xi
          \partial_p M'_{p+t,q}(\lambda)\dd t\\
 &=\int_0^\xi \lim_{\lambda\downarrow0}
          \partial_p M'_{p+t,q}(\lambda)\dd t .
 \end{aligned}
\]
Divide by $\xi$ and let $\xi\to0$. Integration along the $q$-axis
gives the same result for the other exponent. Thus
\[
 \partial_pD(p,q)=\lim_{\lambda\downarrow0}\partial_p M'_{p,q}(\lambda),
 \qquad
 \partial_qD(p,q)=\lim_{\lambda\downarrow0}\partial_q M'_{p,q}(\lambda).
\]
Both limits are continuous in $(p,q)$.

\end{proof}

The sign of the initial slope and the positive curvature now determine
where the mass has a minimum at small frequency.

\begin{proposition}\label{prop:birth-local}
Near $(p_\dagger,2p_\dagger-1)$ the right derivative $M_{p,q}'(0+)$ is
$C^1$ in $(p,q)$. With $h=q-2p+1$, its zero set is a curve $p=P(h)$
with $P(0)=p_\dagger$. The slope is positive for $p<P(h)$ and negative
for $p>P(h)$.

There is a common $\delta>0$ such that $M''_{p,q}(\lambda)>0$ for
$0<\lambda\le\delta$ and $M_{p,q}(\delta)>a_0$. If $p\le P(h)$, then
$M'_{p,q}>0$ on $(0,\delta]$. If $p>P(h)$, there is exactly one zero
$\lambda_{\min}$ of $M'_{p,q}$ in $(0,\delta)$ and exactly one
$\lambda_+\in(\lambda_{\min},\delta)$ with $M_{p,q}(\lambda_+)=a_0$.
As $(p,q)\to(p_\dagger,2p_\dagger-1)$ with $p>P(h)$, we have
$\lambda_{\min}\to0$.
\end{proposition}
\begin{proof}[\normalfont\proofname]
Write $D(p,q)=M'_{p,q}(0+)$ and $d(p)=D(p,2p-1)$.
Lemma~\ref{lem:initial-slope-regularity} gives $D\in C^1$.
By Proposition~\ref{prop:threshold} and \eqref{eq:explicit-slope},
\[
 d(p_\dagger)=0,\qquad d'(p_\dagger)<0 .
\]
Fix $h=q-2p+1$. Differentiating the slope with $h$ fixed gives
\[
 \frac{d}{dp}D(p,2p-1+h)
 =(\partial_p D+2\partial_q D)(p,2p-1+h).
\]
At $(p,h)=(p_\dagger,0)$ this derivative is $d'(p_\dagger)<0$,
so it remains negative nearby. The implicit function theorem gives
a unique $C^1$ function $P$ such that
\[
 P(0)=p_\dagger,\qquad D(P(h),2P(h)-1+h)=0 .
\]
For each fixed $h$ the slope decreases with $p$. It is therefore
positive for $p<P(h)$ and negative for $p>P(h)$.

We now count the stationary frequencies near zero. For the current
pair, abbreviate $D=M'_{p,q}(0+)$ and $a_0=M_{p,q}(0)$.
By Lemma~\ref{lem:birth-curvature}, $M''_{p,q}(\lambda)\to B>0$ jointly at the
reference pair and zero frequency. Choose a neighborhood and
$\delta>0$ such that $M''>B/2$ on $0<\lambda\le\delta$. After fixing
$\delta$, shrink the exponent neighborhood so that $|D|<B\delta/8$.
Since $M'(\lambda)\to D$ and $M(\lambda)\to a_0$, integration from
zero gives
\begin{equation}\label{eq:local-mass-bounds}
 M'(\lambda)\ge D+B\lambda/2,\qquad
 M(\lambda)-a_0\ge D\lambda+B\lambda^2/4 .
\end{equation}
In particular $M'(\delta)>0$ and $M(\delta)>a_0$ for every pair in
this neighborhood. If $D\ge0$, then $M'>0$ for $0<\lambda\le\delta$,
including the case $D=0$. If $D<0$, $M'$ increases strictly from $D$
to a positive value and has exactly one zero $\lambda_{\min}$.
The mass decreases from $a_0$ to $M(\lambda_{\min})<a_0$ and then
increases to $M(\delta)>a_0$. Thus it returns to $a_0$ at exactly one
$\lambda_+\in(\lambda_{\min},\delta)$.

At the minimum, \eqref{eq:local-mass-bounds} gives
$0<\lambda_{\min}\le 2|D|/B$. Since $D\to0$ as the exponents tend
to the transition point, $\lambda_{\min}\to0$.
\end{proof}

The endpoint formula also gives the tangent to the transition curve.
For each $h=q-2p+1$ near zero, $p=P(h)$ is the unique nearby solution of
\begin{equation}\label{eq:zero-slope-implicit}
 D(p,2p-1+h)=0,\qquad
 D(p,q)=8\pi\int_0^\infty r^2u_0(r)w_0(r)\dd r .
\end{equation}
Here $u_0$ is the zero-frequency profile at $(p,q)$, and $w_0$
solves \eqref{eq:zero-linearized} at the same exponents. Differentiating
the equation for $P(h)$ gives
\[
 0=\frac{d}{dh}D(P(h),2P(h)-1+h)
   =(\partial_pD+2\partial_qD)P'(h)+\partial_qD .
\]
At $h=0$, the coefficient of $P'(0)$ is $d'(p_\dagger)<0$. Therefore
\begin{equation}\label{eq:zero-slope-tangent}
 P'(0)=-\frac{\partial_q D(p_\dagger,2p_\dagger-1)}{d'(p_\dagger)},
 \qquad d(p)=D(p,2p-1).
\end{equation}

To describe the two sides of the transition curve, use
\[
 (p,q)\longmapsto(h,D)=(q-2p+1,D(p,q)).
\]
At $(p_\dagger,2p_\dagger-1)$ the Jacobian determinant is
\[
 \det\begin{pmatrix}-2&1\\ \partial_pD&\partial_qD\end{pmatrix}
 =-(\partial_pD+2\partial_qD)=-d'(p_\dagger)>0 .
\]
The inverse function theorem gives a $C^1$ change of coordinates.
In these coordinates $D=0$ is the transition curve and $h=0$ is the
line $q=2p-1$, as shown in Figure~\ref{fig:mass-frequency-local}.

\subsection{An explicit frequency above all decreasing intervals}
\label{sec:explicit-upper-frequency}
We first show that $\partial_\lambda u_\lambda(r)$ changes sign
exactly once as $r$ increases.

\begin{lemma}\label{lem:w-one-node}
Suppose $1<p<q<5$ and $\lambda>0$. The frequency derivative
$w=\partial_\lambda u_\lambda$ has a unique zero $r_w>0$.
It satisfies
\[
 w(0)>0,\qquad w>0\text{ on }(0,r_w),\qquad
 w'(r_w)<0,\qquad w<0\text{ on }(r_w,\infty).
\]
\end{lemma}
\begin{proof}[\normalfont\proofname]
Write $u=u_\lambda$, $U=ru$, $W=rw$, and
\[
 \mathcal L=-\frac{d^2}{dr^2}+V(r),\qquad
 V=\lambda+pu^{p-1}-qu^{q-1}.
\]
The radial operator $L$ satisfies $\mathcal L(rf)=rLf$ and
$\mathcal LW=-U$ by \eqref{eq:reverse-linearized}. The Dirichlet
operator $\mathcal L$ on $(0,\infty)$ has one negative eigenvalue
and zero kernel. Indeed $\langle Lu,u\rangle<0$, and the quadratic
form is nonnegative on the codimension-one space in
\eqref{eq:tangent}. Lemma~\ref{lem:rad} excludes a kernel. These
properties transfer to the half-line under $f\mapsto rf$.

\noindent Step 1. We construct a decaying solution $\Psi$ of $\mathcal L\Psi=0$ with exactly one zero.
Let $E(s)$ be the infimum of the Dirichlet
Rayleigh quotient on $(s,\infty)$. Translation to $(0,\infty)$
gives
\[
 |E(s)-E(t)|\le\|V(s+\cdot)-V(t+\cdot)\|_\infty.
\]

The potential is uniformly continuous and tends to $\lambda>0$. Hence
$E$ is continuous, $E(0)<0$, and $E(s)>0$ for large $s$. Choose $R$
with $E(R)=0$. Multiplication by $V-\lambda$ is relatively compact
with respect to the Dirichlet Laplacian by local compactness and decay
at infinity. The essential spectrum on $(R,\infty)$ starts at
$\lambda$, so zero is the lowest eigenvalue there. Taking the absolute
value of a lowest eigenfunction and using uniqueness for the ODE gives
\[
 \Psi(R)=0,\qquad \Psi'(R)>0,\qquad \Psi(r)>0\quad(r>R).
\]

Extend $\Psi$ to zero by the ODE. If it vanished at $R_1\in[0,R)$,
restrict $\Psi$ to $(R_1,R)$ and $(R,\infty)$ and extend both
restrictions by zero. Both have zero quadratic form by integration
by parts. Since the interiors of their supports are disjoint, every
linear combination also has zero form. Their span has dimension two
and therefore contains a nonzero function orthogonal to the single
negative eigenfunction. The form is strictly positive on that
orthogonal complement because zero is below the essential spectrum
and is not an eigenvalue on $(0,\infty)$. This contradiction shows
that $\Psi<0$ on $[0,R)$.

\noindent Step 2. We prove that $\int_r^\infty\Psi(s)U(s)\,ds>0$ at every radius.
Recall $S(s)=s^{q-1}-s^{p-1}$ and put $Z=r^2u'$. The radial equation
\eqref{eq:radial} gives
\[
 L(ru')=-2\Delta u=2(S(u)-\lambda)u,\qquad
 \mathcal LZ=2(S(u)-\lambda)U .
\]
At zero $Z=Z'=0$. To check the boundary at infinity, the profile estimate
\cite[(2.4)]{FHnew} gives $U=O(e^{-\sqrt\lambda r})$, and
$U''=(\lambda-S(u))U$ gives the same bound on $U'$ by integration.
Thus $Z=rU'-U$ and $Z'=rU''$ tend to zero. The eigenfunction equation
gives the same limits for $\Psi,\Psi'$.
Green's formula on $(0,\infty)$ and on $(R,\infty)$ now gives
\[
 \int_0^\infty\Psi U(S(u)-\lambda)\dd r=0,\qquad
 2\int_R^\infty\Psi U(S(u)-\lambda)\dd r
       =-Z(R)\Psi'(R)>0 .
\]
Consequently $S(u)>\lambda$ somewhere beyond $R$.
The function $S(s)$ is nonpositive for $s\le1$ and increases
strictly for $s\ge1$. Since $u$ decreases strictly,
$S(u(R))>\lambda$. Moreover
$S(u)-S(u(R))$ is positive before $R$ and negative after $R$.
Its product with $\Psi U$ is strictly negative on both intervals.
Using $\int_0^\infty\Psi U(S(u)-\lambda)\dd r=0$ gives
\[
 (S(u(R))-\lambda)\int_0^\infty\Psi U\dd r
 =-\int_0^\infty\Psi U[S(u)-S(u(R))]\dd r>0 .
\]
Define $J(r)=\int_r^\infty\Psi(s)U(s)\dd s$. The last inequality gives
$J(0)>0$. On $(0,R)$ we have $J'=-\Psi U>0$, and for $r\ge R$ the
integral defining $J$ is positive. Hence $J(r)>0$ at every radius.

\noindent Step 3. We locate the zero of $w$ by following the quotient $W/\Psi$ on each side of $R$.
The equations for $W$ and $\Psi$ give
\[
 \Psi W'-W\Psi'=-J,\qquad
 (W/\Psi)'=-J/\Psi^2<0\quad(r\ne R).
\]
The Wronskian vanishes at infinity, since $W\in H^2(0,\infty)$
by $w\in H^2_{\rm rad}$ and $\mathcal LW=-U$.
At zero, $W(0)=0$ and $\Psi(0)W'(0)=-J(0)<0$.
Thus $w(0)=W'(0)>0$. On $(0,R)$ the quotient $W/\Psi$ decreases
from zero, giving $W>0$. At $R$ the Wronskian gives $W(R)=J(R)/\Psi'(R)>0$.

On $(R,\infty)$ the same quotient decreases from $+\infty$.
Choose $r_0>R$ so large that $\lambda\le V\le\lambda+1$ and
\[
 \mathcal L U
 =\big[(p-1)u^{p-1}-(q-1)u^{q-1}\big]U\ge0
       \qquad(r\ge r_0).
\]
Such an $r_0$ exists because $u\to0$ and $q>p$. Both $U$ and
$\Psi$ tend to zero. Comparison on $[r_0,\infty)$ therefore gives
$U\ge c\Psi$, with $c=U(r_0)/\Psi(r_0)>0$.
For $s\ge r\ge r_0$ the function
$\Psi(r)e^{-\sqrt{\lambda+1}(s-r)}$ has nonpositive image under
$\mathcal L$, agrees with $\Psi$ at $r$ and tends to zero.
A second comparison and integration give
\[
 J(r)\ge c\int_r^\infty\Psi(s)^2\,ds
       \ge \frac{c\Psi(r)^2}{2\sqrt{\lambda+1}}.
\]
For $r\ge r_0$ it follows that
\[
 \left(\frac W\Psi\right)'=-\frac{J}{\Psi^2}
 \le-\frac{c}{2\sqrt{\lambda+1}} .
\]
Thus $W/\Psi$ decreases from $+\infty$ at $R$ to $-\infty$ and has
exactly one zero $r_w>R$. At that zero,
$W'(r_w)=-J(r_w)/\Psi(r_w)<0$. Since $w=W/r$ for $r>0$, the sign
and the direction of crossing are the same for $w$.
\end{proof}

We use the opposite signs of $\partial_\lambda u_\lambda$ inside
and outside the ball of radius $r_w$ to bound every frequency at which
$M'(\lambda)\le0$.

\begin{proposition}\label{prop:explicit-upper-frequency}
For $1<p<q<7/3$, set
\begin{equation}\label{eq:explicit-upper-frequency}
 \Lambda_{p,q}=
 \frac{3(q-p)}{7-3q}
 \left(\frac{7-3p}{7-3q}\right)^{(p-1)/(q-p)} .
\end{equation}
Then $M'(\lambda)>0$ for every $\lambda\ge\Lambda_{p,q}$.
Let $m=\min_{\lambda\ge0}M(\lambda)$ and include the zero-frequency
solution at its mass $a_0$. A solution exists exactly when $a\ge m$.
Every sufficiently large mass admits exactly one solution up to
translations. If $m<a_0$, each $m<a\le a_0$ admits at least two.
\end{proposition}
\begin{proof}[\normalfont\proofname]
Assume $M'(\lambda)\le0$ and set $I_s=\int_{\R^3}u^s w\,dx$,
so $I_1=M'(\lambda)/2\le0$.
From $Lu=-(q-1)u^q+(p-1)u^p$ and $L(ru')=-2\Delta u$ we get
\begin{align*}
 L\left(\frac{u}{q-1}+\frac{ru'}2\right)
   &=-\lambda u-\frac{q-p}{q-1}u^p,\\
 L\left(\frac{u}{p-1}+\frac{ru'}2\right)
   &=-\lambda u-\frac{q-p}{p-1}u^q .
\end{align*}
Pair each identity with $w$, use $Lw=-u$ and
$\int u\,ru'=-3M/2$. The result is
\begin{align}
 I_p&=\frac{7-3q}{4(q-p)}M-\frac{q-1}{q-p}\lambda I_1>0,
 \label{eq:stationary-moments}\\
 I_q&=\frac{7-3p}{4(q-p)}M-\frac{p-1}{q-p}\lambda I_1 .
 \notag
\end{align}
In particular
\[
 I_q\le \rho I_p,\qquad \rho=\frac{7-3p}{7-3q}>1,
\]
since $I_q-\rho I_p=4\lambda I_1/(7-3q)$.

Let $B$ be the ball of radius $r_w$ from Lemma~\ref{lem:w-one-node}
and $E=\R^3\setminus\overline B$. Put $X=u^{p-1}$,
$\gamma=(q-1)/(p-1)>1$, and denote the average of $X$ over
the positive part of $uw$ by
\[
 b=\frac{\int_B Xuw}{\int_B uw}>0 .
\]
On $B$, $uw>0$, so Jensen's inequality gives
\[
 \int_B X^\gamma uw\ge b^\gamma\int_B uw
 =b^{\gamma-1}\int_B Xuw .
\]
On $E$ we have $X<X(r_w)<b$ and $uw<0$. Thus
$(X^\gamma-b^{\gamma-1}X)uw>0$ on $E$. Adding the two integrals
gives $I_q>b^{\gamma-1}I_p$. Since $I_q\le\rho I_p$ and $I_p>0$, we obtain
\begin{equation}\label{eq:node-amplitude-limit}
 b<\rho^{1/(\gamma-1)} .
\end{equation}

Write $h=-w>0$ on $E$. From $\int_B uw\le\int_E uh$ we obtain
\begin{equation}\label{eq:exterior-moment}
 I_p=b\int_B uw-\int_E Xuh\le\int_E(b-X)uh .
\end{equation}
Suppose $\lambda\ge(\rho-1)b$. On $E$ we have $X<b$ and
$b^{\gamma-1}<\rho$ by \eqref{eq:node-amplitude-limit}. Therefore
\[
 \lambda+X-X^\gamma-(\rho-1)(b-X)
 =\lambda-(\rho-1)b+X(\rho-X^{\gamma-1})\ge0 .
\]
Take $f=-ru'/2>0$. Then $Lf=(\lambda+X-X^\gamma)u$,
while $Lh=u$ on $E$. At its inner boundary $h(r_w)=0$ and
$h'(r_w)>0$. Green's formula gives
\begin{align*}
 (\rho-1)I_p
 &\le\int_E hLf
   =\int_E fu-4\pi r_w^2 f(r_w)h'(r_w)\\
 &<\int_{\R^3}fu=\tfrac34 M .
\end{align*}
All boundary terms at infinity vanish by the decay used in the lemma.
But \eqref{eq:stationary-moments} gives
$(\rho-1)I_p\ge3M/4$, a contradiction.
Thus every $\lambda$ with $M'(\lambda)\le0$ satisfies
\[
 \lambda<(\rho-1)b<(\rho-1)\rho^{1/(\gamma-1)}=\Lambda_{p,q}.
\]

It remains to count the masses covered by these estimates. By
Proposition~\ref{prop:reverse-ends}, $M$ is positive and continuous on
$[0,\infty)$ and tends to infinity. It attains a positive minimum $m$
and its range is $[m,\infty)$. For
\[
 a>\max_{0\le\lambda\le\Lambda_{p,q}} M(\lambda)
\]
every frequency with mass $a$ lies above $\Lambda_{p,q}$, where
strict increase gives exactly one. Finally, if $m<a_0$, any minimizing
frequency $\lambda_m$ is positive. For $m<a\le a_0$ the intermediate
value theorem gives one frequency in $[0,\lambda_m)$ and another in
$(\lambda_m,\infty)$. At $a=a_0$ the first may be taken to be zero.
\end{proof}

\subsection{Positive frequencies on the explicit line and near the transition}
On $q=2p-1$, convexity of a power of $u_\lambda$ gives an upper
bound for the moment used to exclude a stationary frequency.

\begin{lemma}\label{lem:line-convex}
Let $1<p\le p_\dagger$, $q=2p-1$, $k=1/(p-1)$ and $\lambda>0$.
Write $u=u_\lambda$ and $z(r)=u(r)^{-1/k}$. Then
$t\mapsto z(\sqrt t)$ is convex, and
\[
 \frac{\int u^2z}{M(\lambda)}<
 \frac{(2k-1)^2}{(k+1)(4k-5)} .
\]
\end{lemma}
\begin{proof}[\normalfont\proofname]
Since $p_\dagger<3/2$, we have $k>2$. Set
\[
 z(r)=u(r)^{-1/k},\qquad Z(t)=z(\sqrt t),\qquad t=r^2 .
\]
Primes on $Z$ denote derivatives in $t$. The radial profile is
smooth in $r^2$ near zero. For $r>0$, radial monotonicity gives $u'\le0$.
If $u'(r)=0$ at an interior radius, then $u'$ attains a maximum
there, so $u''(r)=0$. The radial equation~\eqref{eq:radial} then
has the constant solution with the same value and derivative there.
Uniqueness for the radial equation away from zero would make $u$
constant, contrary to its decay. Thus $u'<0$ and $Z'(t)>0$ for
$t>0$. At the origin, \eqref{eq:center-strict} gives
\[
 -3u''(0)=u(0)[u(0)^{q-1}-u(0)^{p-1}-\lambda]>0.
\]

Differentiating $Z(r^2)=u(r)^{-1/k}$ twice at $r=0$, using $u'(0)=0$, gives
\[
 2Z'(0)=-\frac1k u(0)^{-1-1/k}u''(0)>0 .
\]

Substituting $u=z^{-k}$ in \eqref{eq:radial} gives
\begin{equation}\label{eq:global-pressure}
 z\Delta z-(k+1)|\nabla z|^2=(1-z-\lambda z^2)/k,
 \qquad \Delta z=6Z'+4tZ'' .
\end{equation}
Differentiate the equation for $Z$ in $t$. With
$T=(4k-2)Z'-(1+2\lambda Z)/k$ this gives
\begin{equation}\label{eq:global-convex}
 4tZ Z'''+[10Z-4(2k+1)tZ']Z''=Z'T,
 \qquad T'=(4k-2)Z''-(2\lambda/k)Z' .
\end{equation}
We claim $Z''\ge0$. On a component $(s,t_1)$ of $\{Z''<0\}$ one
has $T(s)\le0$. For $s>0$ use $Z''(s)=0$, $Z'''(s)\le0$ in
\eqref{eq:global-convex}. For $s=0$ use
$10Z(0)Z''(0)=Z'(0)T(0)$ and $Z''(0)\le0$.
Now $T'<0$ on the component. At a finite right endpoint one would
have $Z''(t_1)=0$, $Z'''(t_1)\ge0$, while
$4t_1Z(t_1)Z'''(t_1)=Z'(t_1)T(t_1)<0$.
Thus $t_1=\infty$ and eventual concavity would give $Z(t)=O(1+t)$.
This is impossible. For large $r$ we have $u\le1$ and
$(-\Delta+\lambda)u\le0$. Comparison outside a ball with a multiple
of $e^{-\sqrt\lambda r}/r$ gives $Z(r^2)\ge c r^{1/k}
e^{\sqrt\lambda r/k}$. This proves the claim.

At the origin \eqref{eq:global-convex} now gives $T(0)\ge0$, so
\[
 Z'(t)\ge Z'(0)\ge\frac{1+2\lambda Z(0)}{k(4k-2)}
 >\frac1{k(4k-2)} .
\]
Multiply the bound for $Z'(0)$ by $6kZ(0)$ and use
$6kZ(0)Z'(0)=1-Z(0)-\lambda Z(0)^2$ from \eqref{eq:global-pressure}.
Collecting the terms with $Z(0)$ and $\lambda Z(0)^2$ gives
\begin{equation}\label{eq:global-center-bound}
 (2k+2)Z(0)+(2k+5)\lambda Z(0)^2\le 2k-1 .
\end{equation}
In particular, $Z(0)<(2k-1)/(2k+2)$ for $\lambda>0$.

We now integrate the convexity bound to estimate the moment of $z$.
Convexity gives
\[
 r z'(r)=2t Z'(t)\ge2\bigl(z(r)-z(0)\bigr).
\]
Integrate $\operatorname{div}(x z^{1-2k})$ on $\R^3$.
The boundary flux tends to zero since $z^{1-2k}=u^{2-1/k}$ decays
exponentially. Hence
\[
 3\int z^{1-2k}
  =(2k-1)\int z^{-2k}rz'(r)
  \ge 2(2k-1)\left(\int z^{1-2k}-z(0)\int z^{-2k}\right).
\]
Since $z^{-2k}=u^2$, inequality \eqref{eq:global-center-bound} gives
\begin{equation}\label{eq:global-upper-moment}
 \frac{\int u^2z}{M}\le\frac{4k-2}{4k-5}\,z(0)
 <\frac{(2k-1)^2}{(k+1)(4k-5)} .
\end{equation}

\end{proof}

At a stationary frequency the same moment satisfies a lower bound.
Comparing the two bounds determines the sign of $M'$ on the line.

\begin{lemma}\label{lem:birth-line}
Let $q=2p-1$ and $1<p\le p_\dagger$. Then
$M'_{p,q}(\lambda)>0$ for every $\lambda>0$.
\end{lemma}
\begin{proof}[\normalfont\proofname]
Use $k,u,z$ from Lemma~\ref{lem:line-convex} and set
$Z(t)=z(\sqrt t)$.
Suppose $M'(\lambda)=0$. With $w=\partial_\lambda u_\lambda$ we have
$\int uw=0$. The scaling identity
\eqref{eq:stationary-moments} gives
\begin{equation}\label{eq:global-first-moment}
 \int u^{1+1/k}w=(k-3/2)M .
\end{equation}
We test $Lw=-u$ with $uz=u^{1-1/k}$. The product rule first gives
\[
 \frac{L(uz)}u=\frac1k-\frac{2}{kz}
             +2k\frac{|\nabla z|^2}{z}-\Delta z .
\]
Eliminate $|\nabla z|^2/z$ by \eqref{eq:global-pressure} to obtain
\begin{equation}\label{eq:global-power-test}
 \frac{L(uz)}u
   =G(r)+\frac{3k+1}{k(k+1)}
          -\frac{4k+2}{k(k+1)}\frac1{z(r)},\qquad
 G(r)=\frac{k-1}{k+1}\Delta z(r)+\frac{2\lambda}{k+1}z(r).
\end{equation}
The function $G$ increases strictly with $r$. Indeed,
$\Delta z=6Z'+4tZ''$ and \eqref{eq:global-convex} give
\[
 \frac{k+1}{k-1}\frac{d}{dt}G(\sqrt t)
 =\frac{Z'}Z
  \left[4(2k+1)tZ''+(4k-2)Z'-\frac1k
                 +\frac{2\lambda Z}{k(k-1)}\right]>0 .
\]
Indeed $Z''\ge0$ and $(4k-2)Z'\ge(1+2\lambda Z(0))/k>1/k$
by the derivative bound in the proof of Lemma~\ref{lem:line-convex}.

Let $r_w$ be the unique zero of $w$ from Lemma~\ref{lem:w-one-node}.
For $0<r<r_w$ we have $w>0$ and $G(r)-G(r_w)<0$. For $r>r_w$
we have $w<0$ and $G(r)-G(r_w)>0$. Since $\int uw=0$,
\[
 \int uwG=\int uw\,[G-G(r_w)]<0 .
\]
All these integrals converge. At fixed positive frequency the profile
estimate \cite[(2.4)]{FHnew} gives $ru\asymp e^{-\sqrt\lambda r}$.
The radial equation for $ru$ gives $(ru)'=O(e^{-\sqrt\lambda r})$
by integration, so $|u'|/u$ is bounded for large $r$.
Together with \eqref{eq:radial}, this gives $uz\in H^2(\R^3)$ and
$uG=O(u^{1-1/k})$ at infinity. Multiply \eqref{eq:global-power-test}
by $uw$, integrate, and use $Lw=-u$ and \eqref{eq:global-first-moment} to get
\[
 -\int u^2z=\int w L(uz)
   =\int uwG-\frac{4k+2}{k(k+1)}(k-3/2)M .
\]
Since $\int uwG<0$, this identity and \eqref{eq:global-upper-moment} give
\[
 \frac{(2k+1)(2k-3)}{k(k+1)}
 <\frac{\int u^2z}{M}
 <\frac{(2k-1)^2}{(k+1)(4k-5)} .
\]
All denominators are positive. Multiplying out gives
\[
 (2k+1)(2k-3)(4k-5)-k(2k-1)^2
 =12k^3-32k^2+7k+15<0.
\]
By \eqref{eq:explicit-slope} this would give $M'(0+)<0$, whereas
Proposition~\ref{prop:threshold} gives $M'(0+)\ge0$ for
$p\le p_\dagger$. Thus $M'$ has no positive zero, including when
$p=p_\dagger$. Proposition~\ref{prop:explicit-upper-frequency} gives
$M'>0$ at large frequencies. Continuity gives the same sign for all
$\lambda>0$.
\end{proof}

\begin{proof}[\normalfont Proof of Theorem~\ref{thm:global-birth}]
On $q=2p-1$, Lemma~\ref{lem:birth-line} gives strict increase for
$1<p\le p_\dagger$. For $p_\dagger<p<5/3$,
Proposition~\ref{prop:threshold} gives $M'(0+)<0$, so
$M(\lambda)<M(0)$ for all sufficiently small $\lambda>0$.
Since $q<7/3$, $M(\lambda)\to\infty$ as $\lambda\to\infty$
by Proposition~\ref{prop:reverse-ends}. Continuity on $[0,\infty)$ now
gives a global minimum below $M(0)$ at a positive frequency. This
proves the assertions on the line.

We turn to exponent pairs near $(p_\dagger,2p_\dagger-1)$.

We next prove that the slope stays positive beyond the small-frequency interval.
Fix the common $\delta>0$ from Proposition~\ref{prop:birth-local}.
The frequency bound $\Lambda_{p,q}$ in
\eqref{eq:explicit-upper-frequency} is continuous for $1<p<q<7/3$.
Choose a closed rectangle about $(p_\dagger,2p_\dagger-1)$ inside the
neighborhood given by Proposition~\ref{prop:birth-local}. Choose $A>\delta$ larger than
the maximum of $\Lambda_{p,q}$ on the rectangle. Then
Proposition~\ref{prop:explicit-upper-frequency} gives
\[
 M'_{p,q}(\lambda)>0 \qquad(\lambda\ge A)
\]
for every exponent pair in this rectangle.

On $[\delta,A]$, Lemma~\ref{lem:birth-line} gives a positive minimum
at the critical exponent pair,
\[
 c_0=\min_{\delta\le\lambda\le A}
 M'_{p_\dagger,\,2p_\dagger-1}(\lambda)>0 .
\]
By Lemma~\ref{lem:rad}, the function
\[
 (p,q,\lambda)\longmapsto M'_{p,q}(\lambda)
 =2\int u_\lambda\,\partial_\lambda u_\lambda
\]
is continuous for $\lambda>0$. Uniform continuity on the closed
exponent rectangle times $[\delta,A]$ gives, after shrinking the rectangle,
\[
 \sup_{\delta\le\lambda\le A}
 \big|M'_{p,q}(\lambda)-M'_{p_\dagger,\,2p_\dagger-1}(\lambda)\big|
 < c_0/2 .
\]
Thus $M'_{p,q}>0$ on $[\delta,A]$ as well. Any positive stationary
frequency for an exponent pair in the smaller rectangle must lie in
$(0,\delta)$.

It remains to count the frequencies at each mass. The endpoint values are
\[
 M_{p,q}(0)=a_0,\qquad \lim_{\lambda\to\infty}M_{p,q}(\lambda)=\infty
\]
by Proposition~\ref{prop:reverse-ends}, since $q<7/3$.
If $p\le P(h)$, Proposition~\ref{prop:birth-local} gives $M'_{p,q}>0$
also on $(0,\delta]$, even at $p=P(h)$ where $M'_{p,q}(0+)=0$.
Hence the mass increases from $a_0$ to infinity, taking each value
in $[a_0,\infty)$ exactly once.

If $p>P(h)$, Proposition~\ref{prop:birth-local} gives
$M'_{p,q}<0$ on $(0,\lambda_{\min})$,
$M''_{p,q}(\lambda_{\min})>0$, and $M'_{p,q}>0$ on
$(\lambda_{\min},\delta]$. Since the derivative stays positive for
$\lambda\ge\delta$, this is the unique global minimum.
The decreasing branch $[0,\lambda_{\min}]$ takes each mass in
$[m,a_0]$ once. The increasing branch
$[\lambda_{\min},\infty)$ takes each mass in $[m,\infty)$ once.
They share the single frequency $\lambda_{\min}$ at $a=m$.
For $a=a_0$ the two frequencies are $0$ and $\lambda_+$.
Fixed-frequency uniqueness now gives the counts in the theorem.
The convergence $\lambda_{\min}\to0$ follows from
Proposition~\ref{prop:birth-local}.
\end{proof}

\section*{Declarations}
\noindent\textbf{Conflict of interest.} The author declares no conflict of interest.

\medskip\noindent\textbf{Data availability.} The figures are drawn from the formulas and qualitative results in the text.

\medskip\noindent\textbf{Use of AI-assisted tools.}
The author used GPT-6 Astra to help check proof steps and edit the English.
The author developed and verified the mathematical ideas independently,
reviewed all AI outputs, and corrected errors in the wording and arguments.
The author takes full responsibility for the accuracy and rigor of the paper.

\end{document}